\documentclass[11pt,reqno]{amsart}
\usepackage{graphicx,verbatim,lineno,titletoc}
\usepackage{amssymb,mathrsfs}
\usepackage{amsmath}
\usepackage{adjustbox}
\usepackage{calc}
\usepackage{ytableau}
\usepackage{array}
\usepackage{tikz}
\usetikzlibrary{shapes.multipart,patterns,arrows}
\usepackage[font=footnotesize]{caption}
\usepackage[T1]{fontenc} 
\usepackage{fourier} 
\usepackage[english]{babel} 
\usepackage{appendix}
\usepackage[all]{xy}
\usepackage{enumerate}
\usepackage[english]{babel}
\usepackage{multirow}
\usepackage{float}
\usepackage{enumitem}
\usepackage{accents}
\usepackage[numbers]{natbib}
\usepackage{hyperref}
\hypersetup{colorlinks}
\usepackage{algorithm}
\usepackage[noend]{algpseudocode}
\usepackage{mathtools}
\mathtoolsset{showonlyrefs}
\usepackage{wrapfig}

\usepackage[top=1.0in, left=1.3in, right=1.3in, bottom=1.0in]{geometry}
\hypersetup{colorlinks,linkcolor={red},citecolor={olive},urlcolor={red}}

\newcolumntype{M}[1]{>{\centering\arraybackslash}m{#1}}
\newcolumntype{N}{@{}m{0pt}@{}}

\newtheorem{theorem}{Theorem}[section]
\newtheorem{prop}[theorem]{Proposition}

\newtheorem{lemma}[theorem]{Lemma}
\newtheorem{corollary}[theorem]{Corollary}
\newtheorem{conjecture}[theorem]{Conjecture}

\def\Z{\mathbb{Z}}

\def\eps{\varepsilon}

\def\depth{\mathtt{dep}}

\def\root{\mathtt{r}}
\def\depth{\mathtt{dep}}
\def\mod{\textup{mod\, }}

\newcommand{\floor}[1]{\left\lfloor#1\right\rfloor}

\definecolor{hancolor}{rgb}{0.1, 0.0, 0.9}

\def\Z{\mathbb{Z}}

\def\eps{\varepsilon}

\def\1{\boldsymbol{1}}
\def\2{\boldsymbol{2}}

\newlength\myindent
\newtheorem{innercustomthm}{Theorem}
\newenvironment{customthm}[1]
{\renewcommand\theinnercustomthm{#1}\innercustomthm}
{\endinnercustomthm}

\newtheorem{innercustomdef}{Definition}
\newenvironment{customdef}[1]
{\renewcommand\theinnercustomdef{#1}\innercustomdef}
{\endinnercustomdef}

\theoremstyle{definition}

\allowdisplaybreaks

\makeatletter
\newcommand{\addresseshere}{%
	\enddoc@text\let\enddoc@text\relax
}
\makeatother

\usepackage{color}

\def\boxit#1{\vbox{\hrule\hbox{\vrule\kern6pt\vbox{\kern6pt#1\kern6pt}\kern6pt\vrule}\hrule}}

\begin{document}

	\title[Discrete pulse-coupled oscillators on trees]{ Time complexity of Synchronization \\ of discrete pulse-coupled oscillators on trees}

	\author{Hanbaek Lyu}
	\address{Hanbaek Lyu, Department of Mathematics, University of Wisconsin - Madison, WI 53709, USA}
	\email{\texttt{hlyu@math.wisc.edu}}

\begin{abstract}
	A major open question in the study of synchronization of coupled oscillators is to find necessary and sufficient condition for a system to synchronize on a given family of graphs. This is a difficult question that requires to understand exactly how the nonlienar interaction between local entities evolves over the underlying graph. Another open question is to obtain bounds on the time complexity of synchronization, which has important practical implications in clock synchronization algorithms. We address these questions for one-parameter family of discrete pulse-coupled inhibitory oscillatorscalled the $\kappa$-color firefly cellular automata (FCA). 
	Namely, we show that for $\kappa\le 6$, recurrence of each oscillator is a necessary and sufficient condition for synchronization on finite trees, while for $\kappa \ge 7$ this condition is only necessary. As a corollary, we show that any non-synchronizing dynamics for $\kappa\le 6$ on trees  decompose into synchronized subtrees partitioned by `dead' oscillators. Furthermore, on trees with diameter $d$ and maximum degree at most $\kappa$, we show that  the worst-case number of iterations until synchronization is of order $O(\kappa d)$ for $\kappa\in \{3,4,5\}$,  $O(\kappa d^{2})$ for $\kappa=6$, and infinity for $\kappa \ge 7$.  
	Lastly, we report simulation results of FCA on lattices and conjecture that on a finite square lattice, arbitrary initial configuration is synchronized under $\kappa$-color FCA if and only if $\kappa=4$. 
\end{abstract}



\maketitle

\section{Introduction}
Many physical and biological complex systems consist of levels of hierarchies of locally interacting dynamic units, whose internal dynamics are induced by non-linear aggregation of local interactions between units at lower levels. Top levels are forced to have a certain macro-behavior suitable for survival, which is miraculously supplied by the right micro-level local interactions, forged by the evolutionary process. This chain of emergent dynamics is at the heart of the challenge we are facing in understanding not only biological systems but also many other complex systems in our society as well as in designing cooperative control protocol of large networked systems \cite{mesbahi2010graph, strogatz2001exploring}. 

Consisting of only two levels of hierarchies with simple internal dynamics for units at the bottom level, a system of coupled oscillators has been a central subject in non-linear dynamical systems literature for decades \cite{strogatz2000kuramoto}. As populations of blinking fireflies \cite{buck1938synchronous} and circadian pacemaker cells \cite{enright1980temporal} do, two neighboring oscillators are coupled so that they tend to synchronize their phase or frequency, and the question is whether the such local tendency to synchrony does lead to global synchronization in the entire network. Despite their simplicity, they exhibit many fundamental difficulties which repel our traditional reductionist approach based on linear methods, and yet our enhanced knowledge of such systems is finding fruitful applications, ranging from robotic vehicle networks \cite{nair2007stable} to electric power networks \cite{dorfler2012synchronization}, and more recently, to distributed control of wireless sensor networks \cite{hong2005scalable, pagliari2011scalable, wang2012energy, wang2013increasing}.

A major open question in the study of synchronization of coupled oscillators is to find necessary and sufficient condition for a system to synchronize on a given family of graphs. This is a difficult question that requires to understand exactly how the nonlienar interaction between local entities evolves over the underlying graph. Another open question is to obtain bounds on the time complexity of synchronization, which has important practical implications in clock synchronization algorithms. 
These questions are exteremely difficult to address for traditional smooth oscillators such as Kuramoto oscillators \cite{kuramoto1975self}, where the time complexity should be understood as the time required to reach an `$\eps$-synchronization'. In this work, we address these questions for a simple model of discrete coupled oscillators. 

Discrete coupled oscillators provide a simple framework for coupled oscillators and clock synchronization problems: A vertex coloring $X_{t}:V\rightarrow \mathbb{Z}_{\kappa}=\{0,1,\cdots, \kappa-1\}$ on a given graph $G=(V,E)$ is updated in discrete time according to a fixed local transition rule. To model cyclic phase update of oscillators, it is natural to assume $X_{t+1}(v) = X_{t}(v)$ or $X_{t}(v)+1\, (\text{mod}\,\kappa)$ depending on local data. One of the main advantages of this discrete setting is that the nonlinear dynamics of the oscillators over spatial networks become more tractable and it can be a testbed for developing general proof techniques. An immediate advantage in considering discrete models for coupled oscillators is that one can easily compute examples of synchronizing as well as non-synchronizing examples on various graphs. Moreover, one can simply observe how and why a particular trajectory synchronizes or converges to a periodic orbit (see Figures \ref{fig:ex1} and \ref{fig:ex2}). Another advantage in a such discrete setting is that by drawing the initial coloring $X_{0}$ from some probability measure and asking how the probability $\mathbb{P}(\text{$X_{t}$ has property $Q$})$ behaves in $t$, one can adopt various probabilistic techniques to study the synchronization problem \cite{lyu2016synchronization, gravner2018limiting, lyu2018persistence}.


For synchronous systems with discrete oscillators taking $\kappa$ distinct phase values, a number of algorithms that are self-stabilizing on trees with constant memory per node are known: e.g., for $\kappa=3$ by Herman and Ghosh \cite{herman1995stabilizing}, for all odd $\kappa \ge 3$ by Boulinier et al. \cite{boulinier2006toward}. Upper bounds of $O(\texttt{diameter})$ for time complexity of synchronization are known for such algorithms. The problem of designing a discrete model for coupled oscillators that have the capacity to synchronize arbitrary $\kappa$-coloring on a class of finite graphs has been known as the `digital clock synchronization problem' in distributed algorithms literature. If one allows $\kappa$ to grow with the size of $G$, then there is such a solution that works on arbitrary finite graphs (e.g., see Dolev \cite{dolev2000self} or Arora et al. \cite{arora1992maintaining}). Roughly speaking, the idea is that if $\kappa$ is large enough, then one can let every vertex adapt the local maximum color within distance 1 at each time step in parallel; then the globally maximum color would propagate and ``eat up'' all vertices. In fact, this idea of ``tuning toward maximum'' dates back to a famous consensus algorithm by Lamport \cite{lamport1978time}. One can readily see that such an algorithm relies on some notion of global total ordering among colors of vertices, which is not the case in our case due to the cyclic nature of the color space $\Z/\kappa\mathbb{Z}$. In fact, this issue arising from the cyclic hierarchy between colors is fundamental to our problem, and in fact, is a key source that generates interesting emergent behavior in the system. Hence we may restrict ourselves to GCA models with $\kappa$ independent of $G$. 

Dolev \cite{dolev2000self} showed that no GCA model using a fixed $\kappa$ is able to synchronize arbitrary $\kappa$-coloring on all connected finite graphs. Roughly speaking, for any such given GCA model, one can construct a symmetric configuration on a cycle of some length so that the vertices have no way to break the symmetry by blindingly following a homogeneous local rule. On trees, however, such a construction is topologically prohibited so one may hope that there exists a $\kappa$-color GCA model which synchronizes all initial $\kappa$-colorings on any finite trees. Indeed, a 3-color GCA model was studied by Herman and Ghosh \cite{herman1995stabilizing}, and odd $\kappa\ge 3$ models by Boulinier, Petit, and Villain \cite{boulinier2006toward}. When $\kappa=3$, the latter model coincides with another well-known GCA model called the cyclic cellular automaton, which was introduced by Bramson and Griffeath \cite{bramson1989flux} as a discrete-time analog of the cyclic particle systems. In Gravner, Lyu, and Sivakoff \cite{gravner2016limiting}, the limiting behavior of 3-color cyclic cellular automaton together with the 3-color Greenberg-Hastings model \cite{greenberg1978spatial} on infinite trees were investigated using probabilistic methods.


\subsection{The Firefly Cellular Automata}

The $\kappa$-color firefly cellular automata (FCAs) is a discrete model for pulse-coupled inhibitory oscillators, which was first proposed in  \cite{lyu2015synchronization}. The model is defined for each integer $\kappa\ge 3$, which equals the number of distinct states that each oscillator assumes. A special state $b(\kappa) = \floor{\frac{\kappa-1}{2}}$ is designated as the `blinking' color. In a network of $\kappa$-state identical oscillators, each oscillator increments from state $i$ to $i+1(\text{mod $\kappa$})$ unless it has a neighbor of blinking color and its color is post-blinking (specifically, $i \in \{ b(\kappa)+1, b(\kappa)+2,\dots, \kappa-1 \}$), in which case it waits for one iteration without an update. More precisely, the transition map $\tau:X_{t}\mapsto X_{t+1}$ for the $\kappa$-color FCA is given as follows:
\begin{equation*}
	\text{(FCA)} \qquad X_{t+1}(v)=\begin{cases}
		X_t(v) & \text{if $X_t(v) > b(\kappa)$ and $|\{u\in N(v)\,:\, X_t(u)=b(\kappa)\}| \ge 1 $}\\
		X_t(v)+1 (\,\text{mod} \kappa)& \text{otherwise}\end{cases}
\end{equation*}
where $N(v)$ denotes the set of all neighbors of $v$ in $G$. We say a vertex $v$ \textit{blinks} at time $t$ if $X_{t}(v)=b(\kappa)$, is \textit{pulled} at time $t$ if $X_{t+1}(v)=X_{t}(v)$, and \textit{pulls} its neighbor $u$ at time $t$ if $u$ is pulled at time $t$ and $v$ blinks at time $t$. Given a $\kappa$-color FCA trajectory $(X_{t})_{t\ge 0}$ on a graph $G=(V,E)$, we say $X_{t}$ (or $X_{0}$) \textit{synchronizes} if there exists $N\ge 0$ such that $X_{t}\equiv Const.$ for all $t>N$.

\begin{figure*}[h]
	\centering
	\includegraphics[width=0.96 \linewidth]{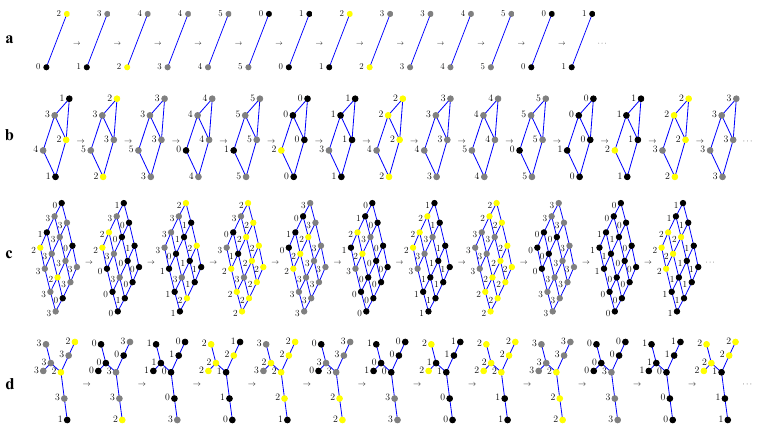}
	\caption{Four examples of synchronizing 6-color FCA trajectories. Dots and blue lines represent nodes (oscillators) edges in the underlying graph, respectively.  Nodes with pre-blinking color (in $\{0,1\}$), blinking color $b(6)=2$, and post-blinking color (in $\{ 3,4,5 \}$) are marked with colors black, yellow, and gray, respectively.  For \textbf{c} and \textbf{d}, the last configurations shown consist of two colors so they are guaranteed to synchronize by \cite[Cor. 3.3]{lyu2015synchronization}.}
	\label{fig:ex1}
\end{figure*}

\begin{figure*}[h]
	\centering
	\includegraphics[width=0.96 \linewidth]{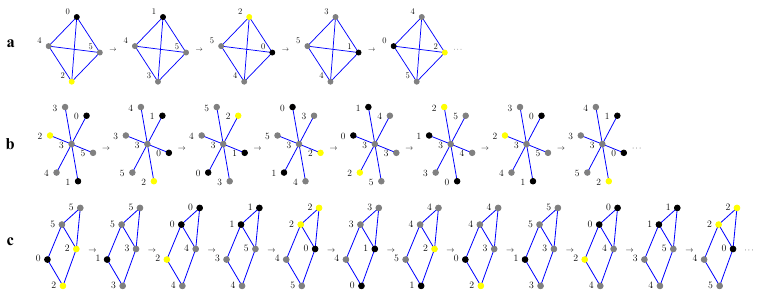}
	\caption{Three examples of non-synchronizing 6-color FCA trajectories.  Dots and blue lines represent nodes (oscillators) edges in the underlying graph, respectively. Nodes with pre-blinking color (in $\{0,1\}$), blinking color $b(6)=2$, and post-blinking color (in $\{ 3,4,5 \}$) are marked with colors black, yellow, and gray, respectively. In \textbf{a}, the last configuration is symmetric to the initial configuration. }
	\label{fig:ex2}
\end{figure*}

Being a deterministic dynamical system with finite state space for each vertex, any $\kappa$-color FCA trajectory $(X_{t})_{t\ge 0}$ on any finite graph $G=(V,E)$ must converge to a periodic limit cycle.  Limit cycles can be an either synchronous or asynchronous periodic orbits, as illustrated in the examples of 6-color FCA trajectories in Figures \ref{fig:ex1} and \ref{fig:ex2}. Note that $b(6)=2$ is the blinking state in this case, so every vertex of state $3,4,$ or $5$ with a state 2 neighbor stops evolving for one iteration and all the other vertices evolve to the next state. For instance, in the examples in Figure \ref{fig:ex1}, the `post-blinking' colors are shown in gray, and whenever they are adjacent to a vertex of blinking color (yellow), they maintain the same color for one iteration. While all examples in Figure \ref{fig:ex1} synchronize, the three examples in Figure \ref{fig:ex2} do not. 

\subsection{FCA on 2D lattices : stable periodic objects and the 4-color criticality}

The one-dimensional integer lattice with nearest neighbor coupling is an extremely restricted environment, where local dynamics on any finite interval may be easily perturbed by external excitation. However, on higher dimensional lattices $\mathbb{Z}^{d}$ for $d\ge 2$, each site is likely to be in contact with many different colors, which enables the emergence of \textit{stable periodic objects (SPOs)}. This drives the characteristic local periodicity by spontaneous spiral formations.  

SPOs are defined to be configurations on a finite region $\Omega_{0}\subset \mathbb{Z}^{2}$ whose dynamics are independent of their surroundings. This requires that all sites on the inside boundary $\partial \Omega_{0}$ must excite with maximal rate since otherwise, they could get further excitation from outside. For example, in CCA, a \textit{clock} is a configuration on a directed cycle $\vec{C}$ on which the color increases by  $1\mod \kappa$ as one traverses around it. Such a configuration is clearly invariant under CCA dynamics since all sites must advance in each step. Hence all sites on $\vec{C}$ excite with maximal rate. Now any simply connected region $\Omega_{0}$ whose boundary has a clock configuration is an SPO, as clocks separate the dynamics of their interior from exterior. Moreover, arbitrary SPO in CCA is given by a non-disjoint union of a finite number of clocks where the outer boundary forms the boundary of SPO. Furthermore, a moment's thought reveals that the same combinatorial architecture also holds for SPOs in GHM dynamics. Note that in both dynamics, it is a simple combinatorial fact that SPOs grow in space, driving their surroundings into the same periodicity. Hence any CCA or GHM dynamics with at least one SPO is locally periodic, where multiple SPOs grow and partition the entire space into near-Voronoi cells. 

Now a standard argument for local periodicity proceeds as follows. Consider a random initial configuration $\eta_{0}$ drawn from uniform product measure on $\mathbb{Z}_{\kappa}^{\mathbb{Z}^{d}}$, and suppose that some finite region $\Omega_{0}\subset \mathbb{Z}^{d}$ has an SPO. Then with high probability, there exists an SPO on a translation of $\Omega_{0}$, so $\eta_{t}$ is locally periodic a.s. Indeed, for arbitrary $\kappa\ge 3$, $\mathbb{Z}^{2}$ contains a cycle of length some multiple of $\kappa$ so clocks, and hence SPOs, do exist with positive probability. This determines limiting behaviors of CCA and GHM on $\mathbb{Z}^{2}$.

\begin{figure*}[h]
	\centering
	\includegraphics[width=1 \linewidth]{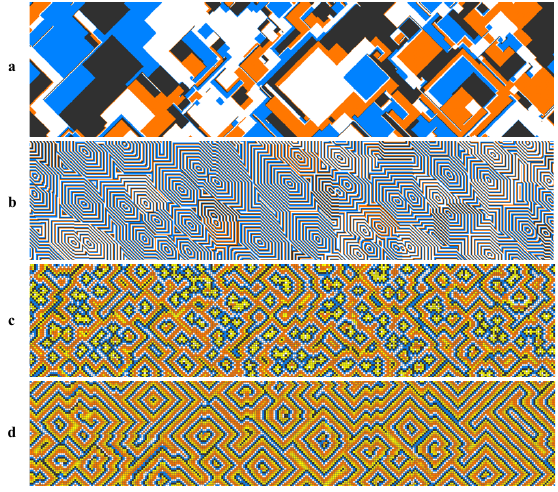}
	\caption{ Snapshots of FCA on two-dimensional lattices at iteration $770$ with random initial configurations (drawn from the uniform distribution). \textbf{a} 4-color FCA on the square lattice $\Z^{2}$; \textbf{b} 4-color FCA on a hexagonal lattice obtained by adding edges between $(i,j)$ and $(i+1,j-1)$ for all $i,j\in \mathbb{Z}$; \textbf{c} 5-color FCA on the square lattice $\Z^{2}$; and \textbf{d} 6-color FCA on the square lattice $\Z^{2}$. We observe global synchronization for case \textbf{a} and stable wave generators in all other cases. Experiments show that 4-color FCA on higher-dimensional square lattices $\Z^{d}$ for $d\ge 3$ synchronizes arbitrary initial configuration. 
	}
	\label{fig:FCA4_Z2}
\end{figure*}

However, SPOs do not completely determine the global dynamics of discrete excitable media in higher dimensions in at least two ways. Firstly, computer simulation suggests SPOs "pop up" out of a chaotic region rather than residing in the very initial configuration, so one needs to understand the mechanism of this spontaneous emergence of SPOs. A nice illustration of the global dynamics of CCA on $\mathbb{Z}^{2}$ in four distinct stages is given in \cite{fisch1991cyclic}.

Secondly, and more crucially, SPOs may not exist after all. Indeed, it can be easily shown that in FCA with arbitrary $\kappa$, there exist no SPOs on $\mathbb{Z}^{2}$. To see this, note that sites with maximal excitation in FCA are sites that never blink and this requires them to be adjacent to at least one blinking neighbor in each instant. In fact, since all sites blink at most once in every $\kappa$ times, such non-blinking sites should be adjacent to at least $\kappa\ge 3$ blinking sites. Now suppose a finite region $\Omega_{0}\subset \mathbb{Z}^{2}$ has an SPO so that whose inside boundary consists of ``non-blinking" sites. Observe that any site on the boundary is adjacent to at most two sites in the interior of $\Omega_{0}$. Hence boundary sites must be excited from outside of $\Omega_{0}$, a contradiction.

Nevertheless, simulation suggests that for $\kappa\ne 4$, FCA on $\mathbb{Z}^{2}$ and $\kappa=4$ on the hexagonal grid should have similar locally periodic limiting behavior driven by spontaneous spiral formations (see figure \ref{fig:FCA4_Z2}). The ``centers" of these spirals \textit{can} be affected by external perturbations, but it seems that they are dynamically stable in the sense that spiral centers in a random environment eventually settle down and organizes the surroundings.

On the other hand, the 4-color FCA shows the most distinctive behavior on $\mathbb{Z}^{d}$ on square lattices. It shows clustering via spontaneously generated cancelative rectilinear waves. This distinguishes the 4-color FCA on square grids from most 2-dimensional excitable media, where spontaneous spiral formation eventually governs the dynamics and hence clustering does not occur (see for instance \cite{fisch1991cyclic} for CCAs and GHMs, \cite{mikhailov2006control} for BZ reaction, and \cite{barkley1990spiral} for continuous excitable media). Even more strikingly, the simulation suggests that this behavior does not depend on the dimension $d$. 

We should mention that the mechanism of clustering in higher dimensions should be significantly different from that in one dimension, which is essentially annihilating particle systems behavior. Here we give a sketch of the dynamics. When more than two different colors meet around a site, they tend to form a spiral but the center then moves with the same speed as the winding arms, resulting in a semi-spiral traveling across the space. These are called \textit{singularities} for a reason that will be explained in the later sections. Hence the dynamic is governed by traveling singularities and \textit{cyclic boundaries} between regions of two consecutive colors $i$ and $i+1$ $\mod 4$, which moves toward the color $i$. Singularities coalesce or annihilate each other upon collision, and regions without singularities are eventually taken over by one color and keep synchronized, at least until it is invaded by another singularity or cyclic boundary.            

We propose the following conjecture about 4-color FCA on $\mathbb{Z}^{d}$:

\begin{conjecture}
	Let $\kappa\ge 3$ and $d\ge 2$, and let $\mathbb{P}$ be a uniform product probability measure on $\kappa$-colorings on $\mathbb{Z}^{d}$. Let $X_{0}$ be a random $\kappa$-coloring drawn from $\mathbb{P}$. Then 
	\begin{description}[noitemsep]
		\item[(i)] If $\kappa=4$, then $\mathbb{P}$-a.s. $X_{t}$ clusters, i.e., for any finite region $\Omega_{0}\subset \Z^{2}$, we have 
		\begin{equation*}
			\lim_{t\rightarrow \infty} \mathbb{P}(\text{$X_{t}\equiv Const.$ on $\Omega_{0}$}) = 1
		\end{equation*}
		
		\item[(ii)] If $\kappa\ne 4$, then $\mathbb{P}$-a.s. $X_{t}$ is uniformly locally periodic with period $\kappa+1$, i.e., 
		
		for each site $x\in \mathbb{Z}^{2}$, 
		\begin{equation*}
			\lim_{t\rightarrow \infty} \mathbb{P}( X_{t}(x)=X_{t+\kappa+1}(x) )= 1
		\end{equation*}
	\end{description}
\end{conjecture}

\subsection{Statement of main results}

In \cite{lyu2015synchronization}, it was shown that for any $\kappa\ge 3$ and any finite path $P$ of $d$edges, arbitrary $\kappa$-color FCA trajectory on $P$ converges in $d\kappa(\kappa+4)/2$ iterations. However, as we can see from the non-synchronizing example in Figure \ref{fig:ex2}\textbf{b}, not every initial configuration synchronizes on finite trees. The reason that the example in Figure \ref{fig:ex2}\textbf{b} does not synchronize is that the center of the star could be pulled by the leaves (degree-1 nodes) constantly. One can easily generalize this observation to construct a non-synchronizing example on trees. Namely, let $v$ be a vertex in finite tree $T$ with degree $\ge \kappa$, and let $T_{1},\cdots,T_{m}$ be the connected components of $T-v$, the graph obtained from $T$ by deleting $v$ together with edges incident to it. Note that $m\ge \kappa$. Assign state $i(\text{mod $\kappa$} )$ to every vertex of $T_{i}$, and assign any state $>\kappa/2$ to vertex $v$. Then $v$ never blinks and each component $T_{i}$ never gets pulled by $v$, which is essentially the counterexample in Figure \ref{fig:ex2}\textbf{b}.

The non-synchronizing example of $\kappa$-color FCA on trees with maximum degree at least $\kappa$ shows that, if every $\kappa$-coloring on $T$ synchronizes, then necessarily $T$ has maximum degree $<\kappa$. A natural question is then the converse of this statement: If the maximum degree of a tree $T$ is less than $\kappa$, does every initial configuration of $\kappa$-color FCA on $T$ synchronize? In \cite{lyu2015synchronization}, it was proved that such a necessary local condition to synchronize arbitrary $\kappa$-coloring on a tree is also sufficient for $\kappa\in \{3,4,5\}$. Furthermore, when $T$ is a path, the worst-case time until synchronization is shown to be of order $O(d)$, where $d$ is the shortest-path diameter of $T$. For its statement, define 
\begin{align}
	N_{G}(\kappa):=  \max_{X_{0}: V\rightarrow \Z/\kappa\Z } \inf\{ n\ge 1\,|\, \textup{$X_{n}$ is synchronized}  \}. 
\end{align}
We take the convention of $\inf \emptyset =\infty$. Hence $N_{G}(\kappa)<\infty$  if and only if every possible $\kappa$-coloring on $G$ synchronizes. 

The main result in the present work is a complete characterization of the long-term behavior of FCA on finite trees for all $\kappa\ge 3$ as well as giving an explicit bound on the worst-case time-to-synchronization. Roughly speaking we prove that FCA on a tree with diameter $d$ synchronizes arbitrary initial configuration within $O(d)$ iterations if the maximum degree of the tree is $<\kappa$, and there exists some non-synchronizing configuration otherwise. 
\begin{theorem}[Time complexity of synchronization]\label{thm:tree_time}
	Suppose $\kappa\in \{3,4,5,6\}$. Let $T=(V,E)$ be a finite tree with a maximum degree less than $\kappa$ and shortest-path diameter $d$. Then the following hold:
	\begin{description}
		\item[(i)]  For $\kappa\in \{3,4,5\}$, $N_{T}(\kappa) \le C_{\kappa} (d+1)$, where $C_{3}=9$, $C_{4}=30$, and $C_{5}=24$.
		\item[(ii)]  For $\kappa=6$, $N_{T}(\kappa)\le C_{6} d^{2}$ for some constant $C_{6}>0$. 
		\item[(iii)]  For $\kappa\ge 7$, $N_{T}(\kappa)=\infty$.
	\end{description}
\end{theorem}

\noindent An immediate consequence of Theorem \ref{thm:tree_time} is that there is a critical number of colors ``between'' 6 and 7; with fewer colors, maximum degree $<\kappa$ implies synchronization of arbitrary $\kappa$-coloring, and with more colors, there are non-synchronizing examples on trees with maximum degree $\le \kappa/2+1$. This is analogous to the clustering-fixation phase transition of the $\kappa$-color cyclic cellular automaton on $\mathbb{Z}$ (see Fisch \cite{fisch1990one}).

\begin{corollary}[Phase transition of FCA on trees]\label{cor:phase_transition}
	${}$
	\begin{description}
		\item[(i)] If $\kappa\in \{3,4,5,6\}$ and $T=(V,E)$ is any finite tree, then every $\kappa$-coloring on $T$ synchronizes if and only if  $T$ has maximum degree $<\kappa$. 
		\item[(ii)] If $\kappa\ge 7$, then there exists a finite tree $T=(V,E)$ with maximum degree $\le \kappa/2 + 1$ and a non-synchronizing $\kappa$-coloring on $T$.  
	\end{description}
\end{corollary}

Next, we turn our attention to characterizing the long-term dynamics of the 6-color FCA on trees without any degree restriction. In Figure \ref{fig:ex2}, we have seen a non-synchronizing example for the 6-color FCA on a star with $6$ leaves. There, the center node is pulled by one of its 6 neighbors every instance so its phase is `stuck' at a post-blinking color (3 in that example). A natural question to ask then is the following: If we have a non-synchronizing 6-color FCA on a tree, should there be some node that is stuck at some post-blinking color? We answer this question positively in Theorem \ref{thm:blinkingtree} below. 

\begin{theorem}[Local recurrence implies global synchronization] \label{thm:blinkingtree}
	Let $T$ be a finite tree and let $X_{0}$ be a $6$-coloring on $T$. Then the $6$-color FCA trajectory $(X_{t})_{t\ge 0}$ synchronizes if and only if every vertex of $T$ blinks infinitely often in the trajectory $(X_{t})_{t\ge 0}$. 
\end{theorem}

We remark that the analogous result for $\kappa\in \{3,4,5\}$ was shown in \cite{lyu2015synchronization}. For $\kappa\ge 7$, the above statement is false due to Corollary \ref{thm:tree_time} \textbf{(iii)}. However, the proof for $\kappa=6$ is substantially more difficult. Roughly speaking, for $\kappa\in\{3,4,5\}$, assuming that every vertex blinking infinitely often, analyzing local dynamics within a subtree of diameter $\le 3$ easily leads to a `constant-time reduction' of the dynamics on a smaller subtree. However, for $\kappa=6$ case, one needs to analyze the dynamics of the entire tree to see any reduction or contradiction. In order to overcome this difficulty,  we introduce a novel technique of  `fractal dynamics', which recursively combine dynamics on subtrees to form dyanmics on larger trees that maintains the same structure (see Lemma \ref{lem:recursion_fractalbranch}). 

An immediate consequence of Theorem \ref{thm:blinkingtree} (as well as the corresponding result for $\kappa\in \{3,4,5\}$ in \cite{lyu2015synchronization}) is a structure theorem for non-synchronizing FCA on trees. Suppose $(X_{t})_{t\ge 0}$ is a non-sychronizing $6$-color FCA on a finite tree $T$. By Theorem \ref{thm:blinkingtree}, the set $V_{0}$ of nodes that blink only finitely many times is nonempty. All nodes in $V_{0}$ fixate to some post-blinking color so the nodes in different connected components in $T-V_{0}$ (deleting all nodes in $V_{0}$ from $T$) do not interact after some finite time. Since all nodes not in $V_{0}$ blink infinitely often, by Theorem \ref{thm:blinkingtree}, the dynamics restricted on each connected components of $T-V_{0}$ will eventually synchronize. This gives the following corollary. 

\begin{corollary}[Characterization of FCA dynamics on  trees]\label{mainthm_tree_ch}
	Suppose $\kappa\in \{3,4,5,6\}$. Let $T=(V,E)$ be a finite tree and let $(X_{t})_{t\ge 0}$ denote an arbitrary  $\kappa$-color FCA trajectory on $T$. Then there exists a proper subset $V_{0}\subset V$ (possibly empty) and an intger $N\ge 1$ such that the following hold: For all $t\ge N$, 
	\begin{description}
		\item[(i)] $X_{t}$ restricted on each connected component of $T-V_{0}$ is synchronized (constant); and 
		\item[(ii)]  For all $v\in V_{0}$, $X_{t}(v)\equiv c > b(\kappa)$ and $v$ is adjacent to all $\kappa$ colors at time $t$.
	\end{description}
\end{corollary}



\subsection{Organization}

This paper is organized as follows. We establish some preliminary lemmas on FCA in Section \ref{sec:preliminary}. We prove Theorem \ref{thm:tree_time} for $\kappa=4,5$ in Section \ref{sec:pf_tree_time_45}. In the following section, Section \ref{sec:pf_FCA_7}, we prove Theorem \ref{thm:tree_time} for $\kappa\ge 7$. In Section \ref{sec:pf_blinkingtree_6}, we provide proofs of Theorem \ref{thm:tree_time} for $\kappa=6$ and Theorem \ref{thm:blinkingtree} assuming two technical lemmas (see Lemmas \ref{lem:terminal_is_fractal} and \ref{lem:recursion_fractalbranch}). The remaining sections are devoted to proving Lemmas \ref{lem:terminal_is_fractal} and \ref{lem:recursion_fractalbranch}. Namely, we characterize irreducible local dynamics of 6-color FCA on branches in Section \ref{sec:6FCA_local} and then prove Lemmas \ref{lem:terminal_is_fractal} and \ref{lem:recursion_fractalbranch} in Sections \ref{sec:proof_lem_terminal_fractal} and \ref{sec:proof_lemma_recursion_fractal}, respectively. 

\subsection{Notations}

Throughout this paper, we assume every graph is finite, simple, and connected unless otherwise mentioned. If $S,H$ are vertex-disjoint subgraphs of $G$, then $S+H$ is defined by the subgraph obtained from $S\cup H$ by adding all the edges in $G$ between $S$ and $H$. On the other hand, $S-H:=S-V(H)$ denotes the subgraph of $G$ obtained from $S$ by deleting all vertices of $H$ and edges incident to them. If $v\in V(H)$ and $H$ is a subgraph of $G$, then $N_{H}(v)$ denotes the set of all neighbors of $v$ in $H$ and $\deg_{H}(v):=|N_{H}(v)|$. We for each integer $\kappa\ge 1$, denote $\mathbb{Z}_{\kappa}:=\mathbb{Z}/\kappa \mathbb{Z}$.

\section{Preliminary lemmas and discussions}
\label{sec:preliminary}

\subsection{The width lemma (a.k.a. the "half-circle concentration")}

Let $u,v$ be two vertices in $G$ and fix a $\kappa$-coloring $X:V\rightarrow \mathbb{Z}_{\kappa}$ on $G$. The \textit{forward displacement of $v$ from $u$ in $X$} is defined by 
\begin{equation*}
	\delta( X(u), X(v) ):= X(v)-X(u) \,\,(\text{mod $\kappa$}) \in \{0,1,\dots, \kappa-1 \}. 
\end{equation*}
We say $v$ is \textit{ahead of} $u$ if $	\delta( X(u), X(v) )<\kappa/2$ and $v$ is \textit{lagging to} $u$ if $	\delta( X(u), X(v) )>\kappa/2$ (Note that if $\kappa$ is even, it could be that $	\delta( X(u), X(v) )=\kappa/2$). Suppose $u$ and $v$ are adjacent in $G$. The \textit{width} of a $X$ is defined to be the quantity 
\begin{align}
	w(X_{t} ) &:= \min_{v\in V} \,  \max_{u\in V} \,  \delta( X(u), X(v)) \\
	&= \min_{a\in \Z_{\kappa}} \, \min_{v\in V}  \left[ X(v)  + a \,\, (\mod \kappa)  \right].
\end{align}
Note that $w(X)$ is the length of the shortest path on the color space $\mathbb{Z}_{\kappa}$ (viewed as a cycle graph of $\kappa$ edges) that covers all states of the vertices in the configuration $X$. For instance, the first configuration in Figure \ref{fig:ex1}\textbf{b} has width 3, whereas the second to the last Figure \ref{fig:ex1}\textbf{b} has width $1$. For any subgraph $B\subseteq G$, we denote by $w_{B}(X)$ the width of the restricted configuration $X|_{V(B)}$ on $B$.

A widely recognized sufficient condition for synchronization for arbitrary $G$ is that the initial phases are concentrated in a small arc (e.g., an open half-circle of the phase space viewed as the unit circle $S^{1}=[0,1]/\Z$ in the case of Kuramoto oscillators \cite{kuramoto1975self, strogatz2000kuramoto, acebron2005kuramoto}). Such concentration of phases essentially breaks the cyclic symmetry of the oscillation cycle and induces a total ordering (e.g., from the most lagging to the most advanced), and a broad class of couplings respects such ordering and contracts the phase configuration toward synchrony. Below, we state and prove a version of such `half-circle concentration' lemma (Lemma \ref{lem:width}) for the case of FCA: It yields that any $\kappa$-color FCA trajectory on any $n$-node connected graph synchronizes in at most $\kappa^{2} n$ iterations. A similar result without a quantitative upper bound on the number of iterations until synchronization was stated in \cite[Lem. 2.2]{lyu2015synchronization} with only a sketch of a proof. 


\begin{lemma}[Width lemma]
	\label{lem:width}
	Let $G=(V,E)$ be a connected graph with shortest-path diameter $d$. Let $(X_{t})_{t\ge 0}$ denote an arbitrary $\kappa$-color FCA dynamics on $G$. Suppose $w(X_{0})<\kappa/2$. Then the following hold:
	\begin{description}
		\item[(i)] $w(X_{t})$  is non-increasing in $t$ and $w(X_{t}) <\kappa/2$ for all $t\ge 0$. 
		
		\item[(ii)] Then $X_{t}$ is synchronized for all $t\ge \frac{\kappa(\kappa+1)}{2} |V|$. 
	\end{description}
\end{lemma}

\begin{proof}
	For \textbf{\eqref{I}},  choose a node $v^{*}\in V$ such that 
	\begin{align}
		X_{0}(v)\in   \{X_{0}(v^{*}), X_{0}(v^{*})+1  ,\dots, X_{0}(v^{*}) + w(X_{0})  \} \, (\mod \kappa) \quad \textup{$\forall v\in V$}. 
	\end{align}
	That is, viewing the color space $\mathbb{Z}_{\kappa}$ as the cycle graph of $\kappa$ edges with nodes $0, 1, \dots, \kappa-1$ arranged in clockwise order, the clockwise shortest-path from $X_{0}(v^{*})$ to $X_{0}(v^{*})+w(X_{0}) \, (\mod \kappa)$ in $\mathbb{Z}_{\kappa}$ covers all colors $X_{0}(v)$ for $v\in V$. (In words, $v^{*}$ is the `most lagging' at time 0.) We claim that the following holds:
	\begin{align}\label{eq:width_lemma_claim1}
		\textup{$\forall  t\ge 1$}  : \quad 
		\begin{cases}
			w(X_{t}) \le w(X_{t-1}) \\
			X_{t}(v^{*}) = X_{0}(v^{*}) + t \, (\mod \kappa) \\ 
			X_{t}(v)\in   \{X_{t}(v^{*}), X_{t}(v^{*})+1  ,\dots, X_{t}(v^{*}) + w(X_{t})  \} \, (\mod \kappa) \quad \textup{$\forall v\in V$}.
		\end{cases}
	\end{align}
	Ineed, by the hypothesis, $w(X_{0})<\kappa/2$, so all colors in $X$ are ahead of $X_{0}(v^{*})$. This implies that $v$ is not pulled by any of its neighbors during the transition $X_{0}\mapsto X_{1}$, so $X_{1}(v^{*})=X_{0}(v^{*})+1 \, (\mod \kappa)$. This holds for any $v^{**}\in V$ such that $X_{0}(v^{**})=X_{0}(v^{*})$.  Since $X_{1}(v) \in \{ X_{0}(v), X_{0}(v)+1 \} \, (\mod \kappa)$ for all $v\in V$, it follows that \eqref{eq:width_lemma_claim1} holds for $t=1$. Repeating the same argument and using induction in $t$, this shows \eqref{eq:width_lemma_claim1}. Note that \textbf{\eqref{I}} follows directly from \eqref{eq:width_lemma_claim1}. 
	
	Now we show \textbf{(ii)}. Let $v^{*}\in V$ be the node for which \eqref{eq:width_lemma_claim1} holds. Let $V_{t}:=\sum_{v\in V} \delta( X_{t}(v^{*}) , X_{t}(v) )$. Then by \eqref{eq:width_lemma_claim1}, we can write 
	\begin{align}
		V_{t-1} - V_{t} &=\sum_{v\in V} \delta( X_{0}(v^{*}) , X_{t-1}(v) -t ) - \delta( X_{0}(v^{*}) , X_{t}(v) -t ) \\
		&= \sum_{v\in V}  \mathbf{1}( X_{t-1}(v) = X_{t}(v) ).
	\end{align}
	That is, $V_{t-1} - V_{t}$ equals the number of nodes that are pulled (delayed) during the update $X_{t-1}\mapsto X_{t}$. According to the second condition in  \eqref{eq:width_lemma_claim1}, the `most lagging node' $v^{*}$ blinks once every $\kappa$ iterations. Let $t^{*}\in \{0, \dots, \kappa-1\}$ denote the first time that $v^{*}$ blinks. By the third condition in \eqref{eq:width_lemma_claim1}, it follows that all nodes in $G$ blinks exaclty once during every time interval $\{t^{*} + (n-1)\kappa+1 ,\dots,  t^{*}+n \kappa\}$ of consecutive integers for every $n\ge 1$. (During $\{0, 1,\dots, t^{*}\}$,  nodes can be pulled by several times.) 
	
	Now fix $n\ge 1$ and suppose that $X_{s}$ is not synchronized where $s= n t^{*}$. Then we can choose a node $u^{*}$ such that $X_{t-1}(u^{*})$ is `most advanced' in $X_{t-1}$. That is, for each $v\in V$, $X_{s}(v) = X_{s}(v^{*}) - a$  for some $a\in \{0,1,\dots, \lfloor \kappa/2 \rfloor \}$. Since $G$ is connected and $X_{s}$ is not synchronized, there exists nodes $u',w'\in V$ such that $X_{s}(u')=X_{s}(u^{*})$ and $X_{s}(w')$ is lagging to $X_{s}(u')$ (i.e., $X_{s}(w') \in \{X_{s}(u')-\lfloor \kappa/2 \rfloor,\dots, X_{s}(u')-1 ,X_{s}(u') \} $). Recall that $w'$ blinks once in the time interval $\{s,s+1,\dots,s+\kappa  \}$. It follows that during the update $X_{s} \mapsto X_{s+\kappa}$, $u'$ must get pulled at least once (either by $w'$ or some other neighbor of $u'$). Hence, we have $V_{s+\kappa}\le V_{s}-1$. Noting that $V_{t^{*}} \le (\kappa/2) |V|$, it follows that $V_{t}=0$ for some $t\le t^{*} + \kappa(\kappa/2) |V| \le  (\kappa+1)(\kappa/2) |V| $. This shows \textbf{(ii)}. 
\end{proof}


\subsection{The branch-width lemma}

One of the key technial ingredients in analyzing FCA on trees is a local version of the width lemma (see Lemma \ref{lem:width}). Roughly speaking, if a small subtree in a tree has a restricted width less than $(\kappa/2) -1$, which is one less than what is required by the width lemma  (see Lemma \ref{lem:width}), then one can delete such subtree from the whole tree without affecting the dynamics of the rest. This enables inductive arguments. This observation has been used in \cite{lyu2015synchronization, lyu2017global} to analyze FCA and its continuous counterpart called the 4-coupling on trees.

A connected subgraph $S\subseteq G$ is called a \textit{$r$-star} if it has a vertex $v$, called the \textit{center}, such that all the other vertices of $S$ are leaves in $G$. A $r$-star $S$ is called a \text{$r$-branch} if the center of $S$ has only one neighbor in $G-S$, which we may call the \textit{root} of $S$. We may denote a $r$-branch by $B$ rather than by $S$. Note that branches are the smallest induced subgraphs of trees with a single vertex adjacent to their complement. Hence it is the smallest subgraph which gets minimal perturbation from outside and yet it should have simple internal dynamics. We say a dynamic $(X_{t})_{t\ge 0}$ on $G$ \textit{restricts on $H\subset G$} if the restriction $X_{t}\mapsto X_{t}|_{H}$ and transition map $\tau$ commute, i.e., the induced restricted dynamic $(X_{t}|_{H})_{t\ge 0}$ follows the same transition map on $H$. We say the dynamic $(X_{t})_{t\ge 1}$ on $G$ \textit{restricts on $H$ eventually} if there exists $r\ge 0$ such that $(X_{t})_{t\ge r}$ restricts on $H$.

\begin{figure*}[h]
	\centering
	\includegraphics[width=1 \linewidth]{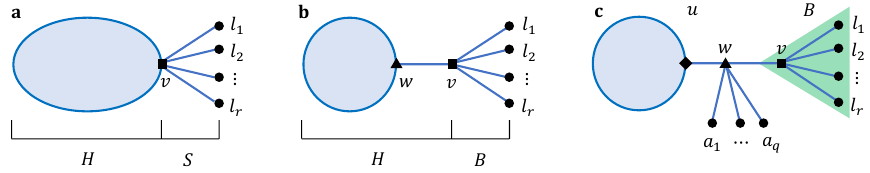}
	\caption{ \textbf{a} The $r$-star $S$ consists of the center $v$ and leaves $v_{1},\dots, v_{r}$. \textbf{b} The $r$-branch $B$ consists of the center $v$, leaves $v_{1},\dots, v_{r}$, and the root $w$. \textbf{c} The root $w$ of the $r$-star $B$ (in green) has leaf neivhbors $a_{1},\cdots,a_{q}$ and one external neighbor $u$. }
	\label{fig:branch_pic}
\end{figure*}         

The following `branch-width' lemma (see Lemma \ref{lem:branchwidth}) statest that, if $B$ is a branch in a connected graph $G$ and if $w_{B}(X_{0})<(\kappa/2)-1$, then $w_{B}(X_{t})<\kappa/2  $ for all $t\ge 0$ and the global dynamics restricts on the subgraph $G-B$. That is, the branch $B$ does not play any role in sustaining the dynamic on $G-B$.

\begin{lemma}[branch-width lemma, Lem. 3.4 in \cite{lyu2015synchronization}]\label{lem:branchwidth}
	Let $G=(V,E)$ be a graph with a $r$-branch $B$ rooted at some vertex $w\in V$. Let $v$ be the center of $B$ and and $l_{1},\cdots,l_{r}$, $r\ge 1$ be its leaves (see Figure \ref{fig:branch_pic}\textbf{b}). Let $H$ be the graph obtained from $G$ by deleting the leaves of this branch. Let $X_{0}$ be a $\kappa$-coloring on $G$ for any $\kappa\ge 3$ with $w_{B}(X_{0})<\kappa/2-1$. Then the following hold:
	\begin{description}[itemsep=0.1cm]
		\item[(i)] $v$ is ahead of all leaves of $B$ at some time $t_{0} \le \kappa(w_{B}(X_{0})+1)$ and $w_{B}(X_{t_{0}})\le w_{B}(X_{0})$;
		\item[(ii)] $v$ is ahead of all leaves of $B$ for all $t\ge t_{0}$, and $w_{B}(X_{t})\le w_{B}(X_{0})+1$ for all $t\ge 0$;
		\item[(iii)]	If $v$ is ahead of all leaves at $t=t_{0}$, then the dynamic $(X_{t})_{t\ge t_{0}}$ restricts on $H$; 
		\item[(iv)] If every $\kappa$-coloring on $H$ synchronizes, then $X_{0}$ synchronizes.
	\end{description}
\end{lemma}

Detailed proof can be found in \cite{lyu2015synchronization}, and here we give a brief sketch through an example. Suppose $\kappa=8$ and $k=3$. Since the coupling is inhibitory, the leaves of $B$ and the root $w$ only pull $u$ until it becomes the most lagging one in $B$. So eventually, we will have a situation where the branch width $w_{B}$ is still strictly less than $(\kappa/2)-1$, and the center $u$ is at most lagging in $B$. Now the root $w$ pulls $v$ at most once in every $\kappa$ iteration, increasing the branch width by 1. But since we have wiggle room on the branch width, the increased branch width is still small ($<\kappa/2$) and the leaves do not pull $u$ until its next blink. Then the center $u$ blinks and pulls all leaves, decreasing the branch width by 1. Hence the original branch width is recovered, and this scenario repeats over and over again. In this cycle, the leaves never pull the center, so the dynamics restrict on $H$.

\subsection{Prohibited configurations on a node with a leaf}

When a node $v$ has a leaf neighbor $x$, then whenever $v$ blinks in an FCA dynamic, $x$ cannot have the `opposite' color, which is $b(\kappa)+\lceil \kappa/2 \rceil \, (\mod \kappa)$, at any time  $t\ge \lceil \kappa/2 \rceil $. 

\begin{lemma}[Opposite color on a leaf is prohibited]
	\label{lem:opposite_leaf}
	Let $G=(V,E)$ be a graph with a node $v$. Suppose $x\in V$ is a neighbor of $v$ and has degree one. Let $(X_{t})_{t\ge 0}$ denote an arbitrary $\kappa$-color FCA trajectoy on $G$ for $\kappa \ge 3$. Then for every $t\ge \lceil \kappa/2 \rceil $, the following implication holds: 
	\begin{align}
		X_{t}(v)=b(\kappa)	 \qquad \Longrightarrow \qquad X_{t}(x)   \ne  b(\kappa)+\lceil \kappa/2 \rceil \, (\mod \kappa).
	\end{align}
\end{lemma}

\begin{proof}
	Suppose for contradiction that $t\ge \lceil \kappa/2 \rceil $, $X_{t}(v)=b(\kappa)$, and $X_{t}(x)  =  b(\kappa)+\lceil \kappa/2 \rceil \, (\mod \kappa)$. Then we back-track the dynamics for $\lceil \kappa/2 \rceil$ iterations from $X_{t}$ restricted on two nodes $v$  and $x$ as follows:
	\begin{align}
		\begin{matrix}
			X_{t}(v)\,|\, X_{t}(x) &\leftarrow	& \cdots & \leftarrow &X_{t-\lceil \kappa/2 \rceil +1 }(v) \,|\, X_{t-\lceil \kappa/2 \rceil + 1}(x) & \leftarrow &X_{t-\lceil \kappa/2 \rceil}(v) \,|\, X_{t-\lceil \kappa/2 \rceil}(x) \\
			b(\kappa)\,|\, X_{t}(x)  &\leftarrow  & \cdots & \leftarrow &0 \,|\,b(\kappa)+1 & \nleftarrow &\kappa-1 \,|\,b(\kappa) \\
		\end{matrix}
	\end{align}
	Note that $(X_{t}(v), X_{t-1}(v),\dots,  X_{t-\lceil \kappa/2 \rceil+1}(v),X_{t-\lceil \kappa/2 \rceil}(v) )=(b(\kappa),b(\kappa)-1,\dots,0,\kappa-1 )$ since nodes with colors between $0$ and $b(\kappa)$ are not pulled and the color prior to $0$ is $\kappa-1$. On the other hand, since $x$ is a leaf and since its only neighbor $v$ does not blink during $\{t-\lceil \kappa/2 \rceil, \dots, t-2,t-1\}$, its color keeps incrementing by one during this interval. This yields the following contradictory configuration: $X_{t-\lceil \kappa/2 \rceil +1 }(v)=0$,  $X_{t-\lceil \kappa/2 \rceil}(v)=\kappa-1$, and  $X_{t-\lceil \kappa/2 \rceil}(x)=b(\kappa)$. Indeed, since $x$ blinks at time $t-\lceil \kappa/2 \rceil$  and since $v$ is adjacent to $x$ and has color $\kappa-1$ at the same time, it's color at the next time $t-\lceil \kappa/2 \rceil+1$ must still be $\kappa-1$, not 0. To see more concrete examples, we provide the above back-tracking dynamics for the following four instances:
	\begin{align}
		\begin{matrix}
			(\kappa=3):\qquad  &	\text{$v$ | leaf $x$}	& \1\,|\, 0 &\leftarrow  & 0 \,|\, 2 &  \nleftarrow & 2\, |\, \1  \\
			(\kappa=4):\qquad & 	\text{$v$ | leaf  $x$}	& \1\,|\, 3 &\leftarrow  & 0 \,|\, 2	 & \nleftarrow & 3\, |\, \1  \\
			(\kappa=5):\qquad & 	\text{$v$ | leaf  $x$}	& \2\,|\, 0 &\leftarrow  & 1\,|\, 4	& \leftarrow  & 0\,|\, 3 & \nleftarrow & 4\, |\, \2  \\
			(\kappa=6):\qquad &	\text{$v$ | leaf $x$}	& \2\,|\, 5 &\leftarrow  & 1\,|\, 4	& \leftarrow  & 0\,|\, 3 & \nleftarrow & 5\, |\, \2.
		\end{matrix}
	\end{align}
\end{proof}

We remark that nothing prevents us to have the configurations prohibited by Lemma \ref{lem:opposite_leaf} at early times $t< \lceil \kappa/2 \rceil$. For instance, the initial configuration in the 5-color FCA example in Figure \ref{fig:branch_pic}\textbf{a} has a blinking node (of color 2) and a leaf neighbor of the `opposite color' $0$.

\subsection{Nodes with a small degree should blink infinitely often.}

In this subsection, we will make an important observation that a node with degree $<\kappa$ in a $\kappa$-color FCA dynamics must blink infinitely often. The following lemma is a quantitative version of this statement.  The idea is that, since FCA is an inhibitory pulse-coupled oscillator system, every node can blink at most once in every $\kappa$ iteration. Hence a node with degree $<\kappa$ can be pulled at most $\kappa-1$ times within each of its oscillation cycles. So such a node must advance its phase at least by one within each of its oscillation cycles.

\begin{lemma}[Blinking gaps and degree]
	\label{lem:blinking_degree}
	Let $(X_{t})_{t\ge 0}$ denote an arbitrary $\kappa$-color FCA trajectoy on a graph $G=(V,E)$ for $\kappa \ge 3$. Suppose $G$ has a node $v$ with degree $\deg(v)<\kappa$. Then $v$ should blink in $(X_{t})_{t\ge 0}$ once at least in $\lfloor \frac{\kappa^{2}}{ 2(\kappa - \deg(v)) } \rfloor + \kappa$ iterations. 
\end{lemma}

\begin{proof}
	Suppose for contradiction that there exists times $t_{0}$ and an integer $N\ge \lfloor \kappa^{2}/(\kappa-\deg(v)) \rfloor$ such that $v$ never blinks during the time interval $\{ t_{0},t_{0}+1 ,\dots, t_{0}+N \}$. We may assume that $v$ has a post-blinking color $>b(\kappa)$ at time $t_{1}:=t_{0}+b(\kappa)$. By the assumption, the color of $v$ must increment by no more than $\kappa-b(\kappa)$ during the time interval $\{ t_{1},t_{0}+1 ,\dots, t_{0}+N-b(\kappa) \}$, since otherwise $v$ will blink during the time interval $\{ t_{1},t_{1}+1 ,\dots, t_{0}+N \}$, contrary to our assumption. 
	
	Since $X_{t+1}(v) \in \{ X_{t}(v), X_{t}(v)+1  \} \, (\mod \kappa)$ for all $v\in V$ and $t\ge 0$, it follows that every node $v$ blinks at most once in every $\kappa$ iterations. Suppose  $\deg(v)<\kappa$. Hence no neighbor of $v$ blinks twice in the interval $\{t_{1},\dots,t_{1}+\kappa\}$. It follows that $v$ can be pulled at most $\deg(v)$ times during the time interval $\{t_{1},\dots,t_{1}+\kappa\}$. Hence the color of $v$ increments by $\kappa-\deg(v)$ $(\mod \kappa)$ during this interval. It follows that
	\begin{align}
		\kappa-b(\kappa)\ge 	\textup{total color increment during \textup{$\{ t_{1},t_{0}+1 ,\dots, t_{0}+N-b(\kappa) \}$}} & \ge \frac{t_{0}+N-b(\kappa)-t_{1} }{\kappa } (\kappa-\deg(v)).
	\end{align}
	It follows that 
	\begin{align}
		N &\le \kappa (\kappa - b(\kappa))/(\kappa-\deg(v))  + 2b(\kappa) \le \frac{\kappa^{2}}{ 2(\kappa - \deg(v)) } + \kappa. 
	\end{align} 
	This shows the assertion.
\end{proof}

\subsection{The relative circular representation}

For the forthcoming analysis, it is convenient to introduce a geometric representation of the FCA dynamics. Let $(X_{t})_{t\ge 0}$ be a $\kappa$-color FCA trajectory on a graph $G=(V,E)$.  Consider the induced dynamics $(Y_{t})_{t\ge 0}$, where $Y_{t}:V\cup\{\alpha\}=\mathbb{Z}_{\kappa}$ is the \textit{relative configuration} given by 
\begin{equation*}
	Y_{t}(x)=\begin{cases}
		b(\kappa)-t \, (\mod \kappa) & \text{if $x=\alpha$}\\
		X_{t}(x)-t+b(\kappa) \, (\mod \kappa) & \text{otherwise}.
	\end{cases}
\end{equation*} 
We refer to the value $Y_{t}(v)$ the \textit{phase} of $v$ at time $t$. Note that in the original dynamics, a node blinks whenever $X_{t}=b(\kappa)$, so in the relative dynamics $(Y_{t})_{t\ge 0}$, a node blinks whenever it has phase $-t \mod \kappa$. In this relative dynamics vertices keep the same phase until they get pulled, in which case they decrease their phase by 1. 

\begin{figure*}[h]
	\centering
	\includegraphics[width=0.9 \linewidth]{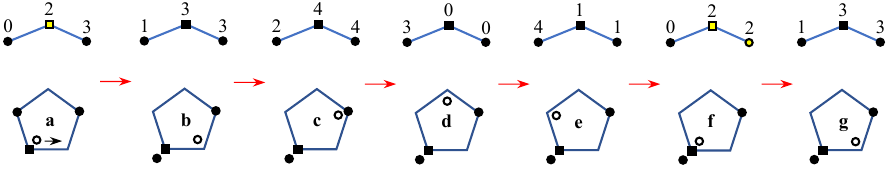}
	\caption{The top row shows the evolution of a $5$-color FCA dynamics on a path of three nodes, where red arrows indicate applications of the FCA transition rule. Nodes with blinking color $b(5)=2$ are indicated in yellow. The bottom row shows the same dynamic in relative circular representation. The pentagon represents the color space $\mathbb{Z}_{5}$, the open circle represents the activator, and the positions of the dots and the square relative to the activator indicate the colors of the corresponding nodes. }
	\label{fig:relative}
\end{figure*}         

A comparison between the original dynamics and the relative dynamics in case $\kappa=5$ is illustrated by example in Figure \ref{fig:relative}. The geometric representation of the relative FCA dynamics in the second row in Figure \ref{fig:relative} is what we call the \textit{relative circular representation}. The pentagon represents the phase space, which is the original color space $\mathbb{Z}_{5}$ modulo rotation, increasing in clockwise orientation. The open circle inside it, called the \textit{activator}, revolves around the phase space counterclockwise at unit speed, whose location at time $t$ is $-t \, (\mod 5)$. Hence a vertex blinks whenever it has the same phase as the activator (e.g., the middle node of color $2$ in Figure \ref{fig:relative}\textbf{a}). 

\section{ Proof of Theorem \ref{thm:tree_time} for $\kappa=4$ and 5 }
\label{sec:pf_tree_time_45}


In this section, we prove Theorem \ref{thm:tree_time} \textbf{\eqref{I}}. The case of $\kappa=3$ is covered by Theorem \ref{thm:tree_time}, so we will only prove the cases for $\kappa\in \{4,5\}$.

\subsection{Proof of Theorem \ref{thm:tree_time} for $\kappa=4$}

The key observation to make is that the 4-color FCA dynamics on a finite tree with maximum degree $\le 3$ `reduces' to a smaller subtree (one diameter less) in every finite number of iterations.

\begin{lemma}[4-color FCA local dynamics on a branch]
	\label{lem:4FCA_branch}
	Let $(X_{t})_{t\ge 0}$ be a 4-color FCA trajectory on a graph $G=(V,E)$, which has a $r$-branch $B$ for $r\ge 1$ with center $v$, root $w$, and leaves $l_{1},\dots,l_{r}$ (see Figure \ref{fig:branch_pic}\textbf{b}). Suppose $v$ blinks at some time $t_{0}\ge 2$. Then the following hold:
	\begin{description}
		\item[(i)] Suppose the dynamic $(X_{t})_{t\ge t_{0}+ 22}$ does not restrict on $G-B$ where the leaves in $B$ have the same colors. Then the local dynamic $(X_{t})_{t_{0}+8 \le t \le t_{0}+16 }$ restricted on the branch $B$ is given by the following sequence:
		\begin{align}\label{eq:lem_4FCA_branch_w}
			\begin{matrix}
				\textup{time} & t_{0}+8 &&&&t_{0}+12&&& t_{0}+15& t_{0}+16 \\
				\text{center $v$ | leaves}		&	\1|\1 0 & 2 | 2 \1 & 2| 3 2 & 3| 0 3 & 3 | \1 0 & 3| 2 \1 & 3| 3 2 & 0 | 0 3 & \1| \1 0 \\
				\text{root $w$\hspace{1.2cm}}& 3 & 3 & 0 &  \1  & 2 &* & *& * &   3
			\end{matrix} 
		\end{align}
		where $*\in \{2,3\}$. In particular $w$ is pulled exactly three times from time $t_{0}+12$ to $t_{0}+16$. In particular, $w$ must have a degree of at least four. 
		
		\item[(ii)] For a fixed integer $m\ge 3$, suppose the dynamic $(X_{t})_{t\ge t_{0}+ 8m}$ does not restrict on $G-B$ where the leaves in $B$  have the same color. Then the local dynamic on $B$ on time intervals $\{ t_{0}+8j,\dots, t_{0}+8(j+1) \}$ for $j=1,\dots,m-2$ are identical to \eqref{eq:lem_4FCA_branch_w}. 
	\end{description}	
\end{lemma}

\begin{proof}
	Suppose the dynamic $(X_{t})_{t\ge t_{0}+ 22}$ does not restrict on $G-B$.  First, by Lemma \ref{lem:opposite_leaf}, there is no leaf in $B$ with color 0 at time $t_{0}$. Then there is a total of seven possible local configurations on $B$ at time $t_{0}$, as shown in Figure \ref{fig:4FCA_branch1}\textbf{a}-\textbf{g} in relative circular representation. Note that the red shaded configurations in Figure \ref{fig:4FCA_branch1}\textbf{h}-\textbf{i} have branch-width 0. 
	The local configurations \textbf{a}-\textbf{d} lead to the ones in \textbf{h}-\textbf{i} in at most six iterations. Hence in this case, 
	by Lemma \ref{lem:branchwidth}, the dynamic $(X_{t})_{t\ge t_{0}+ 6}$ restricts on $G-B$. Next, the ones in \textbf{e}-\textbf{g} may lead to the ones in \textbf{i} in six iterations, in which case by Lemma \ref{lem:branchwidth} the dynamic $(X_{t})_{t\ge t_{0}+ 6}$ restricts on $G-B$. In all other cases, the local configuration on $B$ becomes the one in \textbf{k} in eight iterations. It follows that the local configuration at branch $B$ at time $t_{0}+8$ is the one in \textbf{k}. By applying the same argument, since we are assuming $(X_{t})_{t\ge t_{0}+ 22}$ does not restrict on $G-B$, the branch $B$ must  undergo the transition \textbf{k} $\rightarrow$ \textbf{k} $\rightarrow$ \textbf{k} until time $t_{0}+24$.

	\begin{figure*}[h]
		\centering
		\includegraphics[width=0.7 \linewidth]{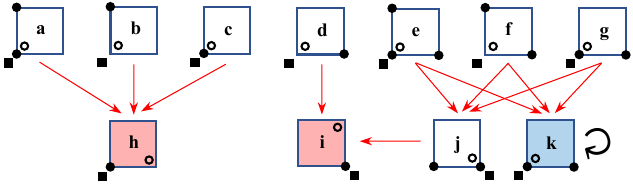}
		\caption{ Possible local relative configurations on the branch $B$ with center $\blacksquare$ and leaves $\bullet$ at times $t\ge 2$ are shown in \textbf{a}-\textbf{g}. The red shaded configurations \textbf{h}-\textbf{i} have branch-width zero. The local configurations \textbf{a}-\textbf{d} lead to the ones in \textbf{h}-\textbf{i} in at most six iterations. The ones in \textbf{e}-\textbf{g} may lead to the ones in \textbf{i} in five iterations or to the ones in \textbf{k} in eight iterations. The one in \textbf{g} may lead to \textbf{k} in eight iterations, which is identical to itself. 
		}
		\label{fig:4FCA_branch1}
	\end{figure*}  
	
	Now we forward-track the local transition \textbf{k} $\rightarrow$ \textbf{k} from time $t_{0}+8$ to time $t_{0}+16$ as shown in Figure \ref{fig:4FCA_branch2}. The transitions \textbf{k}$\rightarrow$\textbf{k}$_{1}$$\rightarrow$\textbf{k}$_{3}$ and \textbf{k}$\rightarrow$\textbf{k}$_{1}$$\rightarrow$\textbf{k}$_{4}$$\rightarrow$\textbf{k}$_{6}$$\rightarrow$\textbf{k}$_{7}$$\rightarrow$\textbf{k}$_{8}$ both contradicts the assumption that the dynamic $(X_{t})_{t\ge t_{0}+ 24}$ does not restrict on $G-B$. Hence the local dynamic during this period should be given by \textbf{k}$\rightarrow$\textbf{k}$_{1}$$\rightarrow$\textbf{k}$_{4}$$\rightarrow$\textbf{k}$_{6}$$\rightarrow$\textbf{k}$_{9}$ in Figure \ref{fig:4FCA_branch2}. This requires the center $v$ (represented as $\blacksquare$ in Figure \ref{fig:4FCA_branch2}) to be pulled by the center $w$ during the transition \textbf{k}$_{1}$$\rightarrow$\textbf{k}$_{4}$. 
	\begin{figure*}[h]
		\centering
		\includegraphics[width=0.7 \linewidth]{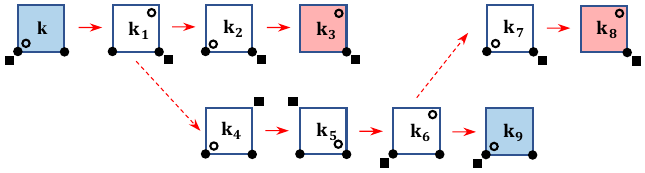}
		\caption{ Possible local dynamics on the branch $B$ with center $\blacksquare$ and leaves $\bullet$ starting from the local configurations in \textbf{a}, which identical to the one in Figure \ref{fig:4FCA_branch1}\textbf{k}. For the transitions indicated by dashed arrows  (\textbf{k}$_{1}$$\rightarrow$\textbf{k}$_{4}$ and \textbf{k}$_{6}$$\rightarrow$\textbf{k}$_{7}$), the center $v$ must be pulled by the root $w$. 
		}
		\label{fig:4FCA_branch2}
	\end{figure*}  
	It follows that the local dynamics on the branch $B$ from time $t_{0}+8$ to $t_{0}+24$ is given by one of the following two sequences in \eqref{eq:4FCA_branch_seq}, depending on when the root $w$ pulls the center $v$. In fact, it must only be the first sequence in \eqref{eq:4FCA_branch_seq} since the second one in \eqref{eq:4FCA_branch_seq}  is not possible:  Back-tracking for two iterations from time $t_{0}+10$ from the third instance leads to a contradiction that $w$ with color $3$ is pulled by $v$ with color $1$ at time $t_{0}+8$ cannot lead to $w$ with color $1$ at time $t_{0}+9$.
	\begin{align}\label{eq:4FCA_branch_seq}
		\begin{cases}
			\begin{matrix}
				\textup{time} & t_{0}+8 &&&&t_{0}+12&&& t_{0}+15& t_{0}+16 \\
				\text{center $v$ | leaves}		&	\1|\1 0 & 2 | 2 \1 & 2| 3 2 & 3| 0 3 & 3 | \1 0 & 3| 2 \1 & 3| 3 2 & 0 | 0 3 & \1| \1 0 \\
				\text{root $w$\hspace{1.2cm}}& 3 & 3 & 0 &  \1  & * &* & *& * &   *\\
				\text{center $v$ | leaves}		&	\1|\1 0 & 2 | 2 \1 & 2| 3 2 & 2| 0 3 & 3 | \1 0 & 3| 2 \1 & 3| 3 2 & 0 | 0 3 & \1| \1 0 \\
				\text{root $w$\hspace{1.2cm}}& 3 & 0 & \1 &  2  & * &* & *& * & * 
			\end{matrix}
		\end{cases}
	\end{align}
	Finally, applying the same analysis as above to the local dynamics on $B$ from time $t_{0}+16$ to $t_{0}+24$, it follows that $w$ must have the same color 3 at time $t_{0}+16$ as it has at time $t_{0}+8$. Hence the local dynamcis on $B$ from  time $t_{0}+8$ to $t_{0}+16$ is as described in \eqref{eq:lem_4FCA_branch_w}, which requires $w$ to be pulled three times from time $t_{0}+12$ to $t_{0}+16$. Since $w$ is not pulled by $v$  during this period in \eqref{eq:lem_4FCA_branch_w}, and since every node blinks at most once in four iterations in any 4-color FCA dynamics, it follows that $w$ must have at least three other neighbors than $v$. This shows \textbf{(ii)}. Repeating the same argument shows \textbf{(iii)}.
\end{proof}

It is now easy to derive Theorem  \ref{thm:tree_time} for $\kappa=4$ by iteratively applying Lemma \ref{lem:4FCA_branch}.

\begin{proof}[\textbf{Proof of Theorem \ref{thm:tree_time} for $\kappa=4$}] 
	We show the assertion by induction on the diameter $d$ of the underlying tree $T=(V,E)$ with maximum degree at most $3$. If $d=1$, then $T$ is a path of two nodes and the assertion can be directly verified. Now assume $d\ge 2$. As an induction hypothesis, suppose that if $T'$ is a finite tree with maximum degree at most $3$ and its shortest-path diameter $d'$ is less than $d$, then arbitrary 4-color FCA trajectory on $T'$ synchronizes in $C(d'+1)$ iterations, where $C>0$ is a constant to be determined. Consider an arbitrary 4-color FCA trajectory $(X_{t})_{t\ge 0}$ on a finite tree $T=(V,E)$ with maximum degree at most $3$ with diameter $d$. By Lemma \ref{lem:blinking_degree}, all nodes in $T$ blink at least once by time 12. Then by Lemma \ref{lem:4FCA_branch} \textbf{\eqref{I}}, it follows that for every branch $B$ in $T$, the dynamic $(X_{t})_{t\ge  t_{0}+22}$ restricts on $G-B$ and the leaves in $B$ have the same color. In particular, the dynamic $(X_{t})_{t\ge   t_{0}+22}$ restricts onto the subtree $T'$ of $T$ obtained by deleting all branch of $T$. Note that $T'$ has a diameter at most $d-2 $. By the induction hypothesis, $(X_{t})_{t\ge   t_{0}+22}$ restricted on $T'$ synchronizes by time $t_{1}:= t_{0}+22+C(d-1)$. In particular, it follows that the roots of all branches in $T$ are synchronized and blink exactly once in four iterations after time $t_{1}$. Furthermore, the leaves in every branch in $T$ have the same color after time $t_{1}$. 
	
	Combining with the analysis of the local dynamics in Figure \ref{fig:4FCA_branch1}, it follows the dynamic $(X_{t})_{t\ge t_{1}+5}$ restricts onto the subtree $T''$ of $T$ obtained by deleting all leaves in $T$. Then by identifying all nodes in $T'$, which have identical dynamics after time $t_{1}$, we can view the dynamic $(X_{t})_{t\ge t_{1}+5}$ restricted on $T''$ as a 4-color FCA dynamic on a star with leaves being the center of every branch in $T$. Again by using the analysis in Figure \ref{fig:4FCA_branch1}, such star dynamic synchronizes by time $t_{1}+15$. By using a similar argument, this time identifying all nodes in $T''$, it follows that the entire configuration $X_{t}$ on $T$ is synchronized by time $t_{1} + 25$. Thus $X_{t}$ is synchronized for $t\ge t_{0}+47+C(d-1)$. Noting that $t_{0}\le 12$, we have that $X_{t}$ is synchronized for $t\ge 59 + C(d-1)$. Now setting $C=30$, we can complete the induction. 
\end{proof}

\subsection{Proof of Theorem \ref{thm:tree_time} for $\kappa=5$}

We prove Theorem  \ref{thm:tree_time} for $\kappa=5$ in this section.  

\begin{lemma}[5-color FCA local dynamics on a branch]
	\label{lem:5FCA_branch}
	Let $(X_{t})_{t\ge 0}$ be a 5-color FCA trajectory on a graph $G=(V,E)$, which has a $r$-branch $B$ for $r\ge 2$ with center $v$, root $w$, and leaves $l_{1},\dots,l_{r}$ (see Figure \ref{fig:branch_pic}\textbf{b}). Suppose $v$ blinks at some time $t_{0}\ge 3$. Then the followings hold:
	\begin{description}[noitemsep]
		\item[(i)]  Suppose the dynamic $(X_{t})_{t\ge t_{0}+ 15}$ does not restrict on $G-B$. Then the local dynamic $(X_{t})_{t_{0}\le t \le t_{0}+ 8}$ restricted on the branch $B$ is given by the following sequence:
		\begin{align}\label{eq:5FCA_local_dynamic_recurrent}
			\begin{matrix}
				\textup{time} & t_{0} &&t_{0}+2&t_{0}+3&&&&t_{0}+7& t_{0}+8 && t_{0}+10  \\
				\text{center $v$ | leaves}		&	\2|14 & 3|\2 4 & 3|30 & 3|41 & 4|0 \2 & 4|1 3 & 0|\2 4 & 1|30 & \2|4 1 & 3 | 4 \2 &  3 | 0 3 \\
				\text{root $w$\hspace{1.2cm}}& 0 & 1 & \2 &  3  & * &* & *& 4 &  0  & 1 & \2
			\end{matrix} 
		\end{align}
		where $*\in \{3,4\}$.  In particular $w$ is pulled exactly three times from time $t_{0}+3$ to $t_{0}+7$.	
		
		\item[(ii)] Suppose the dynamic $(X_{t})_{t\ge t_{0}+ 8m+7}$ does not restrict on $G-B$ for some integer $m\ge 1$. Then the local dynamic on $B$ on time intervals $\{ t_{0} + 8j ,\dots, t_{0} + 8(j+1) \}$ for $j=1,\dots,m$  are identical to the one in \eqref{eq:5FCA_local_dynamic_recurrent}. 
	\end{description}	
\end{lemma}

\begin{proof}
	By Lemma \ref{lem:opposite_leaf}, there is no leaf in $B$ with color 0 at time $t_{0}$. Then there is a total of 15 possible local configurations on $B$ at time $t_{0}$, as shown in Figure \ref{fig:5FCA_branch} in relative circular representation. Note that the red shaded configurations in Figure \ref{fig:5FCA_branch}\textbf{e}-\textbf{o} have branch-width strictly smaller than $1$ in one iteration. Hence if the local configuration on $B$ at time $t_{0}$ is one of \textbf{e}-\textbf{o}, then by Lemma \ref{lem:branchwidth}, the dynamic $(X_{t})_{t\ge t_{0}+ 1}$ restricts on $G-B$. 
	
	\begin{figure*}[h]
		\centering
		\includegraphics[width=0.8 \linewidth]{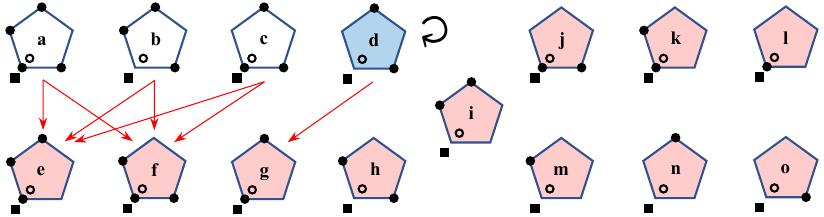}
		\caption{ Possible local relative configurations on the branch $B$ with center $\blacksquare$ and leaves $\bullet$ at times $t\ge 3$. The red shaded configurations \textbf{e}-\textbf{o} have branch-width strictly smaller than $1= \lfloor (\kappa/2)-1\rfloor$ in one iteration. The local configurations \textbf{a}-\textbf{c} lead to \textbf{e}-\textbf{f} (up to rotation). The blue-shaded local configuration \textbf{d} is the only possible recurrent one. 
		}
		\label{fig:5FCA_branch}
	\end{figure*}  
	
	Next, the local configurations in Figure \ref{fig:5FCA_branch}\textbf{a}-\textbf{c} lead to the ones in \textbf{e}-\textbf{f} in at most 11 iterations. Forward tracking the local dynamics from  \textbf{a}-\textbf{c} is shown in Figure \ref{fig:5FCA_branch1}. Hence, if the local configuration on $B$ at time $t_{0}$ is one of the ones in Figure \ref{fig:5FCA_branch}\textbf{a}-\textbf{c}, then the dynamic $(X_{t})_{t\ge t_{0}+ 12}$ restricts on $G-B$. 
	\begin{figure*}[h]
		\centering
		\includegraphics[width=0.8 \linewidth]{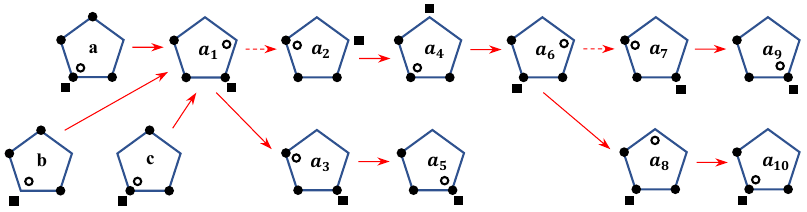}
		\caption{ Possible local dynamics on the branch $B$ with center $\blacksquare$ and leaves $\bullet$ starting from the local configurations in Figure \ref{fig:5FCA_branch}\textbf{a}-\textbf{c}. For the transitions indicated by dashed arrows  ($\mathbf{a}_{1}\rightarrow \mathbf{a}_{2}$ and $\mathbf{a}_{6}\rightarrow \mathbf{a}_{7}$), the center $v$ must be pulled by the root $w$. 
		}
		\label{fig:5FCA_branch1}
	\end{figure*}  
	
	In Figure \ref{fig:5FCA_branch2}, we analyze the detailed local dynamics on the branch $B$ starting from the local configuration in Figure \ref{fig:5FCA_branch}\textbf{d}. One the one hand, the transition $\textbf{d}\rightarrow \textbf{d}_{3}\rightarrow \textbf{d}_{4}$ shows that Figure \ref{fig:5FCA_branch}\textbf{d} could lead to \ref{fig:5FCA_branch}\textbf{g} in 6 iterations. On the other hand, the transitions $\mathbf{d}\rightarrow \mathbf{d}_{i}$ for $i=1,2$ require the center $v$ to be pulled by the root $w$. It then follows the transition $\textbf{d}_{5}\rightarrow \textbf{d}_{8}$ and $\textbf{d}_{8}$ is rotation symmetric to Figure \ref{fig:5FCA_branch}\textbf{d}. 
	\begin{figure*}[h]
		\centering
		\includegraphics[width=0.6 \linewidth]{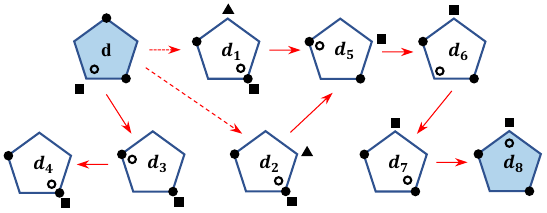}
		\caption{ Possible local dynamics on the branch $B$ with center $\blacksquare$,  leaves $\bullet$, and roote $\blacktriangle$ starting from the local configurations in Figure \ref{fig:5FCA_branch}\textbf{d}. The dashed transitions $d\rightarrow d_{1}$ and $d\rightarrow d_{2}$ requires the center $v$ to be pulled by the root $w$. 
		}
		\label{fig:5FCA_branch2}
	\end{figure*}

	Now suppose that the dynamic $(X_{t})_{t\ge t_{0}+ 14}$ does not restricts on $G-B$. Then by the previous discussion, we need to have the transition $\mathbf{d}\rightarrow \mathbf{d}$ from time $t_{0}$ to time $t_{0}+8$. Furthermore, from time $t_{0}+8$, we cannot have the transition $\mathbf{d}\rightarrow \mathbf{g}$ in Figure \ref{fig:5FCA_branch}, which takes six iterations, since then $(X_{t})_{ t\ge t_{0}+15}$ restricts on $G-B$. Thus we need to have the transition $\mathbf{d}\rightarrow \mathbf{d}$ again from time $t_{0}+8$ to time $t_{0}+16$. Reading off the dynamics on $B$ corresponding to the two transitions $\textbf{d}\rightarrow \textbf{d}_{i}\rightarrow \textbf{d}_{5}\rightarrow \textbf{d}_{8}$ for $i=1,2$  in Figure \ref{fig:5FCA_branch2} show that the local dynamic $(X_{t})_{t_{0}\le  t\le t_{0}+ 16}$ restricted on the branch $B$ is given by one of the  following two sequences:
	\begin{align}
		\begin{cases}
			\begin{matrix}
				\textup{time} & t_{0} &&t_{0}+2&&t_{0}+4&&&& t_{0}+8 & && t_{0}+16 \\
				\text{center $v$ | leaves}		&	\2|14 & 3|\2 4 & 3|30 & 3|41 & 4|0 \2 & 4|1 3 & 0|\2 4 & 1|30 & \2|4 1 & 3 | 4 \2&\cdots & \2 |  1 4 \\
				\text{root $w$\hspace{1.2cm}}& 0 & 1 & \2 &  3 \\
				\text{center $v$ | leaves}		&	\2|14 & 3|\2 4 & 3|30 & 4|41 & 4|0 \2 & 4|1 3 & 0|\2 4 & 1|30 & \2|4 1 & 3 | 4 \2&\cdots & \2 |  1 4 \\
				\text{root $w$\hspace{1.2cm}}&4 & 0& 1 & \2 & 3
			\end{matrix}
		\end{cases}
	\end{align}
	Note that the second sequence is impossible due to the first transition: The root $w$ must be pulled by the center $v$ with color 2 at time $t_{0}$, but it does not. It follows that $w$ must have color 2 at time $t_{0}+2$. By repeating the same argument, we also have that $w$ has color $2$ at time $t_{0}+10$. Then \textbf{\eqref{I}} follows. Lastly, \textbf{(ii)} follows easily by repeating the previous argument. 
\end{proof}

The following lemma in conjunction with Lemma \ref{lem:blinking_degree} shows the `constant-time reduction' of 5-color FCA on finite trees.

\begin{lemma}[5-color FCA local dynamics at the root of a branch]
	\label{lem:5FCA_branch_root}
	Let $(X_{t})_{t\ge 0}$ be a 5-color FCA trajectory on a graph $G=(V,E)$, which has a $r$-branch $B$ for $r\ge 2$ with center $v$, root $w$, and leaves $l_{1},\dots,l_{r}$.  Further assume that the neighbors of the root $w$ are $u,v,a_{1},\dots,a_{q}$ for some leaves $a_{1},\dots,a_{q}$ in $G$ and some node $u$ in $G$ (see Figure \ref{fig:branch_pic}\textbf{c}). Suppose $v$ blinks at some time $t_{0}\ge 3$. Then the dynamic $(X_{t})_{t\ge t_{0}+ 15}$ restrict on $G-B$. 
\end{lemma}

\begin{proof}
	Suppose for that the dynamic $(X_{t})_{t\ge t_{0}+ 15}$ does not restrict on $G-B$. Then by Lemma \ref{lem:5FCA_branch}, the local dynamics on the branch $B$ and its root $w$ from time $t_{0}$ to time $t_{0}+8$ is given by \eqref{eq:5FCA_local_dynamic_recurrent}, which requires the following: (1) $w$ blinks at times $t_{0}+2$ adn $t_{0}+8$; and (2) $w$ is pulled by three times from time $t_{0}+3$ to time $t_{0}+7$; and (3) $w$ is not pulled by the center $v$. In particular, (2) and (3) implies that (4) $w$ is pulled at least twice by its leaf neighbors $a_{1},\dots,a_{q}$  from time $t_{0}+3$ to time $t_{0}+7$. 
	\begin{figure*}[h]
		\centering
		\includegraphics[width=0.9 \linewidth]{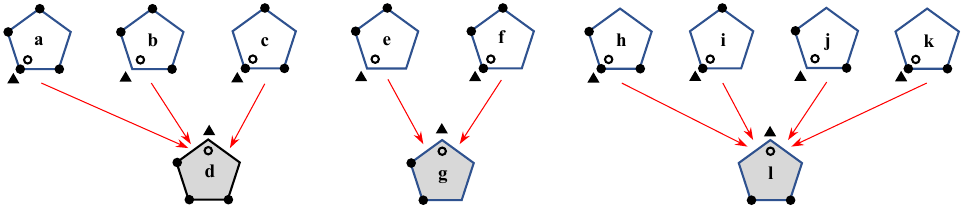}
		\caption{ Possible local dynamics on the root $\blacktriangle$ of the branch $B$ and the leaf neighbors $\bullet$ $a_{1},\dots,a_{q}$ of $w$ starting from the local configurations where $w$ blinks and $a_{1},\dots,a_{q}$ assume at least two distinct phases. Note that $w$ must blink after exactly $8$ iterations. The local dynamics lead to the ones in \textbf{d}, \textbf{g}, and \textbf{l}, which are prohibited to occur after time $t\ge 3$  by Lemma \ref{lem:opposite_leaf}. 
		}
		\label{fig:5FCA_branch_root1}
	\end{figure*}  
	
	In Figure \ref{fig:5FCA_branch_root1},  we analyze all possible local dynamics on $w$ and its leaves $a_{1},\dots,a_{q}$ from time $t_{0}+2$ to time $t_{0}+10$, assuming at least two distinct colors for $a_{1},\dots,a_{q}$ (note that color $4$ prohibited by Lemma \ref{lem:opposite_leaf}). At time $t_{0}+10$, the local configuration on $w$ and its leaf neighbors $a_{1},\dots,a_{q}$ leads to either of the ones in \textbf{d}, \textbf{g}, and \textbf{l}, which are prohibited to occur after time $t\ge 3$  by Lemma \ref{lem:opposite_leaf}. It follows that $a_{1},\dots,a_{q}$ have exacly one distinct state at time $t_{0}+2$. But this contradicts the requirement that (4) $w$ is pulled by its leaf neighbors  $a_{1},\dots,a_{q}$  at least twice from time $t_{0}+3$ to time $t_{0}+7$.  This is a contradiction. Thus we must conclude that the dynamic $(X_{t})_{t\ge t_{0}+ 15}$  restricts on $G-B$. 
\end{proof}

Now we prove Theorem \ref{thm:tree_time} for $\kappa=5$. The argument is similar to the one for Theorem \ref{thm:tree_time} for $\kappa=4$. 

\begin{proof}[\textbf{Proof of Theorem \ref{thm:tree_time} for $\kappa=5$}] 
	We show the assertion by an induction on the diameter $d$ of the underlying tree $T=(V,E)$. The assertion can easily be verified if $d=1$, in which case $T$ is a star of at most four leaves. Now assume $d\ge 1$. As an induction hypothesis, suppose that if $T'$ is a finite tree with maximum degree at most $4$ and its shortest-path diameter $d'$ is less than $d$, then arbitrary 5-color FCA trajectory on $T'$ synchronizes in $C(d'+1)$ iterations, where $C$ is a constant to be determined. Now consider an arbitrary 5-color FCA trajectory $(X_{t})_{t\ge 0}$ on a finite tree $T=(V,E)$ with maximum degree at most $4$ with diameter $d$. Choose the root node $v_{0}$ of $T$ so that viewed as a rooted tree with root $v_{0}$, the maximum depth of the nodes in $T$ is at most half of the diameter of $T$. That is, denoting $d$ to be the diameter of $T$, every node $v$ in $T$ is at most $\lceil (d'+1)/2 \rceil$ shortest-path distance away from $v_{0}$. 
	
	Let $l$ denote any leaf in $T$ that is farthest away from the root $v_{0}$. Let $v$ denote its paraent, and let $w$ denote the parent of $v$. Then it follows that, if $x$ is a neighbor of $w$ that is not the parent of $w$, then $x$ is either a leaf or the center of a branch in $T$. Recall that $w$ blinks at some time $t_{0}\le 17$ by Lemma \ref{lem:blinking_degree}. Let $B_{1},\dots,B_{m}$ denote all branches rooted at $w$. We first suppose that the dynamic $(X_{t})_{t\ge t_{0}+ 15}$ restricts on $T-(B_{1}\cup \dots \cup B_{m})$. If this is not the case, then for each $i\in \{1,\dots, m\}$ , either the dynamic $(X_{t})_{t\ge t_{0}+ 15}$ restricts on $T-B_{i}$ or the local configuration on $B_{i}$ at time $t_{0}$ should all be identical to the first one in \eqref{eq:5FCA_local_dynamic_recurrent}. It then follows that for all branches $B_{i}$ such that the dynamic $(X_{t})_{t\ge t_{0}+ 15}$ does not restirct on $T-B_{i}$, the local dynamic on $B_{i}$ after time $t\ge t_{0}$ are all identical. Hence in the dynamic $(X_{t})_{t\ge t_{0}+15}$, $w$ has at most one branch rooted at itself of distinct local dynamics and leaf neighbors, except the pattern of $w$. Thus by Lemma \ref{lem:5FCA_branch_root}, the dynamic $(X_{t})_{t\ge t_{0}+15}$ restricts on $T-(B_{1}\cup \dots \cup B_{m})$. This holds for all nodes $w$ that are distance two away from leaves in $T$ farthest from $v_{0}$. It follows that the dynamic $(X_{t})_{t\ge t_{0}+15}$ restricts to a subtree $T'$ of $T$ with diameter at most $d-2$. But by the induction hypothesis, $X_{t}|_{T'}$ is synchronized for all $t \ge t_{0}+15+C(d-1)$. Then the dynamic $(X_{t})_{t\ge t_{0}+15+C(d-1)}$ on $T$ is equivalent to a dynamic on a star with at most three leaves. Reusing the analysis in the proof of Lemma \ref{lem:5FCA_branch}, it is easy to see that it takes at most $15$ iterations to synchronize such a dynamic. Altogether, we conclude that $X_{t}$ for $t\ge t_{0}+30+C(d-1)$ is synchronized. This completes the induction. Recalling $t_{0}\le 17$, it follows that $X_{t}$ for $47+C(d-1)$ is synchronized. Then the induction is completed by setting $C=24$. 
\end{proof}

\section{FCA on finite trees when $\kappa\ge 7$.}
\label{sec:pf_FCA_7}

In this section, we prove the following result, which implies Corollary \ref{cor:phase_transition} for $\kappa\ge 7$:

\begin{customthm}{2.1}\label{thm_tree_nonsync}
	Let $\kappa\ge 7$ be an integer. 
	\begin{description}[noitemsep]
		\item[(i)] If $\kappa=2m-1\ge 3$ is odd, there exists a finite tree $T=(V,E)$ with maximum degree $m$ and a  $\kappa$-coloring $X_{0}:V\rightarrow \mathbb{Z}_{\kappa}$ such that $X_{t}$ is non-synchronizing whose period divides $3\kappa^{2}+\kappa$. 
		\item[(ii)] There exists a tree of maximum degree 4 with a non-synchronizing 8-coloring which is 60 periodic. 
		
		\item[(iii)] If $\kappa=2m\ge 10$ is even, there exists a finite tree $T=(V,E)$ with maximum degree $m+1$ and a $\kappa$-coloring $X_{0}:V\rightarrow \mathbb{Z}_{\kappa}$ such that $X_{t}$ is non-synchronizing whose period divides $5\kappa^{2}+\kappa$.   
	\end{description}
\end{customthm}

\begin{proof}
	\begin{description}
		\item{\eqref{I}} Note that $b(\kappa)=m-1\ge 3$. Let $p,q$ be integers such that $p,q\ge 2$ and $p+q=m$. Let $T=(V,E)$ be a star with $m$ leaves $v_{1},\cdots,v_{p},u_{1},\cdots,u_{q}$ and center $w$. Define a relative $\kappa$-coloring $Y_{0}:V\cup\{\alpha\}\rightarrow \mathbb{Z}_{\kappa}$ on $T$ by 
		\begin{equation*}
			Y_{0}(x) = 
			\begin{cases}
				m-1 & \text{if $x\in \{w,\alpha\}$} \\
				i  & \text{if $x=v_{i}$ for some $1\le i \le p$} \\
				m+j  & \text{if $x=u_{j}$ for some $1\le j \le q$}.
			\end{cases}
		\end{equation*}  
		To show the assertion for $Y_{t}$, we claim that 
		\begin{equation*}
			Y_{3\kappa+1}(x) = Y_{0}(x)-1 \quad \text{$\forall x\in V$.} 
		\end{equation*}
		We begin with observing that in the first $\kappa$ iterations, $w$ blinks once at time $0$ to pull each of $u_{j}$'s, and is pulled by each of $v_{i}$'s followed by each of $v_{j}$'s. Hence we have 
		\begin{equation*}
			Y_{\kappa}(x) = 
			\begin{cases}
				-1 & \text{if $x=w$} \\
				i  & \text{if $x=v_{i}$ for some $1\le i \le p$} \\
				m+j-1  & \text{if $x=u_{j}$ for some $1\le j \le q$}.
			\end{cases}
		\end{equation*}  
		Then since $0\le \delta_{\kappa}(v_{i},w) \le m$ for each $1\le i \le p$, $w$ is not pulled by any of the $u_{i}$'s during $[\kappa,\kappa+m]$. Then it pulls all of $u_{i}$'s at time $t=\kappa+m$, and is pulled by each of $u_{j}$'s during $[\kappa+m,2\kappa]$. This yields  
		\begin{equation*}
			Y_{2\kappa}(x) = 
			\begin{cases}
				-1-q & \text{if $x=w$} \\
				i-1  & \text{if $x=v_{i}$ for some $1\le i \le p$} \\
				m+j-1  & \text{if $x=u_{j}$ for some $1\le j \le q$}.
			\end{cases}
		\end{equation*}   
		Finally, $\delta_{2\kappa}(w,v_{p})=m$, so $w$ is pulled by each of $v_{j}$'s during $[2\kappa,2\kappa+m-1]$ so that $Y_{2\kappa+m-1}(w)\equiv -1-q-p \equiv m-2\mod \kappa$. Then $\delta_{2\kappa+m-1}(w,u_{q})\equiv m+a-2\in [m,\kappa-1]$, so $w$ is not pulled by any of $u_{j}$'s during $[2\kappa+m-1,3\kappa]$. Moreover, $Y_{3\kappa+1}(w)=Y_{3\kappa+1}(\alpha)=m-2$. This shows the claim.   
		
		\item{(ii)} Let $T=(V,E)$ be a tree where $V=\{v_{1},\cdots,v_{8}\}$ and edges are determined by  $N(v_{4})=\{v_{1},v_{2},v_{3},v_{5}\}$ and  $N(v_{5})=\{v_{4},v_{6},v_{7},v_{8}\}$. Note that $T$ has a maximum degree 4. Let $X_{0}:V\rightarrow \mathbb{Z}_{8}$ be the initial 8-coloring on $T$ defined by 
		\begin{equation*}
			(X_{0}(v_{1}),\cdots,X_{0}(v_{8})) = (1,5,7,5,6,0,3,6).
		\end{equation*}
		Then it is straightforward to check that $X_{t}$ satisfies the assertion. 
		
		\item{(iii)} Note that $b(\kappa)=m-1$. Let $T=(V,E)$ be a star with center $w$ and $L=V\setminus \{w\}$ the set of all leaves. For each integers $p,q,r\ge 2$ such that $p+q+r=m+1$, let $\mathcal{Y}_{p,q,r}$ be the set of all relative $\kappa$-colorings $Y^{p,q,r}$ on $T$ with the following properties: 
		\begin{description}[noitemsep]
			\item{(a)} $Y^{p,q,r}(w)=Y^{p,q,r}(\alpha)=m-1$ 
			\item{(b)} $Y^{p,q,r}[L] = [1,p]\sqcup [m-q,m-1]\sqcup [m+p,m+p+r-1]$.
		\end{description}  
		Note that $Y^{p,r,q}$ uses exactly $p+q+r=m+1$ colors for the leaves. Fix an initial $\kappa$-coloring $Y_{0}\in \mathcal{Y}^{p,q,r}$. By a similar reasoning as in \textbf{\eqref{I}}, it is easy to check that $Y_{3m+q}-m+q\in \mathcal{Y}^{q,r,p}$. By iterating this observation, we get $Y_{6m+q+r}-q-r\in \mathcal{Y}^{r,p,q}$ and $Y_{9m+p+q+r}-m+p+q+r=Y_{5\kappa+1}+1\in \mathcal{Y}^{p,q,r}$. In fact, it is easy to see that $Y_{5\kappa+1}+1=Y_{0}$ on $V$. As in the proof of \textbf{\eqref{I}}, this shows that $X_{t}$ does not synchronize and its period divides $5\kappa^{2}+\kappa$.
	\end{description}	
\end{proof}

For \textbf{(ii)}, we remark that all 8-colorings on any star synchronize.

\section{Conditional proofs of Theorem \ref{thm:tree_time} for $\kappa=6$ and Theorem \ref{thm:blinkingtree}}

\label{sec:pf_blinkingtree_6}

\subsection{Conditional proof of Theorem \ref{thm:tree_time} for $\kappa=6$}

Our proof of Theorem \ref{thm:tree_time} for $\kappa=6$ will reduce the following `reduction lemma' of the 6-color FCA on trees. We used analogous results for the proof of Theorem \ref{thm:tree_time} for $\kappa\in \{3,4\}$ (see Lemmas \ref{lem:4FCA_branch} and \ref{lem:5FCA_branch_root}), where one could reduce the dynamics on a smaller subtree after a constant time. A crucial difference for the $\kappa=6$ case is that, the time until reduction is not constant anymore but it depends linearly on the diameter of the tree. 

\begin{lemma}[6-color FCA local dynamics on a branch]
	\label{lem:6FCA_branch}
	Let $(X_{t})_{t\ge 0}$ be a 6-color FCA trajectory on a tree $T=(V,E)$ of diameter $d$. Suppose $T$ has  a $r$-branch $B$ for $r\ge 1$ with center $v$, root $w$, and leaves $l_{1},\dots,l_{r}$ (see Figure \ref{fig:branch_pic}\textbf{b}). Suppose $v$ blinks at some time $t_{0}\ge 3$. Then the dynamic $(X_{t})_{t\ge t_{0}+ Cd}$ restricts on $T-B$, where $C>0$ is an absolute constant. 
\end{lemma}

One can conclude Theorem \ref{thm:blinkingtree} for $\kappa=6$ easily from Lemma \ref{lem:6FCA_branch}. 

\begin{proof}[\textbf{Proof of Theorem \ref{thm:blinkingtree} for $\kappa=6$}]
	Suppose the underlying tree $T$ has diameter $d$ and maximum degree $\le 5$. By Lemma \ref{lem:blinking_degree}, every node in $T$ blinks at least once in every 24 iterations. Hence using Lemma \ref{lem:6FCA_branch}, the dynamic on $T$ after time $24+Cd$ reduces to the subtree $T'$ obtained from $T$ by deleting all leaves in $T$. Then one can conclude by induction on the diameter as in the proof of Theorem \ref{thm:tree_time} for $\kappa\in \{3,4\}$. 
\end{proof}

Our proof of  Lemma \ref{lem:6FCA_branch} is significantly more technical than the $\kappa\in\{4,5\}$ case. We will need to develop a number of new techniques so that we can use a recursive argument on the dynamics of the entire tree in order to reduce one step. The key ingredients are stated in the following subsection and we will justify Lemma \ref{lem:6FCA_branch} at the end of this section.

\subsection{Conditional proof of Theorem  \ref{thm:blinkingtree}}

Trees have nice recursive properties when viewed as \textit{rooted trees}. A rooted tree $T=(V,E)$ is a tree with a designated vertex $\mathtt{r}$ called the \textit{root}. For each $v \in V$, let $\vec{P}_{v}$ be the unique path from $\mathtt{r}$ to $v$, and let $v^{-}$ be the unique neighbor of $v$ in $\vec{P}_{v}$. For distinct vertices $u,v$ in $T$, we say $u$ is the \textit{parent} of $v$ and $v$ is a \textit{child} of $u$ if $u=v^{-}$. Write $u \le v$ if $u \in V(\vec{P}_{\tau})$. We say $u$ is an \textit{ancestor} of $v$ and $v$ is a \textit{descendant} of $u$ if $u\le v$. For each $v\in V$, we define the \textit{descendent subtree} $T_{v}$ to be the subtree of $T$ consisting of all descendants of $v$. The \textit{depth} of $T_{v}$, denoted by $\depth(v)$, is the maximum level of vertices in $T_{v}$.  We say a descendent subtree $T_{v}$ is a \textit{terminal branch} if it is a branch and either $T_{v}=T$, or $v^{-}\in V$ and $T_{u}$ is either a singleton or a branch for all children $u$ of $v^{-}$. Note that any rooted tree with depth $\ge 2$ has at least one terminal branch.

Now for each $\kappa\ge 3$, call a pair $(T,X_{0})$ of a finite tree and a $\kappa$-coloring on it a \textit{minimal counterexample} if (1) $X_{0}$ does not synchronize, (2) every vertex in $T$ blinks infinitely often in the dynamic $(X_{t})_{t\ge 0}$, and (3) $|V(T)|\le |V(T')|$ if $(T',X_{0}')$ is another pair with satisfying (1) and (2). We may assume without loss of generality that $(X_{t})_{t\ge 0}$ is periodic, by choosing $X_{0}$ from the periodic limit cycle. Note that by the minimality of $T$ and Lemma \ref{lem:branchwidth}, every branch in a minimal counterexample must have branch width $\ge \kappa/2-1$ for all times. This enforces very specific local dynamics on branches which easily led to a contradiction in the case of $\kappa\in \{3,4,5\}$. 

Let $(T,X_{0})$ be a minimal counterexample for $\kappa=6$, and we fix this notation hereafter throughout later sections. Let $T_{v}$ be a proper descendant subtree in $T$. We say $T_{v}$ is \textit{open} if the induced dynamics on $T_{v}$ requires $v$ to be pulled by its parent $v^{-}$. In particular, the whole tree $T=T_{\root}$ cannot be opened. If $T_{v}$ is open, then the minimality will force particular induced local dynamics on its root $v^{-}$. To represent induced local dynamics on a single vertex concisely,  we introduce the following notion. For each vertex $v$ in $T$, let $\mathtt{bt}_{i}(v)$ be the time of $i$th blink of $v$ in the orbit, for each $i\ge 1$. The $i$th \textit{blinking gap} is given by $\mathtt{g}_{i}(v)=\mathtt{bt}_{i+1}(v)-\mathtt{bt}_{i}(v)$. The \textit{blinking sequence} of $v$ is defined by the sequence $(\mathtt{g}_{i}(v))_{i\ge 1}$ of blinking gaps of $v$. Note that since we are in a periodic orbit, the blinking sequence of any $v\in V$ repeats a finite sequence of blinking gaps. There will be mainly four different types of enforced dynamics on $T_{v}$ as we define below:

\begin{customdef}{4.3}\label{blinkingseq_def}
	Let $(T,X_{0})$ as before and let $T_{v}$ an open descendant subtree of $T$.	
	\begin{description}
		\item{\eqref{I}} We say $T_{v}$ is of \textit{type (a)} iff $\mathtt{g}_{i}(v)\equiv 12$ and $2\in \{X_{t+4}(v^{-}),X_{t+5}(v^{-})\}$ whenever $X_{t}(v)=2$; 
		\item{(ii)} We say $T_{v}$ is of \textit{type (b)} iff $\mathtt{g}_{i}(v)$ alternates  9 and 7, and $2\in \{X_{t+2}(v^{-}),X_{t+4}(v^{-})\}$ whenever $X_{t}(v)=X_{t+9}(v)=2$.
		\item{(iii)} We say $T_{v}$ is \textit{fractal of type 10/9} iff $\mathtt{g}_{i}(v)$  alternates 10 and 9, and $2\in \{X_{t+1}(v^{-}),X_{t+4}(v^{-})\}$ whenever $X_{t}(v)=X_{t+10}(v)=2$; 
		\item{(iv)} We say $T_{v}$ is \textit{fractal of type 11/8} iff $\mathtt{g}_{i}(v)$ alternates  11 and 8, and $2\in \{X_{t+2}(v^{-}),X_{t+4}(v^{-})\}$ whenever $X_{t}(v)=X_{t+11}(v)=2$.
	\end{description}	  
	We say $T_{v}$ is \textit{fractal} if it's a fractal of either type.  
\end{customdef}

Below we give a more direct characterization of type (a), (b), and fractal subtrees in terms of the induced dynamics on $v$ and its parent $v^{-}$.  Suppose a descendant subtree $T_{v}$ is of type (a). From the definition, it is easy to see that the local dynamics on $v$ and $v^{-}$ are given by concatenating the following four sequences:
\begin{align}\label{P}
	\centering{
		\begin{matrix}
			\text{$v$} & \textbf{2}	& 3 & - & - & - & - & - & - & - & 5 & 0 & 1 & \textbf{2} \\
			\text{$v^{-}$} & 5 & 5 & 0 & 1 & \mathbf{2} & 3 & - & - & - & - & - & - & -
		\end{matrix}\tag{P}
	}
\end{align}

\begin{equation}\label{Q}
	\centering{
		\begin{matrix}
			\text{$v$} & \textbf{2}	& 3 & - & - & - & - & - & - & - & 5 & 0 & 1 & \textbf{2} \\
			\text{$v^{-}$} & 4 & 4 & 5 & 0 & 1 & \mathbf{2} & 3 & - & - & - & - & -  & -
		\end{matrix}\tag{Q}
	}
\end{equation}

\begin{equation}\label{R}
	\centering{
		\begin{matrix}
			\text{$v$} & \textbf{2}	& 3 & - & - & - & - & - & - & - & 5 & 0 & 1 & \textbf{2} \\
			\text{$v^{-}$} & 5 & 5 & 5 & 0 & 1 & \mathbf{2} & 3 & - &- & - & - & -  & -
		\end{matrix}\tag{R}
	}
\end{equation}

\begin{equation}\label{S}
	\centering{
		\begin{matrix}
			\text{$v$} & \textbf{2}	& 3 & - & - & - & - & - & - & - & 5 & 0 & 1 & \textbf{2} \\
			\text{$v^{-}$} & 5 & 5 & 0 & 1 & \mathbf{2} & 3 & 4 & 5 & 0 & 1 & \mathbf{2} & 3 & 4
		\end{matrix}\tag{S}
	}
\end{equation}
where the time goes from left to right and none of $-$'s are 2. Each of the above sequences describes dynamics on $v$ and $v^{-}$ for 12 iterations, and the sequence \eqref{P}\eqref{P} obtained by concatenating the sequence \eqref{P} twice describes 24 iterations, for instance. The four possibilities came from considering possible instances of $v^{-}$ blinking after its first blink during each blinking gap 12 of $v$. 

Similarly, if $T_{v}$ is of \textit{type (b)}, then the local dynamics on $v$ and $v^{-}$ are given by concatenating the following six sequences:
\begin{equation}\label{I}
	\begin{matrix}
		\text{$v$} & \textbf{2} & 3 & - & - & - & - & 5 &  0 &  1 &  \textbf{2} &  3 &  - &  - &  5 & 0 & 1  & \mathbf{2}\\
		\text{$v^{-}$} & 0 & 1 & \mathbf{2} & 3  & - & - &- &-& -&-& - &-& - &  - &-& - & -
	\end{matrix}\tag{I}
\end{equation}
\begin{equation}\label{J}
	\begin{matrix}
		\text{$v$} & \textbf{2} & 3 & - & - & - & - & 5 &  0 &  1 &  \textbf{2} &  3 &  - &  - &  5 & 0 & 1  & \mathbf{2}\\
		\text{$v^{-}$} & 5 & 5 & 0  & 1 & \mathbf{2} & 3 & - & - & - & -& -  & - & - & - & - & - & -
	\end{matrix}\tag{J}
\end{equation}

\begin{equation}\label{X}
	\begin{matrix}
		\text{$v$} & \textbf{2} & 3 & - & - & - & - & 5 &  0 &  1 &  \textbf{2} &  3 &  - &  - &  5 & 0 & 1  & \mathbf{2}\\
		\text{$v^{-}$} & 0&1& \textbf{2} &3 &4&5&0&1&\textbf{2}&3&-&-&-&-&-&-&-
	\end{matrix}\tag{X}
\end{equation}
\begin{equation}\label{Y}
	\begin{matrix}
		\text{$v$} & \textbf{2} & 3 & - & - & - & - & 5 &  0 &  1 &  \textbf{2} &  3 &  - &  - &  5 & 0 & 1  & \mathbf{2}\\
		\text{$v^{-}$} & 0&1& \textbf{2} &3 &- & -&5 & 0&1&\textbf{2}&3&-&-&-&-&-&-
	\end{matrix}\tag{Y}
\end{equation}
\begin{equation}\label{Z}
	\begin{matrix}
		\text{$v$} & \textbf{2} & 3& - & - & - & - & 5 &  0 &  1 &  \textbf{2} &  3 &  - &  - &  5 & 0 & 1  & \mathbf{2}\\
		\text{$v^{-}$} & 0&1& \textbf{2} &3 &- & - &- & 5 & 0 &1&\textbf{2}&3&-&-&-&-&-
	\end{matrix}\tag{Z}
\end{equation}\begin{equation}\label{W}
	\begin{matrix}
		\text{$v$} & \textbf{2} & 3 & - & - & - & - & 5 &  0 &  1 &  \textbf{2} &  3 &  - &  - &  5 & 0 & 1  & \mathbf{2}\\
		\text{$v^{-}$} &5 & 5 & 0&1& \textbf{2} &3 &- & - &5 & 0&1&\textbf{2}&3&-&-&-&-
	\end{matrix}\tag{W}
\end{equation}
where none of $-$'s are 2, as before. 

Finally, the same holds for the following two sequences 
\begin{equation}\label{F1}
	\small{\begin{matrix}
			v &\textbf{2}&3&3&-&-&-&-&5&0&1&\textbf{2}&3&-&-&-&-&5&0&1&\mathbf{2}\\
			v^{-}&1& \textbf{2}&3&-&-&-&-&-&a_{1}&a_{2}&a_{3}&a_{4}&a_{5}&-&a_{7}&a_{8}&-&-&- &-
	\end{matrix}}\tag{F1}
\end{equation}
\begin{equation}\label{F2}
	\small{\begin{matrix}
			v &\textbf{2}&3&-&-&-&-&-&5&0&1&\textbf{2}&3&-&-&-&-&5&0&1&\mathbf{2}\\
			v^{-} &5& 5&0&1&\textbf{2}&3&-&-&-&-&b_{3}&b_{4}&b_{5}&-&b_{7}&b_{8}&-&-&-&- 
	\end{matrix}}\tag{F2}
\end{equation}
when $T_{v}$ is fractal of type 10/9, and with the following two sequences for $T_{v}$ fractal of type 11/8:
\begin{equation}\label{F3}
	\small{\begin{matrix}
			v&\textbf{2}&3&-&-&-&-&-&-&5&0&1&\textbf{2}&3&-&-&-&5&0&1&\mathbf{2}\\
			v^{-}&0& 1& \textbf{2}&3&-&-&-&-&-&c_{2}&c_{3}&c_{4}&c_{5}&c_{6}&-&c_{8}&-&-&- &-
	\end{matrix}}\tag{F3}
\end{equation}
\begin{equation}\label{F4}
	\small{\begin{matrix}
			v&\textbf{2}&3&-&-&-&-&-&-&5&0&1&\textbf{2}&3&-&-&-&5&0&1&\mathbf{2}\\
			v^{-}&5& 5&0&1&\textbf{2}&3&-&-&-&-&d_{3}&d_{4}&d_{5}&d_{6}&-&d_{8}&-&-&- &-
	\end{matrix}}\tag{F4}
\end{equation}

\noindent As before, none of $-$'s are 2 in any of the sequences above, but in other instances, $a_{i}$'s, $b_{i}$'s, $c_{i}$'s, and $d_{i}$'s could be 2. We could specify all sequences when $v^{-}$ could blink among those instances, but there would be too many cases of doing so. 

Now we outline the proof of Theorem \ref{thm:blinkingtree}. In a nutshell, we show that every proper descendant subtree $T_{v}$ of depth $\ge 1$ is fractal. In particular, every component in $T_{\root}-\root$ will be either a singleton or fractal. A recursive property of fractal subtrees would then yield that the whole tree is fractal, and in particular, open. This contradiction shows that a minimal counterexample for $\kappa=6$ does not exist. To give more detail, we first show by using Lemma \ref{lem:branchwidth}, that every branch must be open and of type (a) or (b), or fractal of type 10/9. Furthermore, we will show that if $T_{v}$ is a terminal branch, then it cannot be of type (a) or (b), as stated in the following lemma: 
\begin{lemma}\label{lem:terminal_is_fractal}
	Let $(T,X_{0})$ be as before. Then every terminal branch of $T$ is fractal. 
\end{lemma}

Next, the induction step is based on the recursive property of fractal branches stated in the following lemma: 

\begin{lemma}\label{lem:recursion_fractalbranch}
	Let $(T,X_{0})$ be as before. Let $w\in V$ and suppose that each connected component of $T_{w}-w$ is either a singleton or fractal. Then $T_{w}$ is open and fractal. In particular, $w^{-}\in V(T)$.  
\end{lemma}

These two lemmas easily imply the main theorem. 

\begin{proof}[\textbf{Proof of Theorem \ref{thm:blinkingtree} }] It suffices to show the "if" part. On the contrary, suppose there exists a minimal counterexample $(T,X_{0})$ for $\kappa=6$. Choose a vertex $\root$ and view $T$ as being rooted at $\root$. By minimality, the depth of $T=T_{\root}$ is at least $1$. If $\depth(\root)=1$, then $T=T_{\root}$ is a terminal branch, which is open by Proposition \ref{branchorbit_detail}, a contradiction. Hence we may assume that $\depth(\root)\ge 2$. 
	
	It suffices to show that for every non-leaf and non-root vertex $v$, $T_{v}$ is fractal. Indeed, this would yield that $T_{\root}-\root$ is a disjoint union of leaves and fractal subtrees, but by Lemma \ref{lem:recursion_fractalbranch} $T_{\root}$ must be fractal, a contradiction. We proceed by induction on $\depth(v^{-})\ge 2$. For the base step, note that $\depth(v^{-})=2$ means that $T_{v^{-}}-v^{-}$ is a disjoint union of leaves and terminal branches. Since by Lemma \ref{lem:terminal_is_fractal} terminal branches are fractal, Lemma \ref{lem:recursion_fractalbranch} gives that $T_{v}$ is fractal. The induction step follows similarly. If $\depth(v^{-})=d\ge 3$, then $T_{v}-v$ is a disjoint union of leaves and depth $< d-1$ descendant subtrees. By induction hypothesis and Lemma \ref{lem:recursion_fractalbranch}, $T_{v}$ must be fractal. This shows the assertion. 
\end{proof}

\subsection{Conditional proof of Lemma \ref{lem:6FCA_branch}}

In this subsection, we derive Lemma \ref{lem:6FCA_branch} from Lemmas \ref{lem:terminal_is_fractal} and \ref{lem:recursion_fractalbranch}. 

\begin{proof}[\textbf{Proof of Lemma \ref{lem:6FCA_branch}}] 
	In this proof, we will assume Lemmas \ref{lem:terminal_is_fractal} and \ref{lem:recursion_fractalbranch} hold. While these lemmas were stated in terms of a minimal counterexample that does not consider time complexity, one can consider a quantitative version of the minimal counterexample as follows. Fix a tree $T$ of diameter $d$ and maximum degree $\le 5$. Fix a large constant $C_{0}$. Call a pair $(T,X_{0})$ a (quantitative) minimal counterexample if (1) $T$ has maximum degree $\le 5$; and (2) the 6-color FCA dynamic $(X_{t})_{t\ge 0}$ on $T$ does not restrict to any smaller subtree within $C_{0}d$ iterations, where $d=\textup{diam}(T)$; (3) if $(T',X_{0}')$ is another pair satisfying (1) and (2), then $|V(T')|\le |V(T)|$. 
	
	We can also adopt Def. \ref{blinkingseq_def} as well as the definition of fractal branches so that one requires the properties of blinking gap sequences to only hold during the time interval $[1,C_{0}d^{2}]$. Namely, we say an open descendant tree $T_{v}$ of $T$ a \textit{fractal of type (a) or (b) from time $t_{0}$} if the blinking sequence of $v^{-}$ satisfy the conditions in \eqref{I}-(iv) in Def. \ref{blinkingseq_def} for times $t\in [1,C_{0}d]$. Then one can verify that Lemma \ref{lem:terminal_is_fractal} still holds for the quantitative notion of fractal branches by the proof of Lemma \ref{lem:terminal_is_fractal} (see Section \ref{sec:proof_lem_terminal_fractal}) after a very minor modification. Namely, every terminal branch of $T$ is fractal after some time $t_{0}$ independent of $T$. 
	
	Next, one can also extend Lemma \ref{lem:recursion_fractalbranch} quantitatively. Namely, if $T_{w}-w$ is either a singleton or a fractal from time $t_{1}$, then $T_{w}$ is open and fractal from time $t_{1}+C_{1}$, where $C_{1}>0$ is an absolute constant. The constant $C_{1}$ needs to be large enough so that we get the same contradiction as in the proof of Lemma \ref{lem:recursion_fractalbranch} (see Section \ref{sec:proof_lemma_recursion_fractal}) for the blinking sequence of $w$ during time interval $[t_{1}, t_{1}+C_{1}]$. Here $C_{1}$ does not depend on the tree $T$ as the total number of possible induced blinking sequences of $w$ does not depend on $T$. 
	
	Now we show how to conclude. Suppose there exists a quantitative counter example $(T,X_{0})$ and $T$ has diameter $d$. By the quantitative version of Lemma \ref{lem:terminal_is_fractal}, all of its terminal branches are fractal after time $t_{0}$. By the quantitative version of Lemma \ref{lem:recursion_fractalbranch}, it follows that the whole tree $T=T_{\mathtt{r}}$ is open and fractal after time $t_{0} + d C_{2}$, where $C_{2}$ is a possibly larger constant than $C_{1}$. Thus if $C_{0}$ is large enough so that $t_{0}+dC_{2}<C_{0} d$, this is a contradiction. 
\end{proof}

\section{Analysis of enforced local orbits on branches}
\label{sec:6FCA_local}

Throughout this section, $(T,X_{0})$ is a minimal counterexample for $\kappa=6$. For each descendant subtree $T_{v}$, the minimality forces a particular local dynamic on $T_{v}$, which may give some constraints on its root $v^{-}$. In this section, we analyze such enforced local orbits on $T_{v}$ and see how they restrict the dynamics on $v^{-}$ when $T_{v}$ is either a leaf, branch, or a fractal branch. Furthermore, we investigate a possible ensemble of such constraint on the local dynamics of $v^{-}$ when it has multiple descendant subtrees rooted at itself. A conceptual background is a classic technique in dynamical systems literature called the Poincar\'e return map, which is to look at transitions between snapshots of system configuration where a particular vertex takes a particular state. We adapt this concept in a local setting: we consider all possible local configurations on a descendant subtree $T_{v}$ in which $v$ blinks. Since we are assuming that $v$ blinks infinitely often in the dynamic, the global periodic orbit $(X_{t})_{t\ge 0}$ must induce a periodic orbit on such special local configurations, together with constraints on the local dynamics on $v^{-}$.

\begin{figure}[H]
	\centering
	\includegraphics[width=0.3 \linewidth]{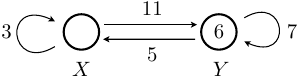}
	\caption{A digraph that can generate finite sequence $3,11,6,7,6,5,3$.}
	\label{fig:digraph_generating_seq_ex1}
\end{figure}

We will rely heavily on diagrammatic analysis to study possible blinking sequences of $v^{-}$ and their ensemble. We shall represent local dynamics on $T_{v}$ often as a weighted digraph, in which edge weights represent blinking gaps of $v$ and nodes could be snapshots of local configurations or a finite sequence of local dynamics. Let us first introduce some terminologies. Let $D=(V,\overline{E})$ be a digraph with vertex and edge weights $\omega:V\sqcup E\rightarrow \mathbb{N}\cup\{0\}$. We say a sequence $(a_{n})$ of positive integers is \textit{generated} by $D$ if there exists a directed walk $P=v_{1}e_{1}v_{2}e_{2},\cdots $ in $D$ such that $(a_{n})$ can be obtained from the sequence $\omega(v_{1}),\omega(e_{1}),\omega(v_{2}),\cdots$ by dropping the  zero terms. For example, consider a digraph $D$ with vertex set $V=\{X,Y\}$ and edge set $E=\{(XX), (XY), (YX), (YY)\}$ with weights given as in Figure \ref{fig:digraph_generating_seq_ex1}. 
Whenever a vertex has a weight of 0, we shall omit the weight in the diagram. Notice that the directed walk $P=X,(XX),X,(XY),Y,(YY),Y,(YX),X,(XX)$ gives the sequence of weights $0,3,11,6,7,6,5,0,3$. By dropping out the zero terms, we see that the given digraph can generate the sequence $3,11,6,7,6,5,3$.

\begin{prop}\label{prop_leaf_dynamics}
	Let $(T,X_{0})$ as before. Let $v$ be a vertex in $T$ with a leaf neighbor $u$. Then the blinking sequence of $v$ is given by $(a_{i}+6k_{i})_{i\ge 1}$ where $(a_{i})_{i\ge 1}$ is generated by the digraph in Figure \ref{localdynamicswithleaf} and $(k_{i})_{i}$ is some sequence of non-negative integers which depend on the dynamics.  
	\begin{figure}[H]
		\centering
		\includegraphics[width=0.5 \linewidth]{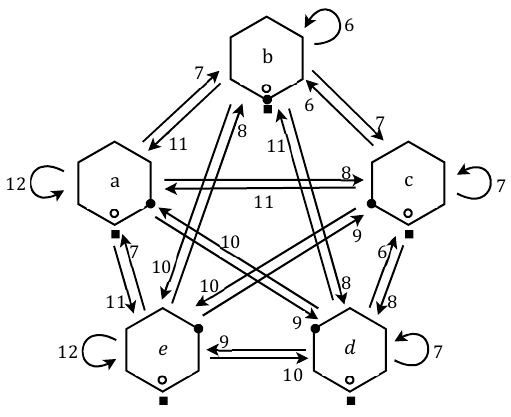}
		\caption{A weighted digraph on five possible non-opposite local configurations on the 1-star $v+u$. $\blacksquare=$phase of center $v$, $\bullet=$phase of leaf $u$, and $\circ = $phase of activator. The blinking gap of $v$ corresponding to each transition is given by $a+6k$ where $a$ is the edge weight and $k$ is some nonnegative integer depending on the structure and dynamics of $G$.}
		\label{localdynamicswithleaf}
	\end{figure}         
\end{prop}

\begin{proof}
	There are five local configurations on the 1-star $v+u$ with center $v$ where $v$ blinks such that $u$ is not opposite to $v$. The above digraph shows every possible transition between such five non-opposite local configurations. For the edge weights, note that since $v$ is the only neighbor of the leaf $u$, once $v$ blinks, $u$ maintains its phase until the next blink of $v$. This determines the blinking gap of $v$ during each transition in the above digraph modulo $6$. 
	
	For example, consider the transition $d\rightarrow e$ in Figure \ref{localdynamicswithleaf}, which is shown in Figure \ref{leaf_transition_ex}.
	\begin{figure}[H]
		\centering
		\includegraphics[width = 0.45 \linewidth]{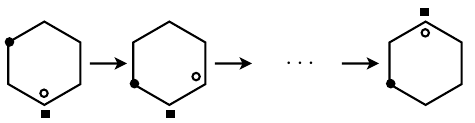}
		\caption{The transition $d\rightarrow e$ in Figure \ref{localdynamicswithleaf}. The blinking gap must be $9+6k$ for some nonnegative integer $k$.}
		\label{leaf_transition_ex}
	\end{figure}         	
	Since the center pulls the leaf initially, the phase of the leaf moves one step clockwise after the first iteration. Now the leaf does not move until the next blink of the center, so to get the bottom left local configuration in Figure 16, the center must be at the top of the hexagon by the time it blinks again for the first time. Hence looking that the initial and terminal phases of the activator, we conclude that the blinking gap of $v$ during this transition is $9$ modulo $6$. Other edge weights are determined in a similar way. This shows the assertion. 
\end{proof}

Next, we analyze forced local dynamics on branches. By Lemma \ref{lem:branchwidth}, if the $T$ has a $r$-branch $B$ and if the dynamic on $T$ does not restrict on $T-B$, then the branch width of $B$ must maintain a value $\ge 2$. This will enforce the root $v$ of the branch $B$  to blink at particular instances.  The following proposition gives how the blinking sequence of a vertex $v$ is restricted if it has multiple leaves, which includes the case when $v$ is a center of a branch in $T$. Its proof is given at the end of this section. 

\begin{prop}[Local dynamics of 6-color FCA on a branch] \label{kstar_orbit}
	Let $(T,X_{0})$ be as before. Suppose $T$ has a $k$-star $S$ for $k\ge 2$ with center $v$. Then we have the following: 
	\begin{description}[noitemsep]
		\item{\eqref{I}} The induced local dynamics on $S$ is given by one of the four digraphs in Figure \ref{starorbit_fig}. In particular, the blinking sequence of $v$ is given by $(a_{i}+6k_{i})_{i\ge 1}$ where $(a_{i})_{i\ge 1}$ is generated by the digraph in Figure \ref{starorbit_fig} and $(k_{i})_{i}$ is some sequence of non-negative integers which depend on the dynamics.  
		
		\item{(ii)} If $S=T_{v}$ is a branch, then the induced local dynamics on $T_{v}$ only uses the five shaded local configuration in Figure \ref{starorbit_fig}. 
		
		\item{(iii)} If $S=T_{v}$ has local dynamics given by Figure \ref{starorbit_fig} (a), (b), or (c), then it is open and of type (a), (b), or fractal of type 10/9, respectively. 
	\end{description}
	\begin{figure}[H]
		\centering
		\includegraphics[width=0.95 \linewidth]{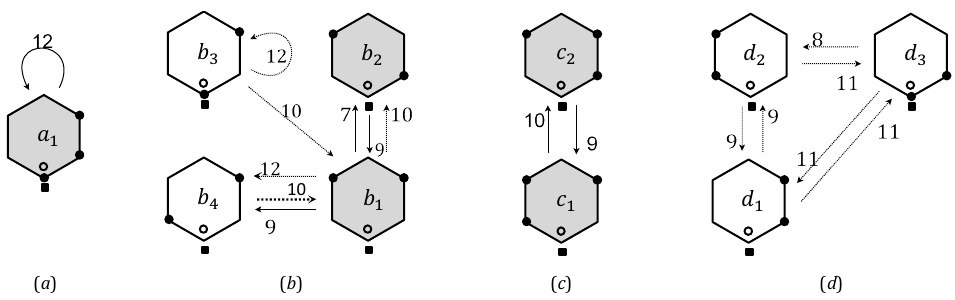}
		\caption{Possible transitions between local configurations on a star with $\ge 2$ leaves. The edge weights indicate the corresponding minimal blinking gaps of the center. $\blacksquare=$phase of center $v$, $\bullet=$phase of leaves, and $\circ = $phase of activator. The dotted transitions are impossible if $B=T_{v}$ is the center of a branch.}
		\label{starorbit_fig}
	\end{figure}         
\end{prop}

Next, we investigate how the three types of closed orbits on a branch restrict the blinking sequence of its root. Let $(a_{n}), (b_{n})$ be two sequences of real numbers. We say the sequence $(b_{n})$ \textit{refines} $(a_{n})$ and $(b_{n})$ is a \textit{coarsening} of $(a_{n})$ if there exists an increasing sequence $(d_{n})$ of natural numbers such that 
\begin{equation*}
	a_{n} = \sum_{d_{n}\le k < d_{n+1}} b_{k}.
\end{equation*}
For instance, the sequence $1,2,3,4\cdots$ refines $3,7,11,15,\cdots$ since $3=1+2$, $7=3+4$, $11=5+6$, and so on.

\begin{prop}\label{branchorbit_detail} 
	Let $(T,X_{0})$ be as before. Let $T_{v}$ be a branch in $T$ with $v^{-}\in T$. Then we have the following: 
	\begin{description}[noitemsep]
		\item[(i)] If $T_{v}$ is of type (a), then the blinking sequence of $v^{-}$ is generated by the digraph (A) in Figure \ref{branchorbit_root_fig}.
		\item[(ii)] If $T_{v}$ is of type (b) then the blinking sequence of $v^{-}$ refines a sequence generated by digraph(B1) in Figure \ref{branchorbit_root_fig}.
		\item[(iii)] In case of \textbf{(ii)}, the blinking sequence $(\mathtt{g}_{i})$ of $w$ is refined by some sequence $(b_{m})$ generated by diagram (B2) in Figure \ref{branchorbit_root_fig}. Furthermore, $(\mathtt{g}_{i})_{i\ge 1}$ can be obtained from $(b_{m})$ by merging some consecutive terms $b_{m},b_{m+1}$ into $b_{m}+b_{m+1}$, where $b_{m}$ is a vertex weight and $b_{m+1}$ is following edge weight. 
		\item[(iv)] In case of \textbf{(iii)}, the sequence $(b_{m})$ cannot be generated by a closed walk on (B2) which only uses nodes $Y$ or $Z$. 
	\end{description}
	\vspace{-0.5cm}
	\begin{figure}[H]
		\centering
		\includegraphics[width=1 \linewidth]{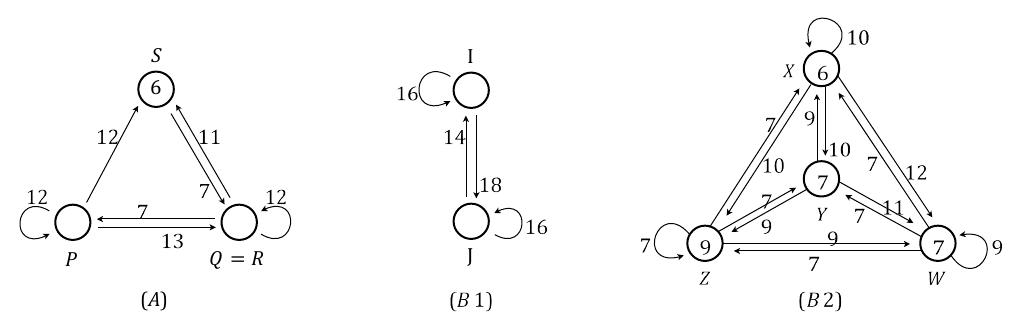}
		\caption{(A) a digraph generating the blinking sequence of the root $v^{-}$ of type (a) branches; (B1) a digraph generating a coarsening of blinking sequence of $v^{-}$ of type (b) branches; (B2) a digraph generating a refinement of blinking sequence of $v^{-}$ of type (b) branches. }
		\label{branchorbit_root_fig}
	\end{figure}         
\end{prop}

\begin{proof}
	The proof follows mostly from definitions. Let $T_{v}$ be of type (a). Then concatenating sequences $\eqref{P}$-$\eqref{S}$ gives a complete description of the blinking sequence of $v^{-}$. For instance, if string $\eqref{P}\eqref{P}\eqref{Q}$ is used in the local dynamics, then $v^{-}$ blinks exactly once in the sequence $\eqref{P}$, and blinks after 12 iterations again in $\eqref{P}$, and then its next blink in $\eqref{Q}$ takes 13 iterations. In digraph $(A)$ in Figure \ref{branchorbit_root_fig}, this is represented as going through the loop at node \eqref{P} twice and then using the edge (PQ). Note that diagram (A) lacks loop at node $S$ and edges from $S$ to $P$ or $Q$, since those sequences cannot be concatenated in such order; the color of $v^{-}$ at the end and beginning does not match. To explain the use of node weight on $S$ in diagram (A) in Figure \ref{branchorbit_root_fig}, consider the string of sequences $\eqref{P}\eqref{S}\eqref{Q}$. After the blink within sequence $\eqref{P}$, $v^{-}$ blinks for the first time in the sequence $\eqref{S}$ after 12 iterations, and then again for the second time after 6 iterations within sequence $\eqref{S}$. Then it takes 7 iterations to blink again within sequence $\eqref{Q}$. In terms of diagrams, we walk through the edge weight 12 of $(PS)$, then node weight 6 of S, and then edge weight 7 of $(SQ)$.  This shows \textbf{\eqref{I}}. 
	
	For type (b) branches, observe that if $v^{-}$ blinks as sparse as possible in the dynamics, then it would only use the ``long periodic'' sequences \eqref{I} and \eqref{J}, in which case its blinking sequence is generated by diagram (B1) in Figure \ref{branchorbit_root_fig}. On the other hand, if $v^{-}$ blinks as often as possible, only those four ``short periodic'' sequences \eqref{X}-\eqref{W} would be used and its blinking sequence is generated by Figure \ref{branchorbit_root_fig}. In general, the actual local dynamics on $v$ and $v'$ could use all combinations, which means that $v^{-}$ could blink within long periodic sequences \eqref{I} and \eqref{J} or could skip the second blinks in short periodic sequences $\eqref{X}$-$\eqref{W}$. Thus the actual blinking sequence of $v^{-}$ refines a sequence generated by diagram (B1), but could be coarser than a sequence generated by diagram (B2); skipping second blinks within short periodic sequences corresponds to merging node weights with the following edge weights in diagram (B2). For example, the string $\eqref{X}\eqref{J}\eqref{Z}$ is represented on diagram (B2) by the directed walk $X$,$(XW)$,$W+(WZ)$,$Z$, which generates the sequence 6, 12, (7+7), 9. This shows \textbf{(ii)} and \textbf{(iii)}.
	
	Lastly, suppose $v^{-}$ only uses sequences \eqref{Y} and \eqref{Z}. Note that the center $v$ does not pull $v^{-}$ in those sequences, since $v^{-}$ has colors $\le 2$ whenever the center has color $2$. Hence if the induced local orbit on $v$ and $v^{-}$ is given by an infinite subsequence of \eqref{Y} and \eqref{Z} only, then the dynamics restricts on $T-T_{v}$, a contradiction.  This shows \textbf{(iv)}. 
\end{proof}

\begin{prop}\label{prop_fb_1}
	Let $(T,X_{0})$ be as before. Let $T_{w}$ be a fractal branch in $T$ with $w^{-}\in V(T)$. Then the blinking sequence of $w^{-}$ refines a sequence generated by digraph (F10-9) or (F11-8) corresponding to the type of $T_{w}$. 
	\begin{figure}[H]
		\centering
		\includegraphics[width=0.8 \linewidth]{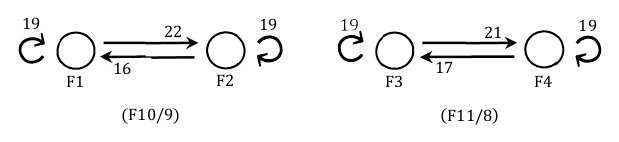}
		\caption{ The blinking sequence of $w^{-}$ refines a sequence generated by (F10/9) or (F11/8) depending on the type of $T_{w}$.
		}
		\label{fractalbranch_coarsening}
	\end{figure}         
\end{prop}

\begin{proof}
	Let $T_{w}$ be of type $10/9$. According to the definition, the dynamic on $w$ and $w^{-}$ during consecutive blinking gaps of 10 and 9 of $w$ is given by concatenations of the four sequences \eqref{F1}-\eqref{F4} in Section 3. Note that $w^{-}$ may or may not blink at some of the $a_{i}$'s, $b_{i}$'s, $c_{i}$'s, or $d_{i}$'s. But if one ignores such blinks within each sequence of 19 iterations (F$i$)'s, the blinking gap of $w^{-}$ must be generated by the digraphs in Figure \ref{fractalbranch_coarsening}. For instance, if sequence \eqref{F1} is followed by \eqref{F2}, then it takes 22 iterations for $w^{-}$ to blink at the beginning of each (F$i$)'s. Thus the actual blinking sequence of $w^{-}$ must refine a sequence generated by digraphs in Figure \ref{fractalbranch_coarsening} depending on the type of $T_{w}$.  
\end{proof}

\begin{lemma}\label{exclusive_branches}
	Let $(T,X_{0})$ as before. Suppose $w\in V(T)$ such that each component of $T_{w}-w$ is either a singleton, branch, or fractal. Then branches of type (a) or (b) or fractal of either type rooted at $w$ are mutually exclusive. 
\end{lemma}

\begin{proof}
	Suppose there are both types (a) and (b) branches rooted at $w$. Then by Proposition \ref{branchorbit_detail}, the blinking sequence of $w$ is generated by the diagram (A) and must refine a sequence generated by diagram (B1) in Figure \ref{branchorbit_root_fig}. It is easy to see that the sum of the edge and vertex weights in any directed walk in diagram (A) cannot be 14 or 16. This means that any sequence generated by (A) cannot refine a sequence that contains a term of 14 or 16. But any sequence generated by (a directed closed walk in) diagram (b1) must contain a term of 14 or 16. Hence this is impossible. 
	
	Next, suppose there is one branch $B$ and a fractal branch $F$ rooted at $w$. Suppose $B$ is of type (a). Then the blinking sequence of $w$ must be generated by diagram (A) in Figure \ref{branchorbit_root_fig} and refine a sequence generated by (F10/9) or (F11/8) in Figure \ref{fractalbranch_coarsening}. First, note that there is no way to refine 16 and 21 using the blinking gaps in diagram (A). Hence the blinking sequence of $w$ must refine the constant sequence of gap 19. It remains to show that sequence $19,19,19$ cannot be refined by any sequence generated by diagram (A). Note that there are 3 ways to refine 19 using weights in diagram (A): $(YX)(XW)$, $(WY)(YY)$, and $(YX)(XX)$(here we may take $Y=Z$). Notice that none of them uses gap 6 inside sequence \eqref{W} or begins with an edge emanating from node $X$ in diagram (A). Hence the first 19 must be refined by $(WY)(YY)$, and the second 19 must be refined by $(YX)(XW)$, but then the following 19 cannot be refined. This shows that the branch of type (a) is exclusive. 
	
	Now suppose the branch $B$ is of type (b). A blinking sequence generated by diagram (B2) in Figure \ref{branchorbit_root_fig} must refine a sequence generated by (F10/9) or (F11/8) in Figure \ref{fractalbranch_coarsening}. To this end, we claim the following: among all directed walks in diagram (B2), 
	\begin{description}[noitemsep]
		\item{(a)} $Z,(ZX),X$ is the only walk which generates a sequence ($9,7,6$) that refines 22;
		\item{(b)} $W,(WY),Y$ is the only walk which generates a sequence ($7,7,7$) that refines 21;
		\item{(c)} $(XY),Y$ is the only walk which generates a sequence ($10,7$) that refines 17;
		\item{(d)} $(XZ),Z$ and $(XW),W$ are the only walks which generate sequences ($10,9$ and $12,7$, respectively) refining 19.
	\end{description}
	To see this, for instance, consider possible ways to refine 22 using diagram (B2). If gap 12 is used, then it must be $22=10+12$, but 12 cannot be preceded or be followed by 10; if 11 is used, it must be $22=11+11$, but this is also impossible; if 10 is used, then 12 must be properly refined, but this is impossible; if 9 is used, then $13=6+7$ is the only way to refine 13, and $Z(ZX)X$ is the only way to generate 6,7, and 9 consecutively. This shows (a), and the other claims can be shown similarly.  
	
	Now we show that we show that any sequence generated by diagram (B2) refines no sequence generated by (F11/8). By (c) and (d), no refinement of 17 can be followed or preceded by any refinement of 19. Since in diagram (F11/8) 19 always follows 17, we see that 17 cannot be refined. This yields that the blinking sequence of $w$ may only refine the constant sequence of 19, but by (d) any refinement of 19 begins with node $X$ and ends with nodes $Z$ or $W$, so 19 cannot be refined repeatedly. 
	
	It remains to show that any sequence generated by diagram (B2) refines no sequence generated by (F10/9). We have seen in the previous paragraph that the constant sequence 19 cannot be refined. So if every 19 is refined, then at some point a refinement of 22 or 16 should follow. But by (a) and (d), the refinement of 22 cannot follow any refinement of 19. This makes that the directed walk in diagram (F10/9) which generates a sequence refined by some sequence generated from diagram (B2) cannot use the loop at node F1, and consequently, also the right-left edge of weight 16; this implies that only a constant sequence 19 from (F10/9) can be refined, which contradicts our earlier observation. This shows the assertion. 
\end{proof}

\begin{proof}[\textbf{Proof of Proposition \ref{kstar_orbit}.}] By the minimality, we may assume that the number of distinct phases occupied by the leaves in $S$ is at least $2$ and constant in time. At each time $t$, by a component we mean the set of consecutive states on the leaves on the hexagon $\mathbb{Z}_{6}$, where the nodes $\{0,1,\dots, 5\}$ are ordered in a clockwise manner; the size of a component is the number of distinct phases in it. Notice that by Lemma \ref{lem:opposite_leaf}, whenever $v$ blinks, every component must lie entirely clockwise or counterclockwise without any leaf opposite to $v$ (color difference of 3). Hence the number of components is a non-increasing function in time, which must be constant in the time since we are in a periodic orbit. Let us call any local configuration in such a closed orbit stable.

	Figure \ref{starorbit_fig} (a) is the only closed local orbit with a single component of size $3$; Figure \ref{starorbit_fig} (d) shows the closed local orbits with a single component of size 2; Figure \ref{starorbit_fig} (c) shows the closed orbits with 2 components of size 1 and 2; Figure \ref{starorbit_fig} (b) for 2 components of size 1 for both. Notice that any configuration of a component of size $\ge 4$ is unstable, likewise, anyone with two components with one component of size $\ge 3$, and anyone with two components with both has size $\ge 2$. So the nine configurations in Figure \ref{starorbit_fig} give all stable local configurations. By the time-invariance types, transitions between local configurations in different types are impossible. Possible transitions within each type and their minimal transition times are investigated similarly in Figure \ref{leaf_transition_ex}. For instance, Figure \ref{star_orbit_fig1} illustrates possible transitions from Figure \ref{starorbit_fig} $b_{1}$ to $b_{2}$. This shows \textbf{\eqref{I}}.
	\begin{figure}[H]
		\centering
		\includegraphics[width=0.4 \linewidth]{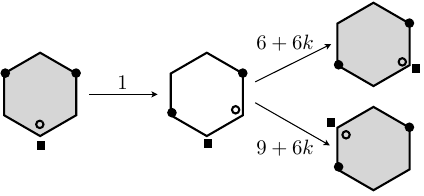}
		\caption{ Transitions from Figure \ref{starorbit_fig} $b_{1}$ to $b_{2}$. Weights on edges indicate a number of iterations, where $k$ is a non-negative integer depending on the dynamics. Due to the rotational symmetry of the targeting configuration, there are two possible phases for $v$ to move into. }
		\label{star_orbit_fig1}
	\end{figure}

	To see \textbf{(ii)}, suppose $S=T_{v}$ is a branch. Note that $v^{-}$ is the only external neighbor of $v$, so $v$ can get at most one external pull from $v^{-}$ in every 6 iterations. This makes the five unshaded local configurations in Figure \ref{starorbit_fig} impossible to appear on branches. For instance, consider possible transitions from Figure \ref{starorbit_fig} $b_{4}$, which is shown in Figure \ref{star_orbit_fig2}. During the second transition of length 5 from the second to the third column, either $v^{-}$ pulls $v$ as in the dotted bottom transition or not as in the solid upper transition. When $v$ blinks for the next time, none of the resulting configurations in the last column is stable. This shows the bottom left local configuration in Figure \ref{starorbit_fig} is impossible on branches. A similar argument applies to the other four unshaded local configurations in Figure \ref{starorbit_fig}. Thus there are exactly three possible types of closed orbit for this branch as asserted in \textbf{(ii)}.
	
	\begin{figure}[H]
		\centering
		\includegraphics[width=0.6 \linewidth]{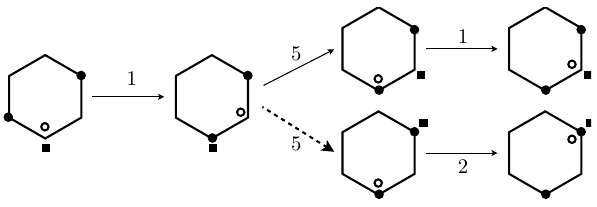}
		\caption{ Assuming $S$ is a branch, local configuration in Figure \ref{starorbit_fig} $b_{4}$ leads to unstable local configurations. 
		}
		\label{star_orbit_fig2} 
	\end{figure}

	Now we show \textbf{(iii)}. First, suppose the local dynamic on $T_{v}$ is given by Figure \ref{starorbit_fig} (a). In such a local orbit, in terms of standard representation, the leaves must have colors $0,1$ and $2$ whenever $v$ blinks. The following sequence shows the first 8 iterations starting from such local configuration (Figure \ref{starorbit_fig} $a_{1}$):
	\begin{equation}	\label{branch_typeA_1}
		\begin{matrix}
			\text{leaves}    &01\textbf{2} & 1\textbf{2}3 & \mathbf{2}34 & 345 & 450 & 501 & 01\mathbf{2} & 1\mathbf{2}3 \\
			\text{$v$}	&\textbf{2} & 3  & 3 & 3 &  a_{1}   &  a_{2}    &   a_{3}  &  a_{4}   \\
			\text{$v^{-}$} &b_{1}&b_{2}&b_{3}&b_{4}&b_{5} & b_{6}
		\end{matrix}
	\end{equation}
	Clearly $a_{3}\ne 2$, and it is easy to check that $a_{3}\ne 5$ leads to a different local configuration at the next blink of $v$: hence we must have $a_{3}=5$. This requires $2\in \{b_{4},b_{5},b_{6}\}$, which in particular yields that $T_{v}$ is open. But $b_{4}=2$ leads to a contradiction since it would yield $b_{1}=5$ and $b_{2}=0$; so $2\in \{b_{5},b_{6}\}$. We extend sequence (\ref{branch_typeA_1}) as follows: 
	\begin{equation}	\label{branch_typeA_2}
		\begin{matrix}
			\text{leaves}    &01\textbf{2} & 1\textbf{2}3 & \mathbf{2}34 & 345 & 450 & 501 & 01\mathbf{2} & 1\mathbf{2}3 & \mathbf{2}34 & 345 & 450 & 501 & 01\mathbf{2} \\
			\text{$v$}	&\textbf{2} & 3  & 3 & 3 &  a_{1}   &  a_{2}    &   5  &  5 & 5 & 5 & x_{1} & x_{2} & x_{3}   \\
			\text{$v^{-}$} &&&&b_{4} & b_{5} &b_{6}&b_{7}&b_{8}&b_{9}&b_{10}&
		\end{matrix}
	\end{equation}
	Note that $2\in \{b_{5},b_{6}\}$ yields $b_{9}\ne 2$, so $x_{1}=0$, $X_{2}=1$, and $X_{3}=2$. This shows a single transition from Figure \ref{starorbit_fig} $a_{1}$ to itself takes exactly 12 iterations, and since $2\in \{b_{5},b_{6}\}$, $T_{v}$ is of type (a) definition.
	
	Next, suppose the local dynamic on $T_{v}$ is given by Figure \ref{starorbit_fig} (b). The argument is similar for type (a). We will show that the transition $b_{2}\rightarrow b_{1}$ and $b_{1}\rightarrow b_{2}$ in Figure \ref{starorbit_fig} takes 9 and 7 iterations, respectively. We look at the first 9 iterations starting from Figure \ref{starorbit_fig} $b_{2}$: 
	
	\begin{equation*}
		\begin{matrix}
			\text{leaves} & 41 & 4\mathbf{2} & 53 & 04 & 15 & \mathbf{2}0 & 31 & 4\mathbf{2}  \\
			\text{$v$} & \textbf{2} & 3 & 3 &x_{1} & x_{2} & x_{3} & x_{4} &  x_{5} \\
			\text{$v^{-}$} & a_{1}&a_{2}& a_{3} &a_{4} &a_{5}&a_{6}&a_{7}&a_{8}
		\end{matrix}
	\end{equation*}
	
	We need to have $2\notin \{a_{3},a_{4},a_{5}\}$ since otherwise $x_{5}=2$ and the resulting local configuration is not Figure \ref{starorbit_fig} $b_{1}$. This makes $x_{3}=5$ and we may extend the sequence further:
	\begin{equation*}
		\small{\begin{matrix}
				\text{leaves} & 41 & 4\mathbf{2} & 53 & 04 & 15 & \mathbf{2}0 & 31 & 4\mathbf{2} & 53 & 04 & 14 & \mathbf{2}5 & 30 & 41 & 5\mathbf{2} & 03 & 14 \\
				\text{$v$} & \textbf{2} & 3 & 3 &* & * & 5 & 5 &  0 &  1 &  \textbf{2} &  3 &  b_{1} &  b_{2} &  b_{3} &b_{4} &b_{5} &b_{6} \\
				\text{$v^{-}$} & a_{1}&a_{2}& a_{3} &a_{4} &a_{5}&a_{6}&a_{7}&a_{8}&a_{9}&a_{10}&a_{11}&a_{12}&a_{13}&a_{14}&a_{15}&a_{16}&a_{17}
		\end{matrix}}
	\end{equation*}
	
	This in particular shows that the transition $b_{2}\rightarrow b_{1}$ in Figure \ref{starorbit_fig} takes 9 iterations. Furthermore, $a_{4}\ne 2$ since it leads to a contradiction by back-tracking in time, so we have  $2\in \{a_{3},a_{5}\}$. Hence by definition, $T_{v}$ would be of type (b) if the transition $b_{1}\rightarrow b_{1}$ in Figure \ref{starorbit_fig} takes 7 iterations, i.e., $b_{6}=2$. To this end, it is enough to show that $2\notin \{a_{11},a_{13},a_{14}\}$. Indeed, $a_{11}\ne 2$ since otherwise $b_{6}=0$ so the local configuration `$v|\text{leaves}$' after two more iterations from the end of the above sequence would be $2|30$, which is not what we should have as in Figure \ref{starorbit_fig} $b_{2}$. Similarly, $2\in \{a_{13},a_{14}\}$ leads to a wrong local configuration $2|30$, so $2\notin \{a_{11},a_{13},a_{14}\}$. Thus $T_{v}$ is of type (b). 
	
	Finally, suppose the local dynamic on $T_{v}$ is given by Figure \ref{starorbit_fig} (c). First five iterations from Figure \ref{starorbit_fig} $c_{1}$ is a follows:
	\begin{equation*}
		\begin{matrix}
			\text{leaves}    &034 & 134 & \mathbf{2}45 & 350 & 401 & 50\mathbf{2} \\
			\text{ $v$ }	&\textbf{2} & 3  &  x_{1}   &  x_{2}    &   x_{3}  &   5\\
			\text{ $v^{-}$ }	&w_{1} & w_{2} &w_{3} &w_{4} &w_{5} &w_{6} 
		\end{matrix}
	\end{equation*}
	In order for this local dynamics lead to Figure \ref{starorbit_fig} $c_{2}$, we need to have $2\in \{w_{2},w_{4},w_{5}\}$. However, $w_{4}=2 2$ would lead to a contradiction by back-tracking upto $w_{1}$, so $2\in \{w_{2},w_{5}\}$. An entirely similar argument for previous cases shows that the transitions $c_{1}\rightarrow c_{2}$ and $c_{2}\rightarrow c_{1}$ in Figure \ref{starorbit_fig} take exactly 10 and 9 iterations, respectively. Thus $T_{v}$ is a fractal of type 10/9.  This shows the assertion. 
\end{proof}

\section{Proof of Lemma \ref{lem:terminal_is_fractal}}
\label{sec:proof_lem_terminal_fractal}

By Proposition \ref{starorbit_fig} \textbf{(iii)} we know that type (c) terminal branches are fractal, so in order to show Lemma \ref{lem:terminal_is_fractal}, it suffices to show that no terminal branches can be of type (a) or (b). We do this in the following two subsections.

\subsection{Terminal branches are not of type (a)}

We begin by ruling out type (a) terminal branches.  

\begin{prop}\label{type_a_prop1}
	Let $(T,X_{0})$ be as before. If there are two types (a) terminal branches $B$ and $B'$ rooted at the same vertex $w$, then one of the two branches must only use the sequence \eqref{P}, and the other must only use \eqref{Q}, which are given in the proof of Proposition \ref{branchorbit_detail}. 
\end{prop}

\begin{proof}
	First note that if the blinking sequence of $w$ ever uses the term 11, then because there is only one weight of 11 in Figure \ref{branchorbit_root_fig} (a), both branches undergo the sequence \eqref{S} in synchrony. Since $w$ fluctuates the centers of $B$ and $B'$ in the same way, the two branches will be in synchrony thereafter, contradicting the minimality. Thus we may assume that $w$ never has a blinking gap 11. Similarly, we may assume that blinking gap 13 never appears for $w$. In general, the same argument applies to any unique sequence generated by diagram (A) in Figure \ref{branchorbit_root_fig}, such as 7-7, 12-6, and 12-7. Once we exclude such segments, the only possible directed closed walk in Figure \ref{branchorbit_root_fig} (A) is the ones that use loops on nodes \eqref{P}, \eqref{Q} or \eqref{R}. Since the induced dynamics on $w$ must coincide, this is possible only if one of the two branches constantly uses sequence \eqref{P} and the other \eqref{R}, as asserted.  
\end{proof}

\begin{prop}\label{type_a_prop2}
	Let $(T,X_{0})$ be as before. Then there is no terminal branch of type (a).  
\end{prop}

\begin{proof}
	Suppose there is a terminal branch $T_{v}$ of type (a). Then each component in $T_{v^{-}}-v^{-}$ is either a leaf or a branch. First, suppose there is another branch, say $T_{u}$, rooted at $v^{-}$. By Lemma \ref{exclusive_branches}, it must be of type (a) as well. By Proposition \ref{type_a_prop1}, we may assume that $T_{v}$ only uses sequence $\eqref{P}$ and $T_{u}$ only $\eqref{R}$.  Any more branches rooted at $w$ will be redundant. So we assume these two are the only branches rooted at $w$. Consider the following sequence, which is obtained by overlapping \eqref{P} and \eqref{R} by matching dynamics on $v^{-}=u^{-}$:
	\begin{equation}
		\centering{
			\begin{matrix}
				u	&\mathbf{2}& 3 & 3 & 3 & 4 & 5 & 5 & 5 & 5 & 5 & 0 & 1 & \mathbf{2} & 3\\
				v	&1&\textbf{2}	& 3 & 3 & 3 & 4 & 4 & 5 & 5 & 5 & 5 & 0 & 1 & \textbf{2} \\
				v^{-}	&-&5 & 5 & 0 & 1 & \textbf{2} & 3 & - & - & - & - & - & 5 & 5
			\end{matrix}\tag{PR}
		}
	\end{equation}
	Note that since $*\ne 2$, this sequence requires $v^{-}$ to be pulled four times in 6 iterations when it goes through the $-$'s. Since any vertex blinks at most once in 6 iterations, this means that $v^{-}$ must have at least 4 external neighbors except $v$ and $u$. Thus except its own parent $v^{--}$, it must have at least three leaves. By Proposition \ref{kstar_orbit}, the local dynamics on this 3-star centered at $v^{-}$  the local dynamic should be given by Figure \ref{starorbit_fig} (a) or (c). However, the latter is not possible since in our circumstance the blinking sequence of $v^{-}$ is the constant sequence $12,12,\cdots$, which is not the form of $10+6k_{1},9+6k_{2},10+6k_{3},\cdots$. Thus the 3-star centered at $v^{-}$ must go through type (a) closed orbit. In particular, whenever $v^{-}$ blinks, its three leaves must have colors $0,1,2$. Adding this to $(P+Q)$, we see that the local dynamics on $T_{v^{-}}$ must be of the concatenation of the following sequence 
	\begin{equation}
		\centering{
			\begin{matrix}
				u	& 3 & 3 & 3 & 4 & 5 & 5 & 5 & 5 & 5 & 0 & 1 & \mathbf{2} & 3\\
				v	&\textbf{2}	& 3 & 3 & 3 & 4 & 4 & 5 & 5 & 5 & 5 & 0 & 1 & \textbf{2} \\
				v^{-}	&5 & 5 & 0 & 1 & \textbf{2} & 3 & 3 & - & - & - & 5 & 5 & 5 
				\\
				\text{leaves} &  \mathbf{2}34 & 345 & 450 & 501 & 01\mathbf{2} & 1\mathbf{2}3 &\mathbf{2}34&345&450&501&01\mathbf{2}&1\mathbf{2}3&\mathbf{2}34
			\end{matrix}
		}\tag{PRl}
	\end{equation}
	
	There are multiple contradictions at this point: $v^{-}$ still needs to be pulled twice from external neighbors when it goes through $-$'s in the above sequence $(PRl)$ but $v^{--}$ is the only remaining external neighbor; whenever $u$ or $v$ pulls $v^{-}$, some leaf pulls $v^{-}$ together, so the branches $T_{v}$ and $T_{u}$ are not contributing anything to the dynamics on $v^{-}$, contradicting the minimality.   
	
	Hence we may assume that there is no other branch rooted at $v^{-}$. Observe that since $v^{-}$ has color 4 at the end of sequence $\eqref{S}$, it must be concatenated with the other three, so sequences $\eqref{P}$-$\eqref{R}$ must be used at least once in the periodic local dynamics. Note that $v^{-}$ must have color 4 or 5 at the end of sequences $\eqref{P}$-$\eqref{R}$ in order to be concatenated by the following one. Since $v^{-}$ does not blink and $v$ does not pull $v^{-}$ within those three sequences, it means that $v^{-}$ must be pulled either by its own parent $v^{--}$ or by its leaves, at least four times during the last six iterations in the three sequences. Since every vertex blinks at most once in every six iterations, this yields that $v^{-}$ needs to have at least three leaves except its own parent. On this 3-star centered at $v^{-}$, Proposition \ref{kstar_orbit} again enforces the local dynamics given by Figure \ref{starorbit_fig} (a) (cf. (c) is not the case as before). Thus every blinking gap of $v^{-}$ should be of the form $12+6k$. Since sequence \eqref{S} contains a blinking gap 6, it cannot be used and only the other three can be. Now the same sequence $(PRl)$ shows that whenever $T_{v}$ goes through sequences \eqref{P} or \eqref{R}, $v^{-}$ is pulled by some leaf whenever it is pulled by $v$. Since $v$ blinks at the exact same time in sequences \eqref{Q} and \eqref{R}, this means $T_{v}$ is redundant to the dynamics of $v^{-}$, contradicting the minimality. This shows the assertion. 
\end{proof}

\subsection{Terminal branches are not of type (b)} 

Next, we rule out type (b) terminal branches. 

\begin{prop}\label{type_b_prop1}
	Let $(T,X_{0})$ as before. Suppose $T$ has a terminal branch $T_{v}$ of type (b). Then the local dynamic on $T_{v}$ does not use the long periodic sequences $\eqref{I}$ and $\eqref{J}$. In particular, the exact blinking sequence of $v^{-}$ is generated by digraph (B2) in Figure \ref{branchorbit_root_fig}.
\end{prop}	

\begin{proof}
	Suppose on the contrary that sequences \eqref{I} and \eqref{J} do appear. This means $v^{-}$ has a blinking gap from one of the four edge weights in Figure \ref{branchorbit_root_fig} (B1). We are going to show that $T_{v}$ is the only branch rooted at $v^{-}$ and $v^{-}$ has at most two leaves. The assertion then easily follows. Indeed, in the last 6 iterations in both sequences \textbf{\eqref{I}} and \eqref{J}, $v^{-}$ must be pulled at least four times. Since $v$ does not pull $v^{-}$ during this period, its parent $v^{--}$ and two other leaves cannot provide this.  
	
	We first show that $v^{-}$ has at most two leaves. Suppose not. By Proposition \ref{kstar_orbit}, the 3-star centered at $v^{-}$ has local dynamics given by Figure \ref{starorbit_fig} (a) or (c). Suppose the former. Then the blinking sequence of $v^{-}$ is of the form $12+6k_{1}, 12+6k_{2},\cdots$. Among the weights in Figure \ref{branchorbit_root_fig} (B1), the edge weight 18 of (IJ) is the only one of that form, and the following blinking gap of $v^{-}$ should be either 16 of (JJ), 14 of (JI), or their refinements. The first two are not of the prescribed form, so they must be their refinements. In Figure \ref{branchorbit_root_fig} (B2), the edge (IJ) in the coarsened diagram (B1) corresponds to the node $Y$ and edge (YW) combined. Thus any refined blinking gap of $v^{-}$ must use the node weight 7 at $W$ in the diagram (B2), which also conflicts with the prescribed form. Hence $v^{-}$ and its leaves cannot have local dynamics given by Figure \ref{starorbit_fig} (a). 
	
	Assuming local dynamics on $T_{v}$ given by Figure \ref{starorbit_fig} (c), the blinking sequence of $v^{-}$ must be of the form $10+6k_{1},9+6k_{2},10+6k_{3},\cdots$. Notice that there is no weight of $9+6k$ for $k\ge 1$ in diagram (B1) and (B2) in Figure \ref{branchorbit_root_fig}, so the blinking sequence of $v^{-}$ must be of the form $10+6k_{1},9,10+6k_{2},\cdots$. The only weights of the from $10+6k$ in (B1) and (B2) in Figure \ref{branchorbit_root_fig} are 10 and 16. This yields that the blinking sequence $v^{-}$ must consist of three terms 9,10, and 16, where 10 and 16 are followed by 9, and 9 must be followed by 10 or 16. We shall see that this is impossible. Note that the sequence 10-9 is uniquely generated by the walk $(XZ),Z$ in Figure \ref{branchorbit_root_fig} (B2), but no edge emanating from node $Z$ in that digraph has weight 10 or 16. Thus 10 is not a blinking gap of $v^{-}$, so the blinking sequence must alternate 9 and 16. But such a sequence cannot refine any sequence generated by Figure \ref{branchorbit_root_fig} (B1), a contradiction. This shows that $v^{-}$ has at most two leaves. 
	
	It remains to show that there is no other branch rooted at $v^{-}$. Suppose on the contrary that another branch $T_{u}$ is rooted at $v^{-}$. By Lemma \ref{exclusive_branches}, $T_{u}$ is of type (b). By the minimality, branches $T_{v}$ and $T_{u}$ must have distinct dynamics.  Now if $v^{-}$ has a blinking gap of 14 or 18, then since those gaps are uniquely generated by Figure \ref{branchorbit_root_fig} (B1), the two branches must be synchronized thereafter, a contradiction. Thus $v^{-}$ never have blinking gaps 14 or 18, but does use gap 16, which are given by the loops at node $I$ or $J$ in Figure \ref{branchorbit_root_fig} (B1). We may assume that when $v^{-}$ has blinking gap 16, $T_{v}$ and $T_{u}$ undergo loops \textbf{(ii)} and (JJ) in Figure \ref{branchorbit_root_fig} (B1), respectively. Note that the loop \textbf{(ii)} is represented by the node $X$ and its loop $(XX)$ combined in the refining digraph Figure \ref{branchorbit_root_fig} (B2), so the blinking gap of $v^{-}$ that follows 16 should be coming from the four edges emanating from node $X$ in the same digraph, which only gives 10 or 12. By the parallel reasoning, loop (JJ) must be followed by an edge emanating from nodes $Z$ or $W$ in Figure \ref{branchorbit_root_fig}, which yields the next blinking gap should be either 7 or 9, a contradiction. Hence $T_{v}$ is the unique branch rooted at $v^{-}$.  This shows the assertion. 
\end{proof}

\begin{prop}\label{type_b_Prop2}
	Let $(T,X_{0})$ as before. Suppose $T$ has a terminal branch $T_{v}$ of type (b). Then $T_{v}$ is the only branch rooted at $v^{-}$. 
\end{prop}	

\begin{proof}
	Suppose for the contrary there is another branch $T_{u}$ rooted at $v^{-}$. By Lemma \ref{exclusive_branches} we know that $T_{u}$ must be of type (b), and by Proposition \ref{type_b_prop1}, they never use long periodic sequences \textbf{\eqref{I}} and \eqref{J} so that the blinking sequence of $v^{-}$ is generated by digraph (B2) in Figure \ref{branchorbit_root_fig}. By minimality, these two branches must not be synchronized. This means that we must be able to find two distinct closed walks in digraph (B2) which generate the same sequence. Since the weights 6, 11, and 12 are unique in the diagram, any such blinking sequence cannot use those numbers. Thus we may delete the node $X$ together with all the indecent edges, and also the edge $(YW)$ of weight 11 from the digraph. The resulting digraph, which generates the blinking sequence of $v^{-}$ in our current situation, is provided below:
	\begin{figure}[H]
		\centering
		\includegraphics[width=0.4 \linewidth]{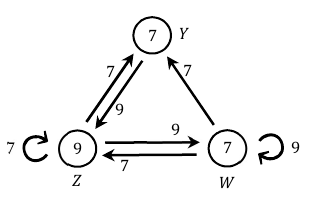}
		\vspace{-0.4cm}
		\caption{ If $T_{v}$ is a terminal branch of type (b), then the blinking sequence of $v^{-}$ is generated by this digraph.
		}
		
		\label{type_b_proof_fig}
	\end{figure} 
	
	Note that by Proposition \ref{branchorbit_detail} \textbf{(iv)}, both branches must use \eqref{W} at least once. We claim that the blinking sequence of $v^{-}$ never repeats $9$ twice. This would yield the assertion as follows. Under this assumption, it would be impossible to use the edge $(ZW)$; thus no edge heading toward \eqref{W} is available, so after a branch uses the node $(W)$, then it must be confined there. Thus both branches use node $W$ only (recall that we are in a periodic orbit), and since they should generate the same blinking sequence for $v^{-}$, they must have synchronized dynamics, a contradiction. 
	
	Thus it suffices to show that the blinking sequence of $w$ cannot repeat 9 twice. Suppose for contrary that $\mathtt{g}_{1}(v^{-})=\mathtt{g}_{2}(v^{-})=9$. Observe that there are only two ways to generate 9-9 from Figure \ref{branchorbit_root_fig} (B2) with node $\eqref{X}$ deleted: $(YZ),Z$ and $Z,(ZW)$. Thus we may assume $T_{v}$ goes through $(YZ)$ and $T_{u}$ goes though $(ZW)$ simultaneously. Since the string 7-7-7 is uniquely generated by $W,(WY),Y$ in the above digraph, it never appears in the blinking sequence of $v^{-}$. This forces $T_{u}$ to be confined at node $W$ after the third blink, forcing $(\mathtt{g}_{i})_{i\ge 3}$ to alternate 7 and 9. This contradicts the periodicity of the blinking sequence, so string 9-9 never appears in the blinking sequence of $v^{-}$. This shows the assertion. 
\end{proof}	

\begin{prop}\label{type_b_Prop3}
	Let $(T,X_{0})$ as before. Then $T$ has no terminal branch of type (b). 
\end{prop}	

\begin{proof}
	Suppose on the contrary that $T_{v}$ is a terminal branch of type (b). By Proposition \ref{type_b_Prop2}, we know that there is no other branch rooted at $v^{-}$. Thus all neighbors of $v^{-}$ except $v$ and its own parant $v^{--}$ are leaves. Furthermore, by Proposition \ref{type_b_prop1}, the local dynamics on $T_{v}$ uses only those short periodic sequences \eqref{X}-\eqref{W} and the blinking sequence of $v^{-}$ is generated by digraph (B2) in Figure \ref{branchorbit_root_fig}.

	We first show that sequence \eqref{X} is never used. To see this, notice that in sequence \eqref{X}, $v^{-}$ is to be pulled at least four times during the last six iterations. Since $v$ does not pull $v^{-}$ during this period, $v^{-}$ must have at least three leaves. By Proposition \ref{kstar_orbit}, $v^{-}$ has exactly three leaves and the 3-star centered at $v^{-}$ has local dynamics given by Figure \ref{starorbit_fig} (a) or (c). This implies that the blinking gap of $v^{-}$ must always be the form of $12+6k_{1}$, $10+6k_{2}$, $9+6k_{3}$ for $k_{1},k_{2},k_{3}\ge 0$. But sequence \eqref{X} forces $v^{-}$ to have blinking gap 6, which is not of the prescribed form, a contradiction. Thus sequence \eqref{X} is never used by the local dynamics on $T_{v}$. In particular, $\mathtt{g}_{i}(v^{-})\in \{7,9,11\}$ for all $i\ge 1$. 
	
	Next, we show that sequences \eqref{Y} and \eqref{Z}  also are not used by the local dynamics on $T_{v}$. Suppose not. Proposition \ref{branchorbit_detail} \textbf{(iv)}, we know that sequence \eqref{W}  is used (periodically). In particular, $v^{-}$ has a blinking gap 7 periodically. Observe that sequence \eqref{Y} itself requires $v^{-}$ to be pulled four times in the last six iterations, and sequence $\eqref{Z}\eqref{W}$ concatenated also requires the same. Thus $v^{-}$ needs to have at least two leaves. Combining Proposition \ref{kstar_orbit} together with the conclusion of the previous paragraph and the fact that $v^{-}$ does have blinking gap 7, we see that the local dynamics on the star centered at $v^{-}$ must only use local configurations Figure \ref{starorbit_fig} $b_{1}$ or $b_{2}$, and blinking sequence of $v^{-}$ alternates 7 and 9. But then once $T_{v}$ uses the node W in digraph (B2) in Figure \ref{branchorbit_root_fig}, it must confine on node $W$, contradicting the periodicity of local dynamics. Thus $T_{v}$ only uses sequence \eqref{W} , and the blinking sequence of $v^{-}$ alternates 7 and 9. 
	
	Note that the concatenated sequence \eqref{W}\eqref{W} requires $v^{-}$ to be pulled at least three times at the end of the first \eqref{W}. Since $v$ does not pull $v^{-}$ during this period, $v^{-}$ needs to have at least one leaf. If it has at least two leaves, then combining sequence \eqref{W} with the local dynamics on $v^{-}$ together with its leaves given by Figure \ref{starorbit_fig} (b), the local dynamics on $T_{v^{-}}$ is given by repeating the following sequence 
	\begin{equation}
		\begin{matrix}
			\text{$v$} & \textbf{2} & 3 & - & - & - & - & 5 &  0 &  1 &  \textbf{2} &  3 &  - &  - &  5 & 0 & 1  & \mathbf{2}\\
			\text{$v^{-}$} &5 & 5 & 0&1& \textbf{2} &3 &- & - &5 & 0&1&\textbf{2}&3&-&-&-&-\\
			\text{leaves} & \mathbf{2}0&31&4\mathbb{2}&53& 04 & 14 & \mathbf{2}5 & 30 & 41 & 5\mathbf{2} & 03 & 14 & \mathbf{2}4 & 35 & 40 & 51 & 0\mathbf{2}
		\end{matrix}\tag{Wl}
	\end{equation}
	Note that in the above sequence, whenever $v$ blinks, one of the two leaves of $v^{-}$ blinks as well. Hence the branch $T_{v}$ is redundant to the dynamics of $v^{-}$, which contradicts minimality. So we may assume $v^{-}$ has exactly one leaf. The local dynamics on this 1-star centered at $v^{-}$ is given by the digraph in Figure \ref{localdynamicswithleaf}. The only compatible closed walk there which generates a sequence that alternates 7 and 9 is $e,(ea),a,(ad),d,(dd),d,(de),e$ (upto choice of starting node). Notice that during the loop $(ee)$, the leaf of $v^{-}$ has color $4$ when $w$ is blinking. If we plug the leaf in \eqref{W}, we get the following sequence: 
	\begin{equation}
		\centering{
			\begin{matrix}
				\text{$v$} & \textbf{2} & 3 & - & - & - & - & 5 &  0 &  1 &  \textbf{2} &  3 &  - &  - &  5 & 0 & 1  & \mathbf{2}\\
				\text{$v^{-}$} &5 & 5 & 0&1& \textbf{2} &3 &- & - &5 & 0&1&\textbf{2}&3&-&-&-&- \\
				\text{leaf} &  0&1&\mathbf{2}&3& 4 & 4 & 5 & 0 & 1 & \mathbf{2} & 3 & 4 & 4 & 5 & 0 & 1 & \mathbf{2}
			\end{matrix}\tag{Wl1}
		}
	\end{equation}
	which should appear in the local dynamics on $T_{v^{-}}$ periodically. Now the last four iterations are conflicting since the leaf of $v^{-}$ does not contribute to extra pull on $v^{-}$. This shows the assertion.
\end{proof}

\section{Proof of Lemma \ref{lem:recursion_fractalbranch}}
\label{sec:proof_lemma_recursion_fractal}

In this section, we show Lemma \ref{lem:recursion_fractalbranch}. Let $(T,X_{0})$ as before and let $T_{v^{-}}$ as stated in Lemma \ref{lem:recursion_fractalbranch}. We say a neighbor of $v^{-}$ \textit{external} if it is either a leaf or its own parent $v^{--}$. By Proposition \ref{kstar_orbit}, $v^{-}$ has at most three leaves. Hence $v^{-}$ can have at most four external pulls during every six iterations. Since large blinking gaps of $v^{-}$ generally require lots of external pulls, it would be not likely under our hypothesis. In fact, blinking gaps of $v^{-}$ can be at most 11, as stated in the following proposition:  

\begin{prop}\label{fractalbranch_recursion_prop1}
	Let $(T,X_{0})$ be as before. Suppose that each connected component of $T_{v^{-}}-v^{-}$ is either a singleton or fractal. Further, assume that at least one such component is fractal. Then the blinking gaps of $v^{-}$ are bounded above by 11. 
\end{prop}

Our strategy for showing the above statement is the following: we collect all possible subsequences arising from the two sequences \eqref{F1}, \eqref{F2} and their eight concatenations (F$i$)(F$j$) for $(i,j)\in \{1,2\}^{2}\cup \{2,3\}^{2}$ with respect to the induced blinking gap of $v^{-}$, and count the number of required external pulls. Detailed proof of this statement is given at the end of this section. 

For further discussions, we give a full list of possible subsequences of (F$i$)(F$j$) generating a fixed blinking gap of $v^{-}$. generating blinking gaps $\le 11$ for $v^{-}$ below: 
\begin{description}
	\item{(11)} Blinking gap 11:
	\begin{equation*}
		\begin{matrix}
			\text{\eqref{F1}($a_{5}$)} & 3&3&-&-&-&-&5&0&1&\textbf{2}&3&-\\ 
			\text{\eqref{F2}($b_{8}$)} & -&-&-&5&0&1&\textbf{2}&3&-&-&-&-\\
			\text{\eqref{F3}($c_{6}$)} & -&-&-&-&-&-&5&0&1&\textbf{2}&3&-\\
			\text{\eqref{F4}($d_{8}$)} & -&-&-&-&5&0&1&\textbf{2}&3&-&-&-\\
			\text{($a_{2}$)\eqref{F1}} & 1& \textbf{2}&3&-&-&-&-&5&0&1&\textbf{2}&3\\
			\text{($a_{5}$ or $b_{5}$)\eqref{F2}} & -&-&-&-&5&0&1&\textbf{2}&3&-&-&-\\
			\text{($c_{3}$ or $d_{3}$)\eqref{F3}} & 1&\textbf{2}&3&-&-&-&5&0&1&\textbf{2}&3&-\\
			\text{($c_{5}$ or $d_{5}$)\eqref{F4}} & 3&-&-&-&5&0&1&\textbf{2}&3&-&-&-\\
			v^{-}	&  \textbf{2}&3&x_{1}&x_2&x_3&x_4&x_5&x_6& 5 &0& 1&\textbf{2}
		\end{matrix}
	\end{equation*}
	
	\item{(10)} Blinking gap 10:
	\begin{equation*}
		\begin{matrix}
			\text{\eqref{F1}($a_{4}$)}& 3&3&-&-&-&-&5&0&1&\textbf{2}&3\\ 
			\text{\eqref{F2}($b_{7}$)}& -&-&-&5&0&1&\textbf{2}&3&-&-&-\\
			\text{\eqref{F3}($c_{5}$)}& -&-&-&-&-&-&5&0&1&\textbf{2}&3\\
			\text{($a_{3}$ or $b_{3}$)\eqref{F1}}&  \textbf{2}&3&-&-&-&-&5&0&1&\textbf{2}&3\\
			\text{($c_{4}$ or $d_{4}$)\eqref{F3}}& \textbf{2}&3&-&-&-&5&0&1&\textbf{2}&3&*\\
			\text{($c_{6}$ or $d_{6}$)\eqref{F4}}& -&-&-&5&0&1&\textbf{2}&3&-&-&-\\
			v^{-}	&  \textbf{2}&3&x_1&x_2&x_3&x_4&x_{5}&5& 0& 1&\textbf{2}
		\end{matrix}
	\end{equation*}
	
	\item{(9)} Blinking gap 9:
	\begin{equation*}
		\begin{matrix}
			\text{\eqref{F1}($a_{3}$)}&  3&3&-&-&-&-&5&0&1&\textbf{2}\\ 
			\text{\eqref{F3}($c_{4}$)}&  -&-&-&-&-&-&5&0&1&\textbf{2}\\
			\text{\eqref{F4}($d_{6}$)}&  -&-&-&-&5&0&1&\textbf{2}&3&-\\
			\text{($a_{4}$ or $b_{4}$)\eqref{F1}}&   3&-&-&-&-&5&0&1&\textbf{2}&3\\
			\text{($a_{7}$ or $b_{7}$)\eqref{F2}}& -&-&5&0&1&\textbf{2}&3&-&-&-\\
			\text{($c_{5}$ or $d_{5}$)\eqref{F3}}& 3&-&-&-&5&0&1&\textbf{2}&3&-\\
			v^{-}	&  \textbf{2}&3&x_1&x_2&x_3&x_4&5& 0& 1&\textbf{2}
		\end{matrix}
	\end{equation*}
\end{description}
\begin{description}
	\item{(8)} Blinking gap 8:
	\begin{equation*}
		\begin{matrix}
			\text{\eqref{F1}($a_{2}$)}& 3&3&-&-&-&-&5&0&1\\ 
			\text{\eqref{F2}($b_{5}$)}& -&-&-&5&0&1&\textbf{2}&3&-\\
			\text{\eqref{F3}($c_{3}$)}& -&-&-&-&-&-&5&0&1&\\
			\text{\eqref{F4}($d_{5}$)}& -&-&-&-&5&0&1&\textbf{2}&3\\
			\text{($a_{5}$ or $b_{5}$)\eqref{F1}}&  -&-&-&-&5&0&1&\textbf{2}&3\\
			\text{($a_{8}$ or $b_{8}$)\eqref{F2}}&  -&5&0&1&\textbf{2}&3&-&-&-\\
			\text{($c_{6}$ or $d_{6}$)\eqref{F3}}& -&-&-&5&0&1&\textbf{2}&3&-\\
			\text{($c_{8}$ or $d_{8}$)\eqref{F4}}& -&5&0&1&\textbf{2}&3&-\\
			v^{-}	&  \textbf{2}&3&x_1&x_2&x_3&5& 0& 1&\textbf{2}
		\end{matrix}
	\end{equation*}
	
	\item{(7 and 6)} Blinking gaps 7 and 6:
	\begin{equation*}
		\begin{matrix}
			\text{\eqref{F1}($a_{1}$)}& 3&3&-&-&-&-&5&0 &\qquad &\text{\eqref{F2}($b_{3}$)}& -&-&-&5&0&1&\textbf{2}\\
			\text{\eqref{F2}($b_{4}$)}& -&-&-&5&0&1&\textbf{2}&3 &\qquad &\text{\eqref{F4}($d_{3}$)}& -&-&-&-&5&0&1\\
			\text{\eqref{F3}($c_{2}$)}& -&-&-&-&-&-&5&0 &\qquad &\text{($a_{1}$)($a_{7}$)}& 0&1&\textbf{2}&3&-&-&-\\
			\text{\eqref{F4}($d_{4}$)}& -&-&-&-&5&0&1&\textbf{2} &\qquad &\text{($a_{2}$)($a_{8}$)}& 1&\textbf{2}&3&-&-&-&-\\
			\text{($a_{1}$)($a_{8}$)}& 0&1&\textbf{2}&3&-&-&-&-&\qquad &\text{($c_{2}$)($c_{8}$)}& 0&1&\mathbf{2}&3&-&-&- \\
			v^{-}	&  \textbf{2}&3&x_1&x_2&5& 0& 1&\textbf{2} &\qquad&v^{-}	 &  \textbf{2}&3&4&5& 0& 1&\textbf{2}
		\end{matrix}
	\end{equation*}	
\end{description}

For instance, \eqref{F1}($a_{5}$) is the subsequence of \eqref{F1} from the first blink of $v^{-}$ to the second blink $a_{5}=2$; ($a_{2}$)\eqref{F1} is the subsequence of \eqref{F1}\eqref{F1} from $a_{2}=2$ to the first blink of $v^{-}$ in the second \eqref{F1}. Note that the last three sequences for gap 6 are contradictory, so only the first two are possible.

\begin{prop}\label{fractalbranch_recursion_prop2}
	Let $(T,X_{0})$ be as before. Suppose that each connected component of $T_{v^{-}}-v^{-}$ is either a singleton or fractal. Further, assume that at least one such component is fractal. Then $v^{-}$ has at least two leaves. 
\end{prop}

\begin{proof}
	Suppose for the contrary that $v^{-}$ has at most one leaf. By Proposition \ref{fractalbranch_recursion_prop1}, blinking gaps of $v^{-}$ are at most $\le 11$. Note that blinking gap 10 is impossible since it requires at least three external pulls during the first six iterations while $v^{-}$ has at most two external neighbors.

	We first claim that $v^{-}$ needs to have at least two leaves in order to have a blinking gap 9. Since gap 9 requires at least one leaf for $v^{-}$, we may assume to the contrary that $v^{-}$ has exactly one leaf. Note that ($a_{7}$ or $b_{7}$)\eqref{F2} is necessary to generate gap 9 with $v^{-}$  since otherwise blinking gap 9 would require at least three external neighbors. Among the sequences which generate gaps $\in\{6,7,8,11\}$, \eqref{F2}($b_{5}$) for gap 8 is the only sequence that can follow. Similarly, the sequence ($b_{5}$)\eqref{F1} for gap 8 can only follow this. Hence we only need to rule out consecutive gaps 8-8. Note that in Figure \ref{localdynamicswithleaf}, $a\rightarrow c\rightarrow d$ and $e\rightarrow b \rightarrow d$ are the only possible walks that generate blinking sequence 8-8 for $v^{-}$. But note that during the second transition in each walk, the center is not pulled by the leaf. This makes the second blinking gap 8 for $v^{-}$ during ($b_{5}$)\eqref{F1} impossible. This shows the second claim. Thus we may assume that $\mathtt{g}_{i}(v^{-})\in \{6,7,8,11\}$ for all $i\ge 1$. 
	
	Our second claim is that the assertion holds assuming $v^{-}$ has blinking gaps $\le 8$ only. By ruling out sequences above which cannot be concatenated by any other to the right or left, we find that the local dynamics of $v$ should be given by repeating the following sequences:
	\begin{equation*}
		(a_{8})\eqref{F2}-\eqref{F2}(b_{5})-(b_{5})\eqref{F1}-\eqref{F1}(a_{1})-(a_{1})(a_{8}) ;\qquad 8-8-8-7-7
	\end{equation*}
	which induce stings of blinking gaps of $v^{-}$ as indicated on the right. By the asymmetry of such strings and minimality, this yields that $T_{v}$ is the unique fractal subtree rooted at $v^{-}$. However, this means that whenever $v^{-}$ has blinking gap 8 induced by sequence \eqref{F2}($b_{5}$), in which $v$ does not pull $v^{-}$, $v^{-}$ needs two external pulls, a contradiction. This shows the second claim.

	Now we show the assertion. If $v^{-}$ has no leaf, then blinking gap $11$ is impossible so the assertion follows from the two claims. Hence we may assume that $v^{-}$ has one leaf. By Proposition \ref{kstar_orbit}, the blinking sequence of $v^{-}$ is generated by the digraph in Figure \ref{starorbit_fig}. In fact, only the edges or loops with weights in $\{6,7,8,11\}$ can be used. Moreover, edges of weight 11 or 8 emanating nodes $b$ or $c$ in the digraph cannot be used, since in the first six iterations starting from those local configurations $v^{-}$ is not pulled by the leaf. Deleting all those edges, we obtain the following digraph which should generate the blinking sequence of $v^{-}$:
	\begin{figure}[H]
		\centering
		\includegraphics[width=0.35 \linewidth]{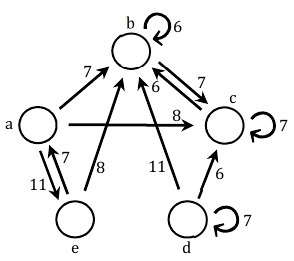}
		\caption{ A digraph generating the blinking sequence of $v^{-}$  when it has a single leaf under the hypothesis of Proposition \ref{fractalbranch_recursion_prop2}.
		}
		\label{type_b_proof_fig}
	\end{figure} 
	
	By the second claim, $v^{-}$ must have blinking gap 11, and the only closed walk in the above digraph containing an edge of weight 11 is the one alternating between node $a$ and $e$. Hence the blinking sequence of $v^{-}$ should alternate 11 and 7, but no sequences for gap 7 can be followed by any sequence for 11. This shows the assertion. 
\end{proof}

\begin{prop}\label{fractalbranch_recursion_prop3}
	Let $(T,X_{0})$ be as before. Suppose that each connected component of $T_{v^{-}}-v^{-}$ is either a singleton or fractal. Further, assume that at least one such component is fractal. Then the blinking sequence of $v^{-}$ alternates 10 and 9 or 11 and 8. 
\end{prop}

\begin{proof}
	By Proposition \ref{fractalbranch_recursion_prop1}, we know that the blinking gaps of $v^{-}$ are at most $11$. We first show that the blinking sequence of $v^{-}$ alternates 8 and 11 or 9 and 10. By Proposition \ref{fractalbranch_recursion_prop2}, $v^{-}$ has at least two leaves. By Proposition \ref{kstar_orbit}, $v^{-}$ never have blinking gap 6. If $v^{-}$ has three leaves, then by Proposition \ref{kstar_orbit} the local dynamics on the 3-star centered at $v^{-}$ should be given by Figure \ref{starorbit_fig} (c), so $v^{-}$ has blinking gap alternating 10 and 9. Hence we may assume that $v^{-}$ has exactly two leaves. Next, we show that string 9-9 also cannot appear in the blinking sequence of $v^{-}$. The two leaves of $v^{-}$ force that the following blinking gap of $v^{-}$ after 9-9 should be either 10 or 11. From the list of subsequences generating blinking gap $9$, we see that $v^{-}$ must blink at $a_{3}$ or $c_{4}$ at the end of the second blinking gap 9. But no subsequence for gap 11 begins with these, and the fourth and fifth subsequence for gap 10 do begin with these but during which $v^{-}$ is not pulled by $v$: this requires $v^{-}$ to have four external pulls during a blinking gap 10, a contradiction.

	Next, we rule out blinking gap 7 for $v^{-}$. Suppose $v{-}$ does have a blinking gap 7. Then the local dynamics on the 2-star centered at $v^{-}$ is given by digraph (b) in Figure \ref{starorbit_fig}. Moreover, since gap $\ge 12$ and 9-9 does not appear, the only possible local configurations for this 2-star are Figure \ref{starorbit_fig} $b_{3}$ and $b_{4}$. This forces that gap 7 is always followed and preceded by gap 9. Observe that for sting 9-7-9 in the blinking sequence of $v^{-}$, by considering possible concatenations of the subsequences in Figure \ref{lem:terminal_is_fractal}, one sees that the second 9 after 7 should be given by $(b_{4})(\text{F1})$. This uniqueness forces that $T_{v}$ is the only fractal subtree rooted at $v^{-}$. Moreover, the second 9 cannot be followed by 7, since otherwise the second 7 is given by $\eqref{F1}(a_{2})$, but no sequence for 9 begins with $(a_{2})$. Thus the blinking sequence of $v^{-}$ must contain the string 9-7-9-10-9-7-9. However, the second 7 in this string must end with $(a_{1})$, but no sequence for gap 9 begins with $(a_{1})$. Thus $v^{-}$ does not have blinking gap 7 if it has at least two leaves.

	Now we may assume that $v^{-}$ has exactly two leaves only the gaps 8,9,10, and 11 appear. If the 2-star centered at $v^{-}$ has local dynamics confined in digraphs (b) or (c) in Figure \ref{starorbit_fig}, then we are done. Hence we may assume that the local dynamics are given by digraph (d) in Figure \ref{starorbit_fig}. We want to show that the local configuration on the 2-star alternates local configurations $d_{2}$ and $d_{3}$ in Figure \ref{starorbit_fig} (d). Since we have shown that 9-9 does not appear in the blinking sequence of $v^{-}$, it suffices to rule out the strings 8-9-11 and 11-11. First, observe that the former is uniquely generated by \eqref{F4}($d_{5}$)-($d_{5}$)\eqref{F3}-\eqref{F3}($c_{6}$). Hence $v^{-}$ has at most one fractal branch rooted at it and needs at least four external pulls during the last blinking gap of 11, a contradiction. To rule out the string 11-11, observe that there are five sequences that generate consecutive blinking gap 11 of $v^{-}$:
	\begin{eqnarray*}
		&&\text{\eqref{F1}($a_{5}$)-($a_{5}$)\eqref{F2}}\\
		&&\text{($a_{2}$)\eqref{F1}-\eqref{F1}($a_{5}$)}\\
		&&\text{($a_{5}$ or $b_{5}$)\eqref{F2}-\eqref{F2}($b_{8}$)}\\
		&&\text{($c_{3}$ or $d_{3}$)\eqref{F3}-\eqref{F3}($c_{6}$)}\\
		&&\text{($c_{5}$ or $d_{5}$)\eqref{F4}-\eqref{F4}($d_{4}$)}
	\end{eqnarray*}
	Since the blinking sequence of $v^{-}$ is generated by the digraph (d) in Figure \ref{kstar_orbit}, the next blinking gap after 11-11 should be either 9 or 11. Note that no subsequences generating those blinking gaps could be concatenated after the last three sequences above. This yields that $v^{-}$ could have at most two fractal subtrees rooted at itself whose local dynamics during the second blinking gap 11 should be given by subsequences ($a_{5}$)\eqref{F2} or \eqref{F1}($a_{5}$). But then the second blinking gap 11 of $v^{-}$ requires at least four external pulls, which is impossible with only two leaves for $v^{-}$. This shows the assertion.
\end{proof}

Now we are ready to prove Lemma \ref{lem:recursion_fractalbranch}.

\begin{proof}[\textbf{Proof of Lemma \ref{lem:recursion_fractalbranch}.}] By Propositions \ref{fractalbranch_recursion_prop1}, \ref{fractalbranch_recursion_prop2}, and \ref{fractalbranch_recursion_prop3}, we may assume that $v^{-}$ has at least two leaves and its blinking sequence alternates 10 and 9 or 11 and 8. To show that $T_{v^{-}}$ is fractal, we need to show that $v^{--}\in V(T)$ and it provides external pulls on $v^{-}$ at right place. First, let us analyze the blinking sequence that alternates 8 and 11. In this case, by Proposition \ref{kstar_orbit} we may assume that $v^{-}$ has exactly two leaves. By Proposition \ref{kstar_orbit}, the 2-star centered at $v^{-}$ should alternate between the two local configurations in Figure \ref{starorbit_fig} $d_{2}$ and $d_{3}$. Hence when $v^{-}$ blinks at the beginning of a gap 11, its two leaves must have colors 3 and 4. Consider the following 11 iterations during a gap 11:
	\begin{equation*}
		\begin{matrix}
			v^{-}			& \mathbf{2} &3 &x_{1}&x_2&x_3&x_4&x_5&x_6&5&0&1&\mathbf{2}\\
			\text{leaves} 	& 34 & 34 & 45 & 50 & 01 & 1\mathbf{2} & \mathbf{2}3 & 45 & 50  & 01 & 1\mathbf{2} & \mathbf{2}3\\
			v^{--}			& y_{1} & y_{2} & y_{3} & y_{4} & y_{5} & y_{6} & y_{7} & y_{8}
		\end{matrix}
	\end{equation*}
	From the list of possible sequences for gap 11, $v^{-}$ is not pulled by any of its internal neighbors (i.e., centers of fractal subtrees rooted at $v^{-}$) during the transition $x_{1}\rightarrow x_{5}$. In order to make gap 11, $v^{-}$ needs to be pulled by an external neighbor during $x^{1}\rightarrow x^{3}$. Thus $v^{--}\in V(T)$ and $2\in \{y_{2},y_{3},y_{4}\}$. Since $y_{4}=2$ yields $y_{1}=5$ which is a contradiction, we have $2\in \{y_{2},y_{4}\}$. This shows $T_{v^{-}}$ is fractal of type 11/8, as desired. 
	
	Now we assume that the blinking sequence of $v^{-}$ alternates 10 and 9. The argument is similar to before. By Proposition \ref{kstar_orbit}, the $k$-star ($k\in \{2,3\}$) centered at $v^{-}$ should alternate local configurations $b_{1}$ and $b_{2}$ or $c_{1}$ and $c_{2}$ in Figure \ref{starorbit_fig}. In the first case, $v^{-}$ has two leaves which have colors 0 and 4 when $v^{-}$ blinks at the beginning of gap 10 as in Figure \ref{kstar_orbit} $b_{2}$; in the second case it has three leaves of color 0,3, and 4 at a such instant as in Figure \ref{starorbit_fig} $c_{1}$. The following sequence shows ten iterations during a blinking gap 10 together with all possible local dynamics on the leaves of $v^{-}$:
	\begin{equation*}
		\begin{matrix}
			v^{-}			& \mathbf{2} &3 &x_{1}&x_2&x_3&x_4&x_5&5&0&1&\mathbf{2}\\
			\text{leaves} 	& 04 & 14 & \mathbf{2}5 & 30 & 41 & 0\mathbf{2} & 13 & \mathbf{2}5 & 30  & 41 & 5\mathbf{2} \\
			\text{leaves} 	& 034 & 134 & \mathbf{2}45 & 350 & 401 & 01\mathbf{2} & 1\mathbf{2}3 & \mathbf{2}35 & 340  & 451 & 50\mathbf{2} \\
			v^{--}			& y_{1} & y_{2} & y_{3} & y_{4} & y_{5} & y_{6} & y_{7} & y_{8}
		\end{matrix} 
	\end{equation*}
	By the list of sequences giving blinking gap 10 for $v^{-}$, there are no internal pulls on $v^{-}$ during the first six iterations in the above sequence. Hence $v^{-}$ still needs one extra external pulls, and this yields $v^{--}\in V(T)$ with $2\in \{y_{2},y_{4},y_{5}\}$. Since $y_{4}\ne 2$ for similar reason this shows that $T_{v^{-}}$ is fractal of type 10/9. This shows the assertion. 
\end{proof}

\begin{proof}[\textbf{Proof of Proposition \ref{fractalbranch_recursion_prop1}.}] By Proposition \ref{prop_fb_1}, the maximum possible blinking gap of $v^{-}$ is 22 generated by \eqref{F1}\eqref{F2} without secondary blink within \eqref{F1}, but this requires at least 5 external pulls during six iterations, so it cannot occur. Blinking gap 20 arises from \eqref{F3}\eqref{F4} but is impossible for a similar reason, and there is no subsequence that gives blinking gap 20 (e.g., see Figure \ref{fractalbranch_coarsening}). We rule out large blinking gaps from 19 to 12 below.

	\begin{description}
		\item{(19)}  There are only four subsequences giving gap 19, namely, (F$i$)(F$i$) for $1\le i\le 4$ without $v^{-}$ blinking more than once in the first sequence (F$i$). Hence if $v^{-}$ has blinking gap 19, then there can be at most four fractal subtrees rooted at $v^{-}$. We overlap all four sequences to see the least number of required external pulls:
		\setcounter{MaxMatrixCols}{25}
		\begin{equation*}
			\small{\begin{matrix}
					\text{\eqref{F1}\eqref{F1}}&3&3&-&-&-&-&5&0&1&\textbf{2}&3&-&-&-&-&5&0&1&\textbf{2}&3\\
					\text{\eqref{F2}\eqref{F2}}&-&-&-&5&0&1&\textbf{2}&3&-&-&-&-&5&0&1&\textbf{2}&3&-&-&-\\
					\text{\eqref{F3}\eqref{F3}}&-&-&-&-&-&-&5&0&1&\textbf{2}&3&-&-&-&5&0&1&\textbf{2}&3&-\\
					\text{\eqref{F4}\eqref{F4}}&-&-&-&-&5&0&1&\textbf{2}&3&-&-&-&5&0&1&\textbf{2}&3&-&-&-\\
					v^{-} & \textbf{2} &3& x_{1} & x_{2} & x_{3} & x_{4} & x_{5} &  &  &  &  &  &  &  &  &  & 5&0&1& \textbf{2}
			\end{matrix}}
		\end{equation*}
		where row (F$i$)(F$i$) above indicates the dynamics of the center of a fractal subtree of sequence (F$i$)(F$i$). Since $v^{-}$ does not blink during gap 19 and we have all possible fluctuations from the fractal subtrees in the above matrix, all missing pulls on $v^{-}$ must be external. Now note that $3\le x_{5} \le 5$, so $v^{-}$ must get at least three external pulls during the first six iterations above. This requires $v^{-}$ to have at least two leaves, but blinking gap 19 is impossible to be generated from digraphs in Figure \ref{starorbit_fig}, contradicting Proposition \ref{kstar_orbit}. This rules out blinking gap 19.  
		\item{(18)} Impossible
		\item{(17)} Gap 17 arises uniquely from \eqref{F4}\eqref{F3}, so in this case $v^{-}$ can have at most one fractal subtree. In the following sequence during blinking gap 17, $v^{-}$ needs to get at least 5 external pulls during the first six iterations, which is impossible.  
		\begin{equation*}
			\begin{matrix}
				\text{\eqref{F4}\eqref{F3}}&-&-&-&-&5&0&1&\textbf{2}&3&-&-&-&5&0&1&\textbf{2}&3&-\\
				v^{-}&\textbf{2}&3& x_{1} & x_{2} & x_{3} & x_{4} & x_{5} &x_{6}&&&&&&&5& 0& 1&\textbf{2}
			\end{matrix}
		\end{equation*}
		
		\item{(16)} Gap 16 arises uniquely from \eqref{F2}\eqref{F1}, so $v^{-}$ can have at most one fractal subtree. Consider the following sequence during blinking gap 16:	
		\begin{equation*}
			\begin{matrix}
				\text{\eqref{F2}\eqref{F1}} &*&*&*&5&0&1&\textbf{2}&3&*&*&*&*&5&0&1&\textbf{2}&3\\
				v^{-}  & \textbf{2}&3&x_{1}&x_{2}&x_{3}&x_{4}&x_{5}&&&&&&&5& 0& 1&\textbf{2} 
			\end{matrix}
		\end{equation*}
		If $x_{5}=5$, then $v^{-}$ needs five external pulls in a row after $x_{5}$, which is impossible. Thus $3\le x_{5}\le 4$. Hence $v^{-}$ needs at least four external pulls for the first six iterations in the above sequence. So it has exactly three leaves. By Proposition \ref{kstar_orbit}, the 3-star centered at $v^{-}$ must have local configuration Figure \ref{starorbit_fig} $c_{1}$ in order to match with blinking gap 16. Inserting the three leaves according to such local configuration, the first six iterations in the above sequence look as follows:
		\begin{equation*}
			\begin{matrix}
				\text{\eqref{F2}\eqref{F1}} &-&-&-&5&0&1&\textbf{2}\\
				v^{-}  & \textbf{2}&3& 3& 3&x_{3}&x_{4}&x_{5}\\
				\text{leaves} & 01\textbf{2} & 1\textbf{2}3 & \mathbf{2}34 & 345 & 450 & 501 & 01\mathbf{2} 
			\end{matrix}
		\end{equation*}
		But since $x_{5}\le 4$, this requires $v^{-}$ to be pulled by its only remaining external neighbor, namely its parent, twice for the last three iterations, a contradiction. 
		\item{(15)} Impossible.
		\item{(14)} Within the sequence \eqref{F1}, we could have blinking gap 14 if $a_{7}=2$. Also possible is $a_{2}=2$ and \eqref{F1} is followed by \eqref{F2}. Last possibility is that $c_{1}=2$ in \eqref{F3} and \eqref{F4} follows. Denote these three cases by \eqref{F1}($a_{7}$), ($a_{2}$)\eqref{F2}, and ($c_{3}$)\eqref{F4}. This gives the following sequence for blinking gap 14: 
		\begin{equation*}
			\begin{matrix}
				\text{\eqref{F1}($a_{7}$)} & 	3&3&-&-&-&-&5&0&1&\textbf{2}&3&-&-&-&-\\
				\text{($a_{2}$)\eqref{F2}}  &  1&\textbf{2}&3&-&-&-&-&5&0&1&\textbf{2}&3&-&-&-\\
				\text{($c_{3}$)\eqref{F4}} & 0&1&\mathbf{2} &-&-&-&-&5&0&1& \mathbf{2}& 3 &-&-&-\\
				v^{-}	&  \textbf{2}&3&x_{1}&x_{2}&x_{3}&x_{4}&x_{5}&x_{6}&x_{7}&x_{8}&x_{9}&5& 0& 1&\textbf{2}
			\end{matrix}
		\end{equation*}
		Hence $v^{-}$ requires at least four external pulls during the five iterations from $x_{2}$ to $x_{7}$. Hence $v^{-}$ has exactly three leaves. But neither digraphs (a) nor (c) in Figure \ref{starorbit_fig} can generate gap 14, a contradiction. 
		
		\item{(13)} There are four possibilities for gap 13 as below, which clearly requires $v^{-}$ to have at least two leaves. But blinking gap 13 is not generated from any digraphs in Figure \ref{starorbit_fig}, contrary to Proposition \ref{kstar_orbit}. 
		\begin{equation*}
			\begin{matrix}
				\text{\eqref{F1}($a_{6}$)}& 3&3&-&-&-&-&5&0&1&\textbf{2}&3&-&-&-\\
				\text{\eqref{F3}($c_{8}$)}& -&-&-&-&-&-&5&0&1&\textbf{2}&3&-&-&-\\
				\text{($c_{3}$)\eqref{F4}} & 1&\textbf{2}&3&-&-&-&5&0&1&\textbf{2}&3&-&-&-\\
				\text{($d_{3}$)\eqref{F4}} & 1 & \textbf{2} &3&-&-&-&5&0&1& \mathbf{2}& 3 & -&-&- \\
				v^{-}	&  \textbf{2}&3&x_{1}&x_{2}&x_{3}&x_{4}&x_{5}&x_{6}&x_{7}&x_{8}&5& 0& 1&\textbf{2}
			\end{matrix}
		\end{equation*}
		
		\item{(12)} $v^{-}$ could have at most six fractal subtrees generating blinking gap 12. First, observe that $v^{-}$ cannot have three leaves. To see this, note that by Proposition \ref{kstar_orbit} its blinking gaps should be of the form  $12+6k_{i}$ for $k_{i}\ge 0$; since we have seen that $v^{-}$ does not have a blinking gap 18, its blinking gap should be 12 constantly. But a fractal subtree rooted at $v^{-}$ makes this impossible (e.g., no subsequence is possible between the six possibilities below). Second, suppose $v^{-}$ has two leaves. Then by Proposition \ref{kstar_orbit}, the 2-star centered at $v^{-}$ must have local configuration Figure \ref{starorbit_fig} $b_{1}$ or $b_{3}$. Their dynamics are inserted in the following matrix:
		\begin{equation*}
			\begin{matrix}
				\text{($a_{1}$)\eqref{F1}} & 0&1&\textbf{2}&3&-&-&-&-&5&0&1&\textbf{2}&3\\
				\text{($a_{4}$)\eqref{F2}} & 3&-&-&-&-&5&0&1&\textbf{2}&3&-&-&-&\\
				\text{($b_{4}$)\eqref{F2}} & 3&-&-&-&-&5&0&1&\textbf{2}&3&-&-&- \\
				\text{($c_{2}$)\eqref{F3}} & 0&1&\mathbf{2} &-&-&-&-&5&0&1&\mathbf{2}& 3 & -\\
				\text{($c_{4}$)\eqref{F4}} & \textbf{2}&3&-&-&-&5&0&1&\textbf{2}&3&-&-&-\\
				\text{($d_{4}$)\eqref{F4}} & \textbf{2}&3&-&-&-&5&0&1&\textbf{2}&3&-&-&-\\
				v^{-}	&  \textbf{2}&3&x_{1}&x_{2}&x_{3}&x_{4}&x_{5}&x_{6}&x_{7}&5& 0& 1&\textbf{2} \\
				\text{leaves ($b_{1}$)} & 04 & 14 & \mathbf{2}5 & 30 & 41 & 5\mathbf{2} & 03 & 14 & \mathbf{2}5 \\
				\text{leaves ($b_{3}$)} & 0\mathbf{2} & 13 & \mathbf{2}4 & 35 & 40 & 51 & 0\mathbf{2} & 13 & \mathbf{2}4 
			\end{matrix}
		\end{equation*}
		Note that we cannot have both of the last two rows at the same time. Considering each case separately, we see that $v^{-}$ still needs at least two external pulls, a contradiction. The above matrix also shows that $v^{-}$ needs at least two leaves, so blinking gap 12 is impossible. 
	\end{description} 
	This shows the assertion. 
\end{proof}

\section*{Acknowledgement}

This work is supported by NSF DMS-2010035 and partially by DMS-2206296. The author gives special thanks to David Sivakoff and Steven Strogatz for valuable discussions.

\section*{Data Availability}

The data sets that we generated in the present study are available in the repository \url{https://github.com/HanbaekLyu/FCA}.

\bibliographystyle{amsalpha}   
\bibliography{mybib}  

\appendix
\renewcommand*{\appendixname}{}
\numberwithin{equation}{section}
\numberwithin{figure}{section}
\numberwithin{table}{section}

\end{document}